\documentclass{article}
\usepackage{amsmath}
\usepackage{amssymb}
\usepackage{amsthm}
\usepackage{hyperref}
\usepackage[left=2.7cm,top=2.7cm,right=2.7cm,bottom=2.7cm]{geometry}
\usepackage{algorithm2e}
\usepackage{enumerate}
\usepackage[normalem]{ulem}
\usepackage[dvipsnames]{xcolor}
\usepackage{ dsfont, bbold }
\usepackage{graphicx}
\usepackage{tikz}
\usetikzlibrary{positioning,chains,fit,shapes,calc}
\usepackage{subcaption}
\usepackage{titlesec}

\titleformat*{\section}{\LARGE\bfseries}
\titleformat*{\subsection}{\Large\bfseries}
\titleformat*{\subsubsection}{\large\bfseries}
\titleformat*{\paragraph}{\large\bfseries}
\titleformat*{\subparagraph}{\large\bfseries}

\newtheorem{thm}{Theorem}[section]
\newtheorem{cor}[thm]{Corollary}
\newtheorem{lem}[thm]{Lemma}

\newtheorem{claim}[thm]{Claim}
\newtheorem{defn}[thm]{Definition}
\newtheorem{obs}[thm]{Observation}
\newtheorem*{lem:appendpreforb}{Lemma \ref{lem:appendpreforb}}
\newtheorem*{lem:union}{Lemma \ref{lem:union}}
\newtheorem*{lem:preforbpowbound}{Lemma \ref{lem:preforbpowbound}}
\newtheorem*{lem:cutcoi}{Lemma \ref{lem:cutcoi}}

\newcommand{\rbrac}[1]{\left(#1\right)}
\newcommand{\sbrac}[1]{\left[ #1\right]}
\newcommand{\cbrac}[1]{\left\{ #1\right\}}

\newcommand{\mc}[1]{\mathcal{#1}}
\newcommand{\tbf}[1]{\textbf{#1}}

\newcommand{\sm}{\setminus}

\newcommand{\pre}{\small\textsc{Pre}}
\newcommand{\coi}{\small{\textsc{Coi}}}
\newcommand{\edgcoi}{\small{\textsc{EdgCoi}}}
\newcommand{\pow}{\small{\textsc{Pow}}}
\newcommand{\maxpow}{\small{\textsc{MaxPow}}}

\newcommand{\free}{\small{\textsc{Free}}}

\newcommand{\prefbd}{\small{\textsc{PreFbd}}}
\renewcommand{\equiv}{\mathrm{equ}}
\def\Ext{\textsc{Ext}}
\def\ext{\mathrm{ext}}
\def\Ava{\textsc{Ava}}
\def\ava{\mathrm{ava}}
\newcommand{\nn}{\nonumber}

\def\Var{\mbox{{\bf Var}}}
\renewcommand{\P}{\mathbb{P}}
\def\E{\mathbb{E}}
\def\D{\Delta}
\def\d{\delta}
\def\eps{\epsilon}
\def\k{\kappa}

\def\a{\alpha}

\def\f{\phi}

\def\t{\tau}

\def\be{{\overline{e}}}
\def\bE{{\overline{E}}}
\def\bF{{\overline{F}}}
\def\bc{{\overline{c}}}
\def\bC{{\overline{C}}}
\def\bS{{\overline{S}}}
\def\bT{{\overline{T}}}
\def\bD{{\overline{D}}}

\def\bV{{\overline{V}}}

\def\bR{{\overline{R}}}

\newcommand{\col}[1]{(\small{\textsc{colored}}, #1)}
\newcommand{\unc}[1]{(\small{\textsc{available}}, #1)}
\newcommand{\wild}{(?)}
\newcommand{\colport}{\mbox{col}}

\def\f{\frac{p-2}{\binom{p}{2} - q + 1}}
\def\m{\frac{1}{\binom{p}{2} - q + 1}}
\def\mf{\frac{p-3}{\binom{p}{2} - q + 1}}
\def\hm{\frac{1/2}{\binom{p}{2} - q + 1}}

\newcommand{\TT}{\Tilde{\Theta}}
\newcommand{\TO}{\Tilde{O}}
\newcommand{\TOM}{\Tilde{\Omega}}

\title{A random coloring process gives improved bounds for the Erd\H{o}s-Gy\'arf\'as problem on generalized Ramsey numbers }
\author{Patrick Bennett\thanks{Dept.~of Mathematics, Western Michigan University (\texttt{Patrick.Bennett@wmich.edu}). Supported in part by Simons Foundation Grant \#426894.} \and Andrzej Dudek\thanks{Dept.~of Mathematics, Western Michigan University (\texttt{Andrzej.Dudek@wmich.edu}). Supported in part by Simons Foundation Grant \#522400.} \and Sean English\thanks{Dept.~of Mathematics and Statistics, University of North Carolina Wilmington (\texttt{EnglishS@uncw.edu}).}}

\date{}

\begin{document}

\maketitle

\begin{abstract}
 The Erd\H{o}s-Gy\'arf\'as number $f(n, p, q)$ is the smallest number of colors needed to color the edges of the complete graph $K_n$ so that all of its $p$-clique spans at least $q$ colors. In this paper we improve the best known upper bound on $f(n, p, q)$ for many fixed values of $p, q$  and large $n$. Our proof uses a randomized coloring process, which we analyze using the so-called differential equation method to establish dynamic concentration. 
\end{abstract}

\section{Introduction}

Ramsey theory is a branch of combinatorics that seeks to understand the emergence of ordered substructures within an otherwise unordered large structure. One of the more well-known and difficult parameters in Ramsey Theory, the \textbf{$k$-color diagonal Ramsey number} $R_k(p)$, is the least integer $n$ such that any $k$-coloring of the edges of the complete graph $K_n$ contains a monochromatic clique on $p$ vertices. The current best general bounds on $R_k(p)$ are still relatively far apart, namely,
\[
2^{kp/4} \le R_k(p) \le k^{kp}.
\]
The lower bound above was proved by Lefmann \cite{L87} and the upper bound follows in a straightforward manner using a ``neighborhood chasing'' technique, perhaps first employed by Erd\H{o}s and Szekeres~\cite{ES35}.
(For better bounds see, e.g.,~\cite{C2009, CF2021, Sah2022, S2022, W2021}).

In 1975, Erd\H{o}s and Shelah \cite{E75} proposed a natural generalization of Ramsey numbers. In particular, given integers $n$, $p$ and $q$, with $N>p\geq 3$ and $q\in [\binom{p}2]$, a \textbf{$(p,q)$-coloring} of $K_n$ is an edge coloring that has the property that every $p$-clique of sees at least $q$ colors, and the \textbf{generalized Ramsey number}, $f(n,p,q)$, denote the least number of colors necessary to give a $(p,q)$-coloring of $K_n$. This generalization did not receive significant attention in the literature until Erd\H{o}s and Gy\'arf\'as began a more systematic study of these Ramsey numbers in 1997~\cite{EG97}. Due to their treatment, $f(n,p,q)$ is sometimes referred to as the \textbf{Erd\H{o}s-Gy\'arf\'as function}. It is worth noting that  
\[
R_k(p) = \min_n \{n: f(n,p,2) > k\}
\]
so the function $f(n,p,q)$ is in general at least as hard to determine as the multicolor diagonal Ramsey numbers.

\subsection{Bounds on \texorpdfstring{$f(n,p,q)$ and our main theorem}{}}

Erd\H{o}s and Gy\'arf\'as proved the following upper bound on $f(n,p,q)$ for fixed $p$ and $q$ in~\cite{EG97} using a short and straightforward application of the Lov\'asz Local Lemma :
\begin{equation}\label{eqn:LLLbound}
    f(n, p, q) = O\rbrac{n^{\frac{p-2}{\binom p2 - q + 1}}}.
\end{equation}
They also showed that this upper bound is sharp at the so-called linear and quadratic thresholds, i.e. if $q=\binom{p}{2}-p+3$, then $f(n,p,q)=\Theta(n)$. Similarly, if $q=\binom{p}2-\lfloor\frac{p}2\rfloor+2$, then $f(n,p,q)=\Theta(n^2)$. There are many pairs $p, q$ for which \eqref{eqn:LLLbound} has been improved, but no improvement has been made in the general case. We summarize such improvements here. Conlon, Fox, Lee and Sudakov~\cite{CFLS15} improved \eqref{eqn:LLLbound} for $q \le p-1$ by showing that $f(n, p, p-1)=n^{o(1)}$. Trivially, \eqref{eqn:LLLbound} can be improved to $O(n^2)$ whenever $q >\binom{p}2-\frac{p}2+2$. Other than the cases mentioned in the last two sentences, there are only a finite number of pairs $p, q$ for which the bound \eqref{eqn:LLLbound} has been improved (of which we will mention some and omit a few). Of particular interest is the diagonal case, when $p=q$. It is straightforward to see that a $(3,3)$-coloring is simply a proper edge coloring, and so $f(n,3,3)=n-1$ or $n$ depending on if $n$ is even or odd. For larger values, the bound in~\eqref{eqn:LLLbound}, along with a simple inductive argument~\cite{EG97} gives us that
\[
n^{\frac{1}{p-2}}-1\leq f(n,p,p)=O(n^{\frac{2}{p-1}})
\]
For $p=4$, Mubayi~\cite{M04} gave an explicit algebraic construction that shows $f(n,4,4)\leq n^{\frac{1}{2}+o(1)}$, matching the lower bound up to sub-polynomial terms. Cameron and Heath extended this idea to also show that $f(n,5,5)\leq n^{\frac{1}{3}+o(1)}$~\cite{CH18}, again matching the lower bound up to subpolynomial factors. Subsequentially, Cameron and Heath were also able to show $f(n,6,6)\leq n^{\frac{1}{3}+o(1)}$ and $f(n,8,8)\leq n^{\frac{1}{4}+o(1)}$~\cite{CH20}, which does not match the lower bound, but is better than \eqref{eqn:LLLbound}. Worth noting, the $(5,5)$, $(6,6)$ and $(8,8)$ colorings beating the lower bound also made use of a $(p,p-1)$-coloring  with $n^{o(1)}$ colors, given by Conlon, Fox, Lee and Sudakov~\cite{CFLS15}. 

We briefly discuss some additional lower bounds. Fox and Sudakov \cite{FS2008} proved that $f(n, 4, 3) = \Omega(\log n)$, which was an improvement on a bound given by Kostochka and Mubayi \cite{KM2008}. Subsequently, Conlon, Fox, Lee and Sudakov~\cite{CFLS15} also showed $f(n, p, p-1) =\Omega(\log n)$ for all $p \ge 5$. Pahoata and Sheffer \cite{PS19} adapted the idea of additive energy to obtain some new lower bounds on $f(n, p, q)$ when we are above the linear threshold. Fish, Pahoata and Sheffer \cite{FPS20} used similar methods to obtain several more lower bounds, as did Balogh, English, Heath and Krueger \cite{BEHK22}.

In this work, we improve upon \eqref{eqn:LLLbound} for a wide range of fixed $p$ and $q$, and large $n$. Our main theorem is the following.

\begin{thm}\label{thm:main}
For fixed $p, q$ with $q \le \frac{p^2 - 26p +55}{4}$, we have
\begin{equation*}
    f(n,p,q) =O\rbrac{ n^{\f} \log^{-\m} n}.
\end{equation*}
\end{thm}

As a corollary we obtain a new bound for an extremal problem first studied by Brown, Erd\H{os} and S\"os~\cite{BES1973}. Let $F^{(r)}(n; k, s)$ be the smallest number of edges $m$ such that every $r$-uniform hypergraph with $n$ vertices and $m$ edges contains some $k$ vertices inducing at least $s$ edges. In terms of classical extremal numbers, $F^{(r)}(n; k, s) = 1+ \mbox{ex}_{r}(n, \mc{G}_{k, s})$ where $\mc{G}_{k, s}$ is the family of all $r$-uniform hypergraphs on $k$ vertices and $s$ edges. Before we state our corollary we discuss the history of the problem. Their original paper \cite{BES1973} contained two unanswered questions that attracted particular interest. First, they asked whether $F^{(3)}(n; 6, 3)=o(n^2)$. This was answered in the positive by Ruzsa and Szemer\'{e}di \cite{RS1976}, and their proof featured the first ever appearance of the triangle removal lemma. The question was generalized to the conjecture that $F^{(3)}(n; s+3, s)=o(n^2)$ for all $s \ge 3$ (often called the Brown-Erd\H{o}s-S\'os conjecture). Very recently Conlon, Gishboliner, Levanzov and Shapira \cite{CGLS2023} gave the best-known result for this problem, which says
\[
F^{(3)}(n; s+O(\log s / \log \log s), s)=o(n^2).
\]
The second conjecture from \cite{BES1973} was that 
\[
\lim_{n \rightarrow \infty} n^{-2}F^{(3)}(n; k, k-2) 
\]
exists for all $k\ge 4$. This conjecture was proved for $k=4$ by Brown, Erd\H{o}s and S\"os \cite{BES1973no2}, for $k=5$ by Glock \cite{G2019}, and for $k=6$ very recently by Glock, Joos, Kim, K\"uhn, Lichev and Pikhurko \cite{GJKKLP2022} and in all of these results the limit was explicitly found. The full conjecture was proved very recently by Delcourt and Postle \cite{DP2022} (without finding the limit). A key ingredient in the proof is a constrained random process, which produces asymptotically optimal constructions.

Using a straightforward alterations approach, Brown, Erd\H{o}s and S\'os \cite{BES1973} proved that 
\begin{equation*}
    F^{(r)}(n; k, s) = \Omega\rbrac{n^{\frac{rs-k}{s-1}}},
\end{equation*}
which for many fixed $r, k, s$ is still the best known lower bound. When $r=2$, as a corollary to our Theorem \ref{thm:main} we obtain an improvement for many fixed $k, s$. This follows from the following observation. In a $\rbrac{k, \binom{k}{2}-s+2}$-coloring of $K_n$, no set of $k$ vertices can have $s$ edges of the same color. Thus $F^{(2)}(n; k, s)$ is at least the size of any color class in such a coloring, so 
\begin{equation*}
    F^{(2)}(n; k, s) \ge \frac{\binom n2 }{f(n, k, \binom{k}{2}-s+2)}.
\end{equation*}
In light of Theorem \ref{thm:main}, we have the following bound.
\begin{cor}
   For fixed $k, s$ with $s \ge \frac{k^2+24k-47}{4}$, we have
   \[
F^{(2)}(n; k, s) = \Omega\rbrac{n^{\frac{2s-k}{s-1}} \log^{\frac{1}{s-1}} n}.
   \]
\end{cor}

\subsection{The coloring procedure and differential equation method}

Now we will give some specifics behind the process which we will need throughout the proof. If $C'$ is a multiset of colors, say $C'=\{c_1 * r_1, \ldots c_k * r_k\}$ with $r_1, \ldots r_k \ge 1$ we say that $C'$ has $\sum_k(r_k-1)$ \tbf{repeated colors} or \tbf{repeats}. Note that an edge colored $p$-clique sees at least $q$ distinct colors if and only if the multiset of colors that appear on the clique has at most $\binom{p}2-q$ repeats.

\subsubsection*{The coloring procedure}

Let $V$ be a set of $n$ vertices. Our coloring procedure will have two phases, each of which uses randomness. Phase 1 will consist of a random coloring process that colors one edge at a time according to a rule we will describe shortly. For Phase 1 we will use a set $C$ of
\[
\k n^{\f} \log^{-\m} n
\]
colors for some large constant $\k$ depending only on $p$.
If Phase 1 is successful (which we will prove it is w.h.p.~\footnote{We say an event dependent on $n$ occurs \emph{with high probability} (abbreviated w.h.p.) if the probability of that event tends to one as $n$ tends to infinity.}), it colors almost all the edges. Phase 2 consists of a much simpler random coloring, where we use a fresh set $C'$ of colors, where $|C'|=|C|$ and each edge that was left uncolored by Phase 1 gets a uniform random color from $C'$ (independently for all edges).

We start our Phase 1 process at step $0$ with all edges uncolored, and at each step $i$ we will randomly color one random edge (details to follow). Let $U_i$ be the set of uncolored edges at the beginning of step $i$. At step $i$, we choose an edge $e_i \in U_i$ uniformly at random. We then assign $e_i$ a color $c_i$ chosen uniformly at random from all colors that are \emph{available} at $e_i$ (definition follows). First, we define
\begin{equation}\label{eqn:Rdef}
  R(s):= \left \lceil \frac{(s-2)\rbrac{\binom{p}{2}-q+1}}{p-2} \right \rceil -1.
\end{equation}
A color $c$ is called \tbf{available} at $e$ if assigning $c$ to $e$ does not cause any set of $s\le p$ vertices to have more than $R(s)$ repeats among its colored edges. Note that $R(p) = \binom p2 - q$ and so our definition of availability enforces that we have a partial $(p, q)$-coloring. If $c$ is not available at $e$ we say $c$ is \tbf{forbidden} at $e$. Of course, $c$ is forbidden at $e$ if and only if there exists some set $S \supseteq e$ with $|S|\le p$ such that the colored edges in $S$ have $R(|S|)$ repeats and $c$ appears among the colors. For such a set $S$ we say that \tbf{$c$ is forbidden at $e$ through $S$.}

We note that Phase 1 of our coloring process uses ideas similar to those used by Guo, Patton and Warnke \cite{GPW20}. In particular they use an  edge-coloring process which colors edges one at a time where each color is chosen randomly and uniformly from all ``available" colors (for some appropriate definition of ``available"). Also, the first two authors together with Cushman and Pralat \cite{BCDP2022} used a strategy somewhat similar to Phases 1 and 2 to bound $f(n, 4, 5)$.

For a formal description of Phase 1, see Algorithm~\ref{algorithm coloring}. Note that for technical reasons having to do with ensuring dynamic concentration of our random variables, we stop our coloring process at a predetermined step $i_{max}$ (defined later), even if the process could continue.

\begin{algorithm}[ht]
	\SetKwInOut{Input}{Input}
    \SetKwInOut{Output}{Output}
 \Input{Vertex set $V$, color set $C$}
 \Output{Partial $(p,q)$-coloring}
 \While{$i\leq i_{max}$ and every edge in $U_i$ has an available color}{
 Choose an edge $e_i$ uniformly at random from $U_i$\;
 Choose a color $c$ uniformly at random from all colors in $C$ available at $e_i$\;
 Assign edge $e_i$ the color $c$\; 
 Update $U_{i+1}=U_i\setminus\{e_i\}$\;
 Increment $i$\;
 }
 \caption{Partial $(p,q)$-Coloring Process}\label{algorithm coloring}
\end{algorithm}
The main work of this paper is to prove that w.h.p.\ Phase 1 is able to run until step $i_{max}$ and outputs a coloring of almost all the edges that has nice properties which will then allow us to very easily show that Phase 2 succeeds with positive probability.

\subsubsection*{Proof methods: The differential equation method, dynamic concentration and the local lemma}

Most of the work necessary to prove Theorem~\ref{thm:main} is in showing that w.h.p.\ Algorithm~\ref{algorithm coloring} does not terminate before step $i_{max}$ (i.e.\ until then, every uncolored edge has some available color), and that the partial $(p,q)$-coloring we obtain at the end of Phase 1 has certain nice properties that allow us to finish the coloring. 

To analyze Phase 1 we use the so-called \tbf{differential equation method} to establish \tbf{dynamic concentration} of a large family of random variables. We say a random variable is dynamically concentrated if there exists a deterministic function (which typically depends on $i$), which we call the \tbf{trajectory} for that random variable, such that w.h.p.\ the random variable stays within a narrow window of its trajectory as the process evolves. See \cite{BD20} for a gentle introduction to the differential equation method. This method evolved from work done at least as early as 1970 (see Kurtz \cite{Kurtz1970}). In the 1990's Wormald \cite{W1995, W1999} developed the method into a quite general tool. Indeed, he proved a ``black box" theorem, which guarantees dynamic concentration under some relatively simple conditions. Warnke \cite{Warnke2020} recently proved a stronger version of Wormald's ``black box" theorem. For our purposes it seems the existing theorems are insufficient, but we are nevertheless able to analyze our process using some standard arguments that resemble previous analyses of other processes. The tools (and style) we use for the differential equation method resemble the work of Bohman \cite{B2009}, Bohman and Keevash \cite{BK2010, BK2021}, and Bohman, Frieze and Lubetzky \cite{BFL10, BFL15}. The gentle introduction in \cite{BD20} uses similar tools and style.

In proving that w.h.p.\ Algorithm~\ref{algorithm coloring} gets to step $i_{max}$, we will obtain some nice properties that the partial edge coloring enjoys when the algorithm terminates. These properties will help us use the Lov\'asz Local Lemma~(see for example \cite{AS}) to show that the final uniform random coloring on the remaining uncolored edges is suitable to complete the $(p,q)$-coloring.

\subsection{Organization of the paper and notation}

In Section~\ref{sec:framework}, we will introduce a framework which allows us to describe a family of random variables relevant to the evolution of our edge coloring process. In Section~\ref{sec:traj}, we heuristically derive the trajectories of our random variables. In Section~\ref{sec:goodevent}, we give the bounds on our random variable which we will prove hold w.h.p.\ throughout the process. In Sections~\ref{sec:crude}-\ref{sec:G}, we prove that our trajectories indeed approximate our random variables. In Section~\ref{sec:platonic}, we prove a few elementary, but technical lemmas which are used in earlier sections, and in Section~\ref{sec:finishing}, we show that we can color the remaining $\binom{n}2-i_{max}$ edges that are not colored by Algorithm~\ref{algorithm coloring} with a small set of new colors, completing a $(p,q)$-coloring. In Section~\ref{sec:concluding} we discuss an extension of Theorem \ref{thm:main} by the first author, Delcourt, Li and Postle \cite{BDLP} which came out shortly after the first draft of this paper.

All asymptotics given throughout the paper will be with $n$ going to infinity and all parameters that do not depend on $n$ will be assumed to be constant unless otherwise specified. Given two functions $f=f(n)$ and $g=g(n)$, we will write $f=O(g)$ if there exists a constant $c$ such that $f\leq cg$, and $f=o(g)$ if $\frac{f}{g}\to 0$. We write $f=\Omega(g)$ if $g=O(f)$, and $f=\Theta(g)$ if $f=O(g)$ and $f=\Omega(g)$. We also will write $f=\TO(g)$ if there exists some constant $c$ such that $f\leq g\cdot \log^c n$, $f=\TOM(g)$ if $g=\TO(f)$, and $f=\TT(g)$ if $f=\TO(g)$ and $f=\TOM(g)$.

When working with a variable $X=X(i)$, where $i$ is a timestep from Algorithm~\ref{algorithm coloring}, we will write $\Delta X$ to mean $X(i+1)-X(i)$ (where the specific value of $i$ we are considering should be clear from context).

\section{Framework}\label{sec:framework}

In this section we introduce the framework we will use to analyze the coloring process. The process keeps going as long as each uncolored edge has some available color, and so we are naturally interested in following random variables.
\begin{defn}
Let $U_i \subseteq \binom V2$ be the set of uncolored edges at step $i$. For each $e \in U_i$, let $\Ava(e)\subseteq C$ be the set of colors available at $e$.
\end{defn}
Of course, the availability of a color $c$ at an edge $e$ depends on sets of vertices $S \supseteq e$, how many repeats each such set $S$ already has, and whether assigning $c$ to $e$ would be an additional repeat. To analyze this process we therefore need a model describing the evolving state of the coloring in terms of sets $S$ of $s \le p$ vertices and how the edges in $S$ may be colored.

From now on we call the elements of $V$ \tbf{real vertices}, and call the elements of $C$  \tbf{real colors}. Let $\bV$ be a set (disjoint from $V$) of $2p$ vertices we will call \tbf{ Platonic\footnote{We use the word Platonic in an allusion to the Platonic theory of forms, which holds that there exist abstract objects that represent real objects in an idealized way, but that are not themselves part of the physical world.} vertices}, and let $\bC$ be a set (disjoint from $C$) of $\binom{2p}{2}$ \tbf{ Platonic colors}. An edge $e \in \binom{V}{2}$ is called a \tbf{real edge}, and an edge $\be \in \binom{\bV}{2}$ is called a \tbf{Platonic edge}. We will use the Platonic vertices to represent sets of real vertices (we are mostly interested in sets of $s \le p$ real vertices but sometimes we will consider the union of two such sets which is why we let $|\bV|=2p$). We will use colorings of the Platonic edges to represent the coloring of the corresponding real edges. In such a coloring of the Platonic edges, we may use some Platonic color on some subset of the edges to indicate that they all have the same real color without specifying what that real color is. This use of Platonic colors is crucial for our analysis for reasons we will explain later. We will also use the framework of Platonic vertices, edges and colors to indicate when a color is available at an uncolored edge. More formally, define the following:

\begin{defn}
For a set $\bS \subseteq \bV$ we define a \tbf{ type on $\bS$} to be a labeling $z$ of $\binom{\bS}{2}$  where each label is of one of the following forms:
\begin{enumerate}[(i)]
    \item $\col{c}$ for some $c \in C \cup \bC$,
    \item  $\unc{c}$ for some $c \in C \cup \bC$,
    \item $\unc{?}$,
    \item $\wild$.
\end{enumerate}
 We say a type $z$ on $\bS$ is \tbf{ legal} if 
for every $\bS' \subseteq \bS$ with $|\bS'|  \le p$, $z$ has at most $R(|\bS'|)$ repeats among the colors of colored edges in $\binom{\bS'}{2}$ (i.e.\ the colors $c$ such that $z$ labels some edge in $\binom{\bS'}{2}$ the label $\col{c}$).
\end{defn} 
\noindent
We abuse notation and write $\bS(z)=\bS$ when $z$ is a type on $\bS$ (i.e.\ when we write $\bS(z)$ we mean the set of vertices such that  $z$ is a labeling of $\binom{\bS}2$).

\begin{defn}\label{def:ext}
For $ \bD  \subseteq \bV$  we call an injection $\phi:\bD \rightarrow V$ a \tbf{ partial embedding}  and say that the \tbf{ order} of $\phi$ is $|\bD|$.  We abuse notation and denote by $\bD(\phi)$ the domain $\bD$ of $\phi$. For a superset $\bS \supseteq \bD$ we say $\phi': \bS \rightarrow V$ is an \tbf{extension} of $\phi$ if $\phi'$ agrees with $\phi$ on $\bD$. For a type $z$ on $\bS$ we say that the extension $\phi'$ \tbf{fits} $z$ if there exists an injection $\psi$ from the set of colors used on colored edges in $z$ to the set of real colors $C$ such that
\begin{enumerate}[(i)]
    \item $\psi(c)=c$ for any real color $c \in C$ that appears on $z$,
    \item for each $\be\in \binom{\bS}2$, if $z(\be)=\col{c}$, then $\phi'(\be)$ is colored $\psi(c)$, and
    \item for each $\be\in \binom{\bS}2$, if $z(\be)=\unc{c}$, then the edge $\phi'(\be)$ has not been assigned a color, and $\psi(c)$ is available at $\phi'(\be)$.
\end{enumerate}
We call $\psi$ the \tbf{color map} for $\phi'$ (note that the color map is unique when $\phi'$ fits $z$). For a type $z$ let $\Ext(z, \phi)= \Ext(z, \phi, i)$ be the set of extensions $\phi'$ of $\phi$ that fit $z$ at step $i$.
We also say that $S \subseteq V$ \tbf{fits} $z$ if there exists some $\phi'$ fitting $z$ whose image is $S$. 
\end{defn}  

We will sometimes write the name of a set when we mean the cardinality of that set. In particular, we will use this convention when we will consider one-step changes, for example, we write $\E[\Delta \Ext(z, \phi)|\mc{H}_i]$ instead of $\E[\Delta |\Ext(z, \phi)||\mc{H}_i]$. 

We now explain the importance of Platonic colors. We will define a family of random variables, including certain variables of the form $\Ext(z, \phi)$, which we will track (i.e.\ give sharp estimates that hold w.h.p.) as the process evolves. In order to show we have the necessary concentration it is important that our tracked variables all be sufficiently large (like a positive power of $n$). We would like to track many extension variables $\Ext(z, \phi)$ including ones for types $z$ assigning colors to almost all their edges. If we used only real colors to represent such a type $z$, then $\Ext(z, \phi)$ may even have expectation going to 0 since extensions with so many specified real colors are unlikely.  By using some Platonic colors in a type $z$ instead of real colors, we get to keep the important information describing the way colors are repeated on the edges while discarding the unimportant information that specifies which real colors correspond to which Platonic colors. In turn there will be more extensions that fit the type $z$ which allows us to track $\Ext(z, \phi)$.

We use labels of the form  $\unc{c}$ because sometimes we would like to specify colors that are available at uncolored edges. However other times we will discard that information. More specifically we define the following:
\begin{defn}
Let $z$ be a type. We define the \tbf{colored portion} of $z$, $\colport (z)$, to be the type on $\bS(z)$ such that for every edge $\be'\in \binom{\bS(z)}2$,
\[
\colport (z)(\be') =\begin{cases}\wild &\text{ if }z(\be')=\unc{c}\text{ for any }c\in C\cup \bC,\\
z(\be')&\text{ otherwise.}
\end{cases}
\]
If $\colport(z)=z$, i.e.\ $z$ does not assign any available colors to uncolored edges, then we say that $z$ is a \tbf{colors-only type}.
\end{defn}

To analyze the evolution of the process, we must describe how an extension $\phi'$ that fits type $z$ at step $i$ may no longer fit type $z$ at step $i+1$ (and how it might now fit some other type). Towards that goal we define the following:

\begin{defn}\label{def:pred}
Let $z$ be a type on $\bS$, suppose $\be\in \binom{\bS}2$, and $z(\be) = \col{c}$. If $c$ is a real color or if $c$ is a Platonic color that is repeated among the colored edges of $z$, then we let $\pre(z, \be)= \{z'\}$ where $z'$ is the unique type on $\bS(z)$ that agrees with $z$ except that $z'(\be) = \unc{c}$. Otherwise (i.e.\ when $c$ is a Platonic color that appears only once among the colored edges of $z$), we let $\pre(z, \be)$ be the set of all types $z'$ on $\bS$ such that for all $\be'\in\binom{\bS}2$,
\[
z'(\be') = \begin{cases} \unc{c'} \mbox{ if $\be'=\be$ or $z(\be')=\unc{c}$}\\
z(\be') \mbox{ otherwise}
\end{cases}
\]
for some real color $c'$. Note that in this case $|\pre(z, \be)|=|C|$ (the number of choices for $c'$). If $z' \in \pre(z, \be)$ for some edge $\be$ we call $z'$ a \tbf{predecessor} of $z$. 
\end{defn}

In our analysis we will consider ways to find extensions $\phi' \in \Ext(z, \phi)$, which generally involves finding suitable real vertices for the image of $\phi'$ such that certain edges are colored appropriately. In such a situation, certain Platonic colors may already be associated with real colors while other Platonic colors may not (in particular, any Platonic color assigned to an edge $\be'\subseteq \bD(\phi)$ must correspond to the real color given to its image $\phi(\be')$). Thus, when we count our extensions $\phi'$, we expect to find fewer of them when there are more requirements stipulating that the colors of edges must be repeats. The next two definitions give us some tools to (heuristically for now, formally later) count extensions $\phi'$ based on the considerations in this paragraph.  

\begin{defn}\label{def:coipow}
Suppose $z$ is a type  on $\bS$ and let $\bE' \subseteq \bE \subseteq \binom{\bS}{2}$ be sets of Platonic edges. Let $M = M(\bE, \bE', z)$ be the set of edges in $\bE \sm \bE'$ that are colored by $z$. Let $\free = \free(\bE, \bE', z)$ be the set of Platonic colors that $z$ assigns to edges of $\bE \sm \bE'$ but does not assign to any edge of $\bE'$. 
Define\footnote{$\coi$ here is short for ``coincidence''}
\begin{equation}\label{eqn:edgcoidef}
    \edgcoi(\bE, \bE', z):= |M| - |\free|.
\end{equation}
 We will also define a version of the above for vertex sets:
\begin{equation*}
    \coi(\bS, \bS', z):= \edgcoi\rbrac{\binom{\bS}{2}, \binom{\bS'} 2 , z}.
\end{equation*}
\end{defn}\noindent
Clearly $|\free| \le |M|$ so we always have that  $\edgcoi(\bE, \bE', z)$ and $\coi(\bS, \bS', z)$ are nonnegative. We often want to calculate $\edgcoi$ and $\coi$ based on repeats and in situations where $\bE'$ (from \eqref{eqn:edgcoidef}) contains no colored edges. In this case, note that $|M|-|\free|$ equals the number of repetitions in $\bE$ plus the number of real colors assigned to edges in $\bE$.

Say we would like a crude heuristic estimate of the number of extensions $\phi' \in \Ext(z, \phi)$ where $z$ is a type on $\bS$ and $\phi$ has domain $\bD$. Assume for now that $z$ does not specify any available colors at uncolored edges. There are $|\bS \sm \bD|$ many Platonic vertices to embed into $V$, and we must embed them in such a way that there are $\coi(\bS, \bD, z)$ many color coincidences among the newly embedded edges (i.e.\ if we inspect the newly embedded edges one by one there are $\coi(\bS, \bD, z)$ that are required to have the same color as a previously embedded or inspected edge). Heuristically, since our coloring should in many ways resemble a uniform random coloring, each color coincidence has a probability of $O(1)/|C| = \TT\rbrac{n^{-\f}}$. Thus we heuristically predict that (ignoring log powers and everything else besides the power of $n$)
\[
\Ext(z, \phi) \approx n^{|\bS \sm \bD| -\f \coi(\bS, \bD, z)}.
\]
This motivates the following definition:
\begin{defn}
Let 
\[
\pow(\bS, \bS', z) := |\bS \sm \bS'| - \f \coi(\bS, \bS', z).
\]
\end{defn}\noindent
Note that  $\pow(\bS, \bS', z)$ can always be written as a rational number whose denominator is $\binom{p}{2}-q+1$. Thus, for example, the smallest positive value $\pow$ can take is $\frac{1}{\binom{p}{2}-q+1}$.
Next we define restrictions of types, which we will use in the next paragraph.

\begin{defn}
Let $z$ be a type on $\bS$, and let $\bS' \subseteq \bS$. The \tbf{restriction of $z$ to $\bS'$}, denoted $z|_{\bS'}$, is the type on $\bS'$ that agrees with $z$ on all edges of $\bS'$. 
\end{defn}

While the $\pow$ function will be very useful, it is not always true that it gives the correct power of $n$ when we count extensions. Indeed, we will have to deal with situations where $\pow$ gives a negative number, which does not necessarily mean that the extensions in question do not exist. This can complicate the situation when a type may have an ``unlikely part" (i.e.\ a subset such that there are probably no extensions to that subset) even if the whole type has a positive $\pow$. A bit more formally, suppose again that $\phi$ has domain $\bD$ and we want to count extensions $\phi' \in \Ext(z, \phi)$ where $z$ is a type on $\bS$. If there is some $\bS'$ with $\bD \subseteq \bS' \subseteq \bS$ such that $\pow(\bS', \bD, z)<0$ then it would seem that there is probably no extension in $\Ext(z|_{\bS'}, \phi)$ and therefore none in $\Ext(z|_{\bS}, \phi)$. However, if there is even one extension in $\Ext(z|_{\bS'}, \phi)$, then it seems there should be about $n^{\pow(\bS, \bS', z)}$ extensions in $\Ext(z|_{\bS'}, \phi)$. This motivates the following definition. 
\begin{defn} For $\bS'' \subseteq \bS \subseteq \bV$ we let
 \[
 \maxpow(\bS, \bS'', z) := \max \cbrac{ \pow(\bS, \bS', z): \bS'' \subseteq \bS' \subseteq \bS }.
 \]
\end{defn}

Note that $\maxpow$ is always nonnegative since $\pow(\bS, \bS, z)=0$. 
To analyze our process, it will be very important to pay attention to how colors become forbidden at edges. To model that using Platonic vertices we define the following. 

\begin{defn}\label{def:ecpreforb}
Let $\be$ and $\be'$ be distinct Platonic edges and $c,c'$ be colors (real or Platonic). A \tbf{$(\be,c,\be',c')$-preforbidder} is a type $y$ on $\bS\supseteq \be\cup\be'$ such that the following all hold.
\begin{enumerate}[(i)]
    \item \label{def:ecpreforb1}$y$ does not assign a color to $\be$ or $\be'$.
    \item \label{def:ecpreforb2} Let $y_\be$ be the type that agrees with $y$ except $y_\be(\be) = \col{c}$, and similarly let $y_{\be'}$ agree with $y$ except $y_{\be'}(\be') = \col{c'}$. Then $y_{\be}$ and $y_{\be'}$ are both legal.
    \item \label{def:ecpreforb3}  Let $y_{\be, \be'}$ be the type that agrees with $y$ except we have both $y_{\be, \be'}(\be) = \col{c}$ and $y_{\be, \be'}(\be') = \col{c'}$. Then $y_{\be, \be'}$ is illegal.
    \item \label{def:ecpreforb4}the restriction of $y_{\be, \be'}$ to any set $\bS'$ with $\be \cup\be'\subseteq \bS'\subsetneq \bS$ is legal, and
    \item \label{def:ecpreforb5}there are at most three real colors appearing on the colored edges of $y$.
\end{enumerate}
We call $y$ a \tbf{preforbidder} if it is a $(\be,c,\be',c')$-preforbidder for some $(\be,c,\be',c')$.
\end{defn}
Note that the above definition has some symmetry, and in particular we have the following:
\begin{obs}\label{obs:preforbswap}
A type $y$ is a $(\be,c,\be',c')$-preforbidder if and only if it is a $(\be',c',\be,c)$-preforbidder.
\end{obs}

Recall the definition \eqref{eqn:Rdef}:
\[
R(s):= \left \lceil \frac{(s-2)\rbrac{\binom{p}{2}-q+1}}{p-2} \right \rceil -1.
\]
Since the argument of the ceiling function above is a fraction with denominator $p-2$, we have
\begin{equation}\label{eqn:Rbounds}
 \frac{(s-2)\rbrac{\binom{p}{2}-q+1}}{p-2} - 1 \le   R(s) \le \frac{(s-2)\rbrac{\binom{p}{2}-q+1}}{p-2} - \frac{1}{p-2}.
\end{equation}

Next we use the above to bound the $\pow$ of preforbidders:

\begin{obs}\label{obs:preforbpow}
Suppose $y$ is a $(\be,c,\be',c')$-preforbidder on $\bS$ and there is at least one real color in $\{c, c'\}$ (including the possibility that $c=c'$ is real). Then 
\[
\pow(\bS, \be, y) \le \f.
\]
\end{obs}
\begin{proof}
Since $y$ does not assign a color to $\be$ or $\be'$ we have $\pow(\bS, \be, y) = \pow(\bS, \be', y)$. Thus WLOG we may assume that $c$ is real. 

Recall that by Definition \ref{def:ecpreforb} \eqref{def:ecpreforb2} and \eqref{def:ecpreforb3} we have that $y_{\be, \be'}$ is illegal, while $y_{\be}$ is legal. Also any proper restriction of $y_{\be, \be'}$ is legal. Therefore $y_{\be, \be'}$ must have exactly $R(|\bS|)+1$ repeats and $y_{\be}$ has exactly $R(|\bS|)$ repeats. Therefore, either the multiset of colors assigned by $y$ has $R(|\bS|)$ repeats and does not include $c$, or else it has $R(|\bS|)-1$ repeats and includes $c$. In either case we have $\coi(\bS, \be, y)\ge R(|\bS|)$ and so 
\[
\pow(\bS, \be, y) \le|\bS|-2 - \f R(|\bS|)\le  \f.
\]
\end{proof}

We will be especially interested in modeling the mechanism through which a particular real color $c$ becomes forbidden at an edge. Thus we define the following:

\begin{defn}\label{def:trackablepreforb}
For a Platonic edge $\be$ and a real color $c$, we let $\prefbd(\be, c)$ be the set of all triples $(z, \be', c')$ such that $z$ is a $(\be, c, \be', c')$-preforbidder with the following properties:
\begin{enumerate}[(i)]
    \item $z(\be')=\unc{c'}$,
    \item either $c'=c$ or $c'$ is Platonic,
    \item other than $c$, all colors used on colored edges of $z$ are Platonic,
    \item other than the edge $\be$, $z$ does not use the labels $\unc{?}$ or $\wild$,
    \item  replacing all of the labels $\unc{c''}$ with $\col{c''}$ except at $\be$ would result in a legal type, and
    \item if $z$ has any edge label $\unc{c''}$ for a Platonic color $c''$ then $z$ also has some edge label $\col{c''}$.
\end{enumerate}
\end{defn}
To estimate the probability that an available color $c$ becomes forbidden at $e$ at step $i$, we would like to estimate the number of ways that can happen. Thus we will track $\Ext(z, \phi)$ for all $z, \phi$ such that $(z, \be', c') \in \prefbd(\be, c)$ for some $\be', c', c$ and $\bD(\phi)=\be$. But in order to track those extension variables, we will also need to track extensions of all their predecessors. More formally we will track the following family of variables:

\begin{defn}\label{def:trackable} We say that $(y, \be)$ is a \tbf{rooted trackable type} if $y$ can be obtained by iteratively taking predecessors of some type $z$ such that $(z, \be', c') \in \prefbd(\be, c)$ for some $\be', c', c$. (This includes the possibility of not taking any predecessors, i.e.\ $y=z$).
We say that $y$ is a \tbf{trackable type} if there exists some $\be$ such that $(y, \be)$ is a rooted trackable type. In that case we may also say that $y$ is a \tbf{trackable type with root} $\be$. We call $\be$ the \tbf{root} of $(y,\be)$.
\end{defn}\noindent
When $y$ is a trackable type, we will abuse notation and write $\be(y)$ to denote an edge $\be$ such that $(y,\be)$ is a rooted trackable type, if there is a choice of edges $\be$, $\be(y)$ will be one chosen arbitrarily. We will track $\Ext(y, \phi)$ for all trackable types $y$ and all $\phi$ of order 2 such that $(y,\bD(\phi))$ is a rooted trackable type. For other pairs $(y, \phi)$ we will use only a crude upper bound on $\Ext(y, \phi)$. We now list some useful properties of trackable types.
\begin{obs}\label{obs:trackable} Suppose  $(y, \be)$ is a rooted trackable type on $\bS$. Then we have the following:
\begin{enumerate}[(i)]
    \item \label{obs:trackable1} replacing all of the labels $\unc{c''}$ with $\col{c''}$ except at $\be$ would result in a legal type,
   \item \label{obs:trackable2}    $\pow(\bS', \be , y)\ge \m$ for all $\be  \subsetneq \bS' \subseteq \bS$,
    \item \label{obs:trackable3} there is at most one real color $c''$ such that $y$ labels any edge $\col{c''}$,
    \item \label{obs:trackable4}  other than possibly at the edge $\be$, $y$ does not use the labels $\unc{?}$ or $\wild$,
    \item \label{obs:trackable5} if $y$ has any edge label $\unc{c''}$ for a Platonic color $c''$ then $y$ also has some edge label $\col{c''}$.
\end{enumerate}
\end{obs}

\begin{proof}
    Parts \eqref{obs:trackable1} and \eqref{obs:trackable3}-\eqref{obs:trackable5} follow from Definitions \ref{def:pred}, \ref{def:trackablepreforb} and \ref{def:trackable}. Indeed, if $(z, \be', c') \in \prefbd(\be, c)$ then by Definition \ref{def:trackablepreforb} $z$ satisfies \eqref{obs:trackable1} and \eqref{obs:trackable3}-\eqref{obs:trackable5}. Furthermore, properties \eqref{obs:trackable1} and \eqref{obs:trackable3}-\eqref{obs:trackable5} are all preserved when we take predecessors by Definition \ref{def:pred}. 
    
    Now we prove \ref{obs:trackable2}. Since this property is also preserved under taking predecessors, we will assume that $(y, \be', c') \in \prefbd(\be, c)$. Suppose $\be  \subsetneq \bS' \subseteq \bS$. Since $y$ is legal, there are at most $R(|\bS'|)$ repeats among the colored edges in $\binom{\bS'}{2}$. If there is no real color appearing on a colored edge in $\binom{\bS'}{2}$ then we have $\coi(\bS', \be, y) \le R(|\bS'|)$ and hence
    \[
    \pow(\bS', \be, y) \ge |\bS'|-2 - \f R(|\bS'|) \ge \m
    \]
    so we are done. If on the other hand there is a real color appearing on a colored edge in $\binom{\bS'}{2}$ then that color must be $c$. We claim that in this case there can only be at most $R(|\bS'|)-1$ repeats among the colored edges in $\binom{\bS'}{2}$. Indeed, if there were $R(|\bS'|)$ repeats then the type $y'$ formed by replacing the label of $\be$ with $\col{c}$ would be illegal, which is a contradiction to Definition~\ref{def:ecpreforb}~\eqref{def:ecpreforb2}. Thus $\bS'$ has at most $R(|\bS'|)-1$ repeats and one real color, so again we have $\coi(\bS', \be, y) \le R(|\bS'|)$ and we are done. 
\end{proof}

Section \ref{sec:platonic} will be devoted to proving four technical lemmas about preforbidders and other types. We state them here since we will use them before proving them.

\begin{lem:preforbpowbound}
Let $y$ be a $(\be',c',\be'',c'')$-preforbidder on $\bS$. Let $\a=\a(\bS')$ be the number of Platonic colors in $\{c', c''\}$ that appear in $\bS$ but not in $\bS'$.  Then for all $\bS'$ with $\be'\cup\be''\subseteq \bS'\subsetneq \bS$,
\begin{equation*}
    \pow(\bS,\bS',y)\leq \f \a  - \m
\end{equation*}
\end{lem:preforbpowbound}

The other two technical lemmas are about unions of types, which we define before stating the lemmas:
\begin{defn}\label{def:compatible}
Let $y_1, y_2$ be types on $\bS_1, \bS_2$ respectively. We say $y_1$ and $y_2$ are \tbf{compatible} if for every edge $\be' \subseteq \bS_1 \cap \bS_2$, we have that 
\[
y_1(\be') = \col{c} \Longleftrightarrow y_2(\be') = \col{c}.
\]
In other words, $\colport(y_1)|_{\bS_1 \cap \bS_2} = \colport(y_2)|_{\bS_1 \cap \bS_2}$.

If $y_1$ and $y_2$ are compatible, we define the \tbf{union of types} $y_1 \cup y_2$ to be the unique colors-only type $y$ such that $y|_{\bS_1} = \colport(y_1)$, $y|_{\bS_2} = \colport(y_2)$, and $y$ only assigns colors to edges within $\bS_1$ and to edges within $\bS_2$.
\end{defn}

\begin{lem:appendpreforb}
Let $y_1$ be a trackable type on $\bS_1$ with root $\be_1$, and let $y_2$ be a $(\be_2,c_2,\be_3,c_3)$-preforbidder on $\bS_2$ that is compatible with $y_1$. Assume $\be_2,\be_3\subseteq \bS_1$ and that $y_1(\be_2)=\unc{c_2}, y_1(\be_3)=\unc{c_3}$. Assume $\bS_2\setminus\bS_1\neq \emptyset$. Then 
\[
\maxpow(\bS_1\cup\bS_2,\be_1 ,y_1\cup y_2)\leq \pow(\bS_1,\be_1,y_1)-\m.
\]
\end{lem:appendpreforb}

\begin{lem:union}
Let $y_1$ be a trackable type on $\bS_1$ with root $\be_1$, and let $y_2$ be a $(\be_2,c_2,\be_3,c_3)$-preforbidder on $\bS_2$ that is compatible with $y_1$. Assume $\be_3\subseteq \bS_1$, $\be_3 \neq \be_1$, and $\be_2 \not \subseteq \bS_1$. Assume $y_1(\be_3)=\unc{c_3}$. Assume  $c_2$ is a real color. Then 
    \[
    \maxpow(\bS_1 \cup \bS_2, \be_1 \cup \be_2, y_1 \cup y_2) \le  \pow(\bS_1, \be_1, y_1) - \m.
    \]
\end{lem:union}

\section{Trajectories of our random variables}\label{sec:traj}

In this section we describe the trajectories for the variables $\Ava(e)$ and $\Ext(z, \phi)$ for trackable $z$ and $\phi$ of order 2. Rather than simply giving formulas we attempt to provide a heuristic derivation for them. 

Say we have colored $i=\binom{n}{2} t$ edges so far. Our heuristic assumptions are as follows: $t$ is the probability that an edge has been colored, $1/|C|$ is the probability that a colored edge has been assigned any particular color, and some number $a=a(t)$ is the probability that a particular color is available at a particular uncolored vertex. We also generally assume that events are approximately independent if we see no significant reason why they should be dependent. Under those heuristic assumptions we anticipate that 
\[
\Ava(e) = \Ava(e, i) \approx \ava(t) := |C| a(t)
\]
for all $e \in U_i$ and for some deterministic function $a(t)$.

 Assume $\phi$ has order 2, and $z$ is a trackable type with $|\bS(z)|=s$ that  has $\ell$ colored edges, $r$ repeats, and $k \in \{0, 1\}$ real colors among its colored edges. Then we heuristically predict
\begin{align}
\Ext(z, \phi) \approx    \ext_z(t) := &n^{s-2}|C|^{\ell-r-k} t^\ell \rbrac{\frac 1{|C|}}^\ell (1-t)^{\binom{s}{2}-\ell-1} a(t)^{\binom{s}{2}-\ell-1} \nonumber \\
& = n^{s-2}|C|^{-r-k}a(t) ^{\binom{s}{2}-\ell-1} t^\ell (1-t)^{\binom{s}{2}-\ell-1} .\label{eqn:exttraj}
\end{align}

We will now derive a differential equation for $a(t)$ by heuristically estimating the one-step change in $\Ava(e)$. We estimate the probability that a fixed real color $c$ becomes forbidden at an edge $e$ in a single step (assuming that $c$ was available at $e$ one step previously). First, note that the probability that any particular edge $e'$ gets any assigned any particular color $c' \in \Ava(e')$ at step $i$ is approximately $1/[(1-t) \binom n2 |C|a]$. So now we count pairs $(e', c')$ such that $c' \in \Ava(e')$ and assigning $c'$ to $e'$ would forbid $c$ at $e$. 
Let $(z, \be', c^*) \in \prefbd(\be, c)$. Each extension  $\phi' \in \Ext(z, \phi)$ counts the pair $(e', c')$ where $e'=\phi'(\be')$ and $c'$ is the color corresponding to $c^*$ (i.e.\ $c'=c^*$ if $c^* \in C$ or otherwise $c'$ is the real color representing $c^*$ in the extension $\phi'$). But each such pair $(e', c')$ is counted about $\equiv(z, \be') \ava(t)^{\binom{s}{2}-\ell-2}$ times, where $\equiv(z, \be')$ counts the number of equivalent ways to represent a clique fitting type $z$ as a preforbidder in $\prefbd (c, \be)$. More precisely, $\equiv(z, \be')$ is the number of ways to permute the Platonic colors in $z$ times the number of bijections from $\bS(z)$ to itself fixing both vertices in $\be$ and fixing the edge $\be'$ (but possibly swapping the endpoints of $\be'$) and preserving the colors of colored edges (not necessarily preserving the available colors at uncolored edges). Thus our estimate for the number of pairs $(e', c')$ is (below, every term is indexed by some type $z$, and $\ell=\ell(z)$ is always the number of edges colored by that $z$, $r=r(z)$ is the number of repeats in $z$, and $k=k(z) \in \{0, 1\}$ is the number of real colors on colored edges in $z$)

\begin{align}
 &  \sum_{\substack{(z, \be', c^*)  \in \prefbd(\be, c)  }} \frac{ \ext_z(t)}{\equiv(z, \be') \ava(t) ^{\binom{s}{2}-\ell-2}} \label{eqn:sumpreforb}\\
     & =\sum_{ 3 \le s \le p} \;\; \sum_{\substack{(z, \be', c^*)  \in \prefbd(\be, c)    \\ |\bS(z)| = s}} \frac{  n^{s-2}|C|^{-r-k}a(t) ^{\binom{s}{2}-\ell-1} t^\ell (1-t)^{\binom{s}{2}-\ell-1}  }{\equiv(z, \be') (|C|a(t)) ^{\binom{s}{2}-\ell-2}} \nonumber\\
        & =\sum_{ 3 \le s \le p} \;\; \sum_{\substack{(z, \be', c^*)  \in \prefbd(\be, c)    \\ |\bS(z)| = s}} \frac{  n^{s-2}|C|^{\ell-r-k-\binom{s}{2}+2}a(t)  t^\ell (1-t)^{\binom{s}{2}-\ell-1}  }{\equiv(z, \be') } \label{eqn:traj1}
\end{align}
Now we claim that the number of triples $(z,\be',c^*)\in\prefbd(\be, c)$ such that $|\bS(z)|=s,  \ell(z)=\ell$ and $\equiv(z,\be',c^*)=\eta$ is $\mu(\ell,\eta, s)\cdot |C|^{\binom{s}{2}-\ell-2}$ for some $\mu(\ell,\eta, s)$ not depending on $n$. Indeed, for fixed $\ell, \eta, s$ there are some constant number of ways to choose the edges that are colored and the colors on those edges (note that the only real color that can show up on a preforbidder in $\prefbd(\be, c)$ is $c$), a constant number of ways to choose $c^*$ (this color must be Platonic or must be $c$), and then for each uncolored edge that is not $\be$ or $\be'$, there are $|C|$ choices for the color available there, giving us the claim. In addition, note that for all $(z, \be', c^*) \in \prefbd(\be, c)$ we have that $r+k = R(s)$. Therefore \eqref{eqn:traj1} becomes 
\begin{align}
& \sum_{\substack{ 3 \le s \le p \\ 1\leq \ell\leq \binom{s}2 -2\\ 1\leq \eta\leq s! |\bC|!}} \mu(\ell,\eta, s)\cdot |C|^{\binom{s}{2}-\ell-2} \cdot  \frac{  n^{s-2}|C|^{\ell- \binom{s}{2}+2-R(s) }a(t)  t^\ell (1-t)^{\binom{s}{2}-\ell-1}  }{\eta } \nn\\
& = \sum_{\substack{ 3 \le s \le p \\ 1\leq \ell\leq \binom{s}2 -2\\ 1\leq \eta\leq s! |\bC|!}}  \frac{ \mu(\ell,\eta, s) n^{s-2}|C|^{-R(s) }a(t)  t^\ell (1-t)^{\binom{s}{2}-\ell-1}  }{\eta } \label{eqn:traj2}
\end{align}
Now for a value $s$ such that $p-2$ divides $(s-2)\rbrac{\binom{p}{2}-q+1}$, we have $R(s)=  \frac{(s-2)\rbrac{\binom{p}{2}-q+1}}{p-2}  -1$ and so  
\[
n^{s-2}|C|^{-R(s)} = n^{s-2}|C|\rbrac{\k n^{\f } \log^{-\m} n}^{-\frac{(s-2)\rbrac{\binom{p}{2}-q+1}}{p-2}} = \k^{-\frac{(s-2)\rbrac{\binom{p}{2}-q+1}}{p-2}} |C| \log^\frac{s-2}{p-2} n .
\]
On the other hand if $p-2$ does not divide $(s-2)\rbrac{\binom{p}{2}-q+1}$, then $R(s) \ge   \frac{(s-2)\rbrac{\binom{p}{2}-q+1}}{p-2}  -1 + \frac{1}{p-2}$ and we get
\[
n^{s-2}|C|^{-R(s)} = \TO\rbrac{n^{s-2}|C| \rbrac{n^{\f } }^{-\frac{(s-2)\rbrac{\binom{p}{2}-q+1}}{p-2}  - \frac{1}{p-2}} } = \TO\rbrac{ n^{-\mf } }.
\]
Thus, \eqref{eqn:traj2} becomes
\begin{align}
&  \sum_{\substack{ 3 \le s \le p \\ (p-2) \; | \; (s-2)\rbrac{\binom p2 - q + 1} \\ 1\leq \ell\leq \binom{s}2 -2\\ 1\leq \eta\leq s! |\bC|!}}  \frac{ \mu(\ell,\eta, s) \k^{-\frac{(s-2)\rbrac{\binom{p}{2}-q+1}}{p-2}}|C| \log^\frac{s-2}{p-2} n \cdot a(t)  t^\ell (1-t)^{\binom{s}{2}-\ell-1}  }{\eta }  + \TO\rbrac{ n^{-\mf } } \nn\\
&= |C| a(t) (1-t)\sum_{\substack{ 3 \le s \le p \\ (p-2) \; | \; (s-2)\rbrac{\binom p2 - q + 1} \\ 1\leq \ell\leq \binom{s}2 -2\\ 1\leq \eta\leq s! |\bC|!}}  \frac{ \mu(\ell,\eta, s) \k^{-\frac{(s-2)\rbrac{\binom{p}{2}-q+1}}{p-2}}    t^\ell (1-t)^{\binom{s}{2}-\ell-2}  }{\eta } \log^\frac{s-2}{p-2} n + \TO\rbrac{ n^{-\mf } } \nn\\
&= |C| a(t) (1-t)\sum_{\substack{ 3 \le s \le p \\ (p-2) \; | \; (s-2)\rbrac{\binom p2 - q + 1} }}  h_s(t) \log^\frac{s-2}{p-2} n + \TO\rbrac{ n^{-\mf } }\label{eqn:hs}\\
&= |C| a(t) (1-t) h(t) + \TO\rbrac{ n^{-\mf } } \label{eqn:h}
\end{align}
where on line \eqref{eqn:hs} we use the definition
\begin{equation*}
h_s(t) :=  \sum_{\substack{  1\leq \ell\leq \binom{s}2 -2\\ 1\leq \eta\leq s! |\bC|!}}  \frac{ \mu(\ell,\eta, s) \k^{-\frac{(s-2)\rbrac{\binom{p}{2}-q+1}}{p-2}}    t^\ell (1-t)^{\binom{s}{2}-\ell-2}  }{\eta } 
\end{equation*}
and on line \eqref{eqn:h} we use 
\begin{equation}\label{eqn:hdef}
 h(t) := \sum_{\substack{ 3 \le s \le p \\ (p-2) \; | \; (s-2)\rbrac{\binom p2 - q + 1} }}  h_s(t) \log^\frac{s-2}{p-2} n.   
\end{equation}

Note that $h_s(t)$ is a polynomial in $t$ with constant coefficients (not depending on $n$).

We derive a differential equation by estimating the one-step change in $\Ava(e)$ in two ways. First, since $\Ava(e) \approx |C|a(t)$ we should have
\[
\D \Ava(e) \approx |C|a'(t) \D t = \frac{1}{\binom n2} |C|a'(t).
\]
Second, we can estimate $\D \Ava(e)$ by considering each available color and the probability that it becomes forbidden in one step. Since there are $\Ava(e) \approx |C|a(t)$ available colors and each one becomes forbidden with probability about 
\[
\frac{1}{(1-t)\binom n2 |C|a(t)} \cdot |C|a(t)(1-t)h(t),
\]
since for each $c \in \Ava(e)$ there are about $|C|a(t)(1-t)h(t)$ pairs $(e', c')$ such that assigning $c'$ to $e'$ would forbid $c$ at $e$, and each such pair has a probability of about $1/\sbrac{(1-t)\binom n2 |C|a(t)}$ of being the choice made in one step. Thus, we heuristically expect
\begin{align*}
    &\frac{1}{\binom n2} |C|a'(t) = -|C|a(t) \cdot \frac{1}{(1-t)\binom n2 |C|a(t)} \cdot |C|a(t)(1-t)h(t)
\end{align*}
or
\begin{equation*}
    a'(t) = - a(t) h(t).
\end{equation*}
Since $a(0)=1$ we have
\begin{equation}\label{equation definition of a}
    a(t) = \exp\cbrac{-H(t) }
\end{equation}
where
\begin{equation}\label{eqn:Hdef}
  H(t) := \int_0^t h(\tau) d\tau 
\end{equation}

Note that we will take \eqref{equation definition of a} as the definition of $a(t)$.

Let
\begin{equation*}
  \eps := 10^{-3}p^{-6}.
\end{equation*}
(this choice will appear arbitrary at the moment. The reader can safely think of $\eps$ as just a small positive constant). 
Note that for $3 \le s \le p$ we have
\[
h_s(t) \le  \k^{-\frac{\binom{p}{2}-q+1}{p-2}} \sum_{\substack{  1\leq \ell\leq \binom{s}2 -2\\ 1\leq \eta\leq s! |\bC|!}}  \frac{ \mu(\ell,\eta, s)    }{\eta } 
\]
and so 
\begin{align}
   h(t)&= \sum_{ 3 \le s \le p} h_{s}(t) \log^\frac{s-2}{p-2}n \nn\\
   & \le \log n \cdot \sum_{3 \le s \le p} h_{s}(t) \nn\\
   &\le \log n \cdot \k^{-\frac{\binom{p}{2}-q+1}{p-2}} \sum_{\substack{3 \le s \le p \\  1\leq \ell\leq \binom{s}2 -2\\ 1\leq \eta\leq s! |\bC|!}}  \frac{ \mu(\ell,\eta, s)    }{\eta } \nn\\
   & \le \eps \log n \label{eqn:hbound}
\end{align}
where on the last line we have used the (as yet unstated) assumption that we will choose $\k$ so that
\begin{equation}\label{eqn:pickkappa}
     \k \ge \rbrac{\frac{\eps}{\displaystyle\sum_{\substack{3 \le s \le p \\  1\leq \ell\leq \binom{s}2 -2\\ 1\leq \eta\leq s! |\bC|!}}  \frac{ \mu(\ell,\eta, s)    }{\eta }}}^{-\frac{p-2}{\binom{p}{2}-q+1}}
\end{equation}
Note that by \eqref{eqn:hbound} and \eqref{eqn:Hdef} we have 
\[
H(t) \le \eps t \log n \qquad \mbox{ for all } 0 \le t \le 1
\]
and so
\[
a(t) \ge \exp\cbrac{-\eps t \log n } \ge n^{-\eps} \qquad \mbox{ for all } 0 \le t \le 1.
\]
\section{The good event}\label{sec:goodevent}

In this section we define  the good event $\mc{E}_i$, which among other things stipulates that every uncolored edge still has plenty of available colors. If we manage to show $\mc{E}_{i^*}$ holds w.h.p.\ then that means the process manages to color at least $i^*$ edges before terminating. More specifically, $\mc{E}_i$ will stipulate that all of our tracked variables are within a small window of their respective trajectories we derived in Section \ref{sec:traj}. $\mc{E}_i$ will also stipulate some crude upper bounds on certain other extension variables.

We define our \tbf{error functions} 
 \begin{align}
     f_\Ava(t) &: = n^{-\hm + 20p^2\eps} (1-t)^{-5p^2} e^{10p^4 \eps t \log n} \label{eqn:pickfava}\\
     f_\Ext(t) &: = n^{-\hm + 20p^2\eps} (1-t)^{-5p^2+1}a(t)^{-1}e^{10p^4 \eps t \log n} \label{eqn:pickfext}\\
     f_z(t) &: = a(t)^{\binom s2 - 1 - \ell}f_\Ext(t) \label{eqn:pickfz}
 \end{align}
 where in the last line $z$ is a type on $s$ vertices with $\ell$ colored edges.

\begin{defn}\label{definition of the good event}
The \tbf{good event} $\mc{E}_i$ to be the event that for all $i' \le i$ we have the conditions below. We let $t'=i'/\binom n2$.
 \begin{enumerate}[(i)]
     \item \label{cond:atraj} For all $e \in U_{i'}$ we have
     \[
      |\Ava(e, i')- \ava(t')| \le |C| f_\Ava(t')
     \]
     \item \label{cond:exttraj} For all rooted trackable types $(z, \be)$ and for all $\phi$ of order 2 mapping $\be$ to an edge $e \in U_{i'}$, we have
     \[
     |\Ext(z, \phi, i')- \ext_z(t')| \le n^{s-2}|C|^{-r-k}f_{z}(t')
     \]
     \item \label{cond:maxpow} Let $z$ be a type. Let $\bS:=\bS(z)$ and $\phi$ be a partial embedding with domain $\bD:=\bD(\phi)$. Let $s:=|\bS|, d:=|\bD|$ and let $m$ be the number of edges not inside $\bD$ that are assigned colors by $z$. Suppose further that we have the property:
\begin{equation}\label{prop:pow}
   \mbox{ for all }\bS' \mbox{ such that }\bD \subsetneq \bS' \subseteq \bS \mbox{ we have that } \pow(\bS, \bS', z) \le 0.
\end{equation}
Then 
\[
\Ext(z, \phi, i') \le n^{\maxpow(\bS, \bD, z) + 10(m+s-d)\eps}.
\]
 \end{enumerate}
\end{defn}

We also define
\[
i_{max} := \rbrac{1-n^{-\eps}} \binom{n}{2}, \quad t_{max}:=i_{max}/ \binom{n}{2} = 1-n^{-\eps}.
\]

\begin{obs}
We have
\[
|C| f_\Ava(t) = o(\ava(t))
\]
for all $t\le t_{max}$. In other words, the bound in Condition \ref{cond:atraj} gives an asymptotically tight estimate of $\Ava(e)$ for all $t\le t_{max}$. 
\end{obs}

\begin{proof} Using the fact that $a(t)\geq n^{-\eps}$ and $1-t \ge n^{-\eps}$,  we have
    \begin{align*}
        \frac{|C| f_\Ava(t)}{\ava(t)} &= \frac{ f_\Ava(t)}{a(t)} \le \frac{n^{-\hm + 20p^2\eps} (1-t)^{-5p^2} e^{10p^4 \eps t \log n }}{a(t)} \le n^{-\hm + (10p^4 + 25p^2   + 1)\eps} = o(1)
    \end{align*}

\end{proof}

We will prove the following:
\begin{thm}
Fix $p$, let $\eps = 10^{-3}p^{-6}$ and choose $\k$ large enough so that \eqref{eqn:pickkappa} holds. Then $\mc{E}_{i_{max}}$ holds w.h.p.. In particular, the coloring process colors at least $i_{max}$ many edges. 
\end{thm}

A proof outline of the theorem is as follows. We will show that the complement of $\mc{E}_{i_{max}}$ is contained in the event that some member of a large family of supermartingales increases by some amount that is very unlikely. The probability that one of our supermartingales misbehaves can be bounded by Freedman's inequality, which we now state:
\begin{lem}\label{lem:Freedman}
Let $Y(i)$ be a supermartingale with $\Delta Y(i) \leq D$ for all $i$, and let \newline $V(i) :=\displaystyle \sum_{k \le i} \Var[ \Delta Y(k)| \mc{H}_{k}]$.  Then,
\[
\P\left[\exists i: V(i) \le b, Y(i) - Y(0) \geq \lambda \right] \leq \displaystyle \exp\left(-\frac{\lambda^2}{2(b+D\lambda) }\right).
\] 
\end{lem}
The probability that $\mc{E}_{i_{max}}$ fails for any reason (i.e.\ due to any member of our family of supermartingales misbehaving) can then be bounded using the union bound. 

\section{The crude bounds on \texorpdfstring{$\Ext$}{2} variables}\label{sec:crude}

In this section we bound the probability that $\mc{E}_{i_{max}}$ fails due to Condition \eqref{cond:maxpow}. We also prove some bounds that are implied by Condition \eqref{cond:maxpow}. In particular, the following lemma will give us that in the good event, the conclusion of Condition~\eqref{cond:maxpow} holds for all extension variables, even those that do not satisfy \eqref{prop:pow}.

\begin{lem}\label{lem:maxpow}
Suppose $z$ is a type on $\bS$ and  $\phi$ is a partial embedding with domain $\bD$. Let $|\bS|=s, |\bD|=d$, and suppose $z$ assigns colors to $m$ edges that are not in $\bD$. Then in the good event we (deterministically) have 
\[
\Ext(z, \phi) \le n^{\maxpow(\bS, \bD, z) + 10(m+s-d)\eps}
\]
\end{lem}

\begin{proof}
We will define a sequence of sets $\bD=\bD_1 \subseteq \bD_2 \subsetneq \ldots \subsetneq \bD_j = \bS$, and let $z_{j'}:=z|_{\bD_{j'}}$. We will then bound $\Ext(z, \phi)$ by multiplying the number of ways to extend an embedding of $\bD_{j'-1}$ fitting $z_{j'-1}$ to an embedding of $\bD_{j'}$ fitting $z_{j'}$ for $j'=2, \ldots, j$. We let
\[
\bD_{2} \in \arg \min_{\bD'} \{\pow(\bD', \bD, z): \bD \subseteq \bD' \subseteq \bS \}
\]
and for $j' \ge 3$
\[
\bD_{j'} \in \arg \min_{\bD'} \{\pow(\bD', \bD_{j'-1}, z): \bD_{j'-1} \subsetneq \bD' \subseteq \bS \}.
\]

 We claim that for all $2 \le j' \le j$ the triple $(\bD_{j'-1}, \bD_{j'}, z_{j'})$ has the property on line \eqref{prop:pow}. Indeed, let $\bD_{j'-1} \subsetneq \bD' \subseteq \bD_{j'}$ and suppose that $\pow(\bD_{j'}, \bD', z_{j'}) = \pow(\bD_{j'}, \bD', z)>0$. Then 
\[
\pow(\bD', \bD_{j'-1}, z) = \pow(\bD_{j'}, \bD_{j'-1}, z) - \pow(\bD_{j'}, \bD', z) <  \pow(\bD_{j'}, \bD_{j'-1}, z) 
\]
contradicting our choice of $\bD_{j'}$. Thus $\pow(\bD_{j'}, \bD', z_{j'}) \le 0$ and so the triple $(\bD_{j'-1}, \bD_{j'}, z_{j'})$ has Property \eqref{prop:pow}. 

We claim that $\maxpow(\bD_2, \bD_1, z_2) = 0$. Indeed, if $\pow(\bD_2, \bD', z_2) = \pow(\bD_2, \bD', z)>0$ for any $\bD_1 \subseteq \bD' \subseteq \bD_2$ then we would have 
\[
\pow(\bD', \bD_1, z) = \pow(\bD_2, \bD_1, z) - \pow(\bD_2, \bD', z)< \pow(\bD_2, \bD_1, z),
\]
contradicting our choice of $D_2$. Now we claim that $\maxpow(\bD_{j'},\bD_{j'-1}, z_{j'}) =  \pow(\bD_{j'}, \bD_{j'-1}, z)$ for all $j' \ge 3.$ Indeed we already showed that if $j' \ge 3$ then $\pow(\bD_{j'}, \bD', z_{j'}) \le 0$ for all $\bD'$ with $\bD_{j'-1} \subsetneq \bD' \subseteq \bD_{j'}$, and so if there were some $\bD'$ with $\pow(\bD_{j'}, \bD', z_{j'}) >  \pow(\bD_{j'}, \bD_{j'-1}, z)$ then $\pow(\bD_{j'}, \bD_{j'-1}, z)<0$, but then
\[
\pow(\bD_{j'}, \bD_{j'-2}, z) = \pow(\bD_{j'}, \bD_{j'-1}, z) + \pow(\bD_{j'-1}, \bD_{j'-2}, z) < \pow(\bD_{j'-1}, \bD_{j'-2}, z)
\]
contradicting our choice of $\bD_{j'-1}$.

Putting this together, we multiply the upper bound given by Condition~\eqref{cond:maxpow} on the number of embeddings of $\bD_{j'}$ given an embedding of $\bD_{j'-1}$ for $j'=2, \ldots, j$. Let $d_{j'}:=|\bD_{j'}|$ and $m_{j'}$ be the number of edges not inside $\bD_{j'-1}$ that $z_{j'}$ assigns colors to. Then $|\Ext(z, \phi)|$ is at most
\begin{align*}
    & \prod_{j'=2}^j n^{\maxpow(\bD_{j'}, \bD_{j'-1}, z_{j'})  + 10(m_{j'} + d_{j'} - d_{j'-1})\eps } \\
    & = n^{0 + \pow(\bD_{3}, \bD_{2}, z) + \ldots + \pow(\bD_{j}, \bD_{j-1}, z) + 10(m_2 + \ldots + m_j + d_j - d_1)\eps} \\
    & = n^{  \pow(\bS, \bD_{2}, z) + 10(m+s-d)\eps}.
\end{align*}
Since by our choice of $\bD_2$, we have that $\maxpow(\bS, \bD, z) = \pow(\bS, \bD_{2}, z)$ we are done. 
\end{proof}

The next corollary follows from the fact that $m + s - d \le \binom p2 + p - 0 < 2p^2.$

\begin{cor}\label{cor:crudepractical}
    Suppose $z$ is a type on $\bS$ and  $\phi$ is a partial embedding with domain $\bD$. Then in the good event we (deterministically) have 
\[
\Ext(z, \phi) \le n^{\maxpow(\bS, \bD, z) + 20p^2\eps}
\]
\end{cor}

\begin{lem}\label{lem:Gpow}
Suppose $z$ is a type on $\bS$ and that $\phi$ is a partial embedding with domain $\bD:=\bD(\phi)\subseteq \bS$ such that $\phi$ fits the type $z|_{\bD}$.  Let $s:=|\bS|, d:=|\bD|$ and let $m$ be the number of edges not inside $\bD$ that are assigned colors by $z$. Suppose further that we have the property:
\begin{equation*}
   \mbox{ for all }\bS' \mbox{ such that }\bD \subsetneq \bS' \subseteq \bS \mbox{ we have that } \pow(\bS, \bS', z) \le 0.
\end{equation*}
Then the probability that there exists a step $i$ such that $\mc{E}_{i-1}$ holds and then at step $i$ we have
\[
\Ext(z, \phi) > n^{\maxpow(\bS, \bD, z) + 10(m+s-d)\eps}
\]
(i.e.\ we do not have $\mc{E}_i$ due to Condition \ref{cond:maxpow} failing for this $z, \phi$) is at most 
\[
\exp\cbrac{-\Omega\rbrac{n^\eps}}
\]
\end{lem}

\begin{proof}
First we handle the trivial case where $m=0$, i.e.\ $z$ does not color any edges outside $\bD$. Then $\maxpow(\bS, \bD, z)=s-d$ and trivially 
\[
\Ext(z, \phi) \le n^{s-d}
\]
so the lemma holds. Assume henceforth that $m \ge 1$.
 We define the variable
\begin{align*}
    \Ext^{+}(z, \phi)&=\Ext^{+}(z, \phi, i)\\
    &:=\begin{cases} 
& \Ext(z, \phi,i) - n^{\maxpow(\bS, \bD, z)-2 + 10 (m + s-d) \eps} \cdot i\;\;\; \mbox{ if $\mc{E}_{i-1}$ holds and $\phi(\bD)$  fits type $z|_{\bD}$},\\
& \Ext^{+}(z, \phi, i-1) \;\;\; \mbox{\hskip0.28\textwidth\relax otherwise}.
\end{cases}
\end{align*}

 We will show that $\Ext^{+}(z, \phi)$ is a supermartingale. Note that  $\Delta \Ext(z, \phi)$ is bounded above by a constant times the maximum possible value of $\Ext(z, \phi')$, where $\phi'$ is an extension of the partial embedding $\phi$ whose domain is $\bD':=\bD \cup \be_i$ (recall $e_i$ is the edge recieving a color in step $i$. Here $\be_i$ is meant to be a Platonic edge correspoinding to $e_i$). If $\be_i \subseteq \bD$ then $\Delta \Ext(z, \phi) $ is $\leq 0$ since no extension of $\phi$ can fit the type $z$ any longer. Thus we assume $\bD' \supsetneq \bD$ and so for any $\bD' \subseteq \bS' \subseteq \bS$ we have that $\bD \subsetneq \bS'$ and so  $\pow(\bS, \bS', z) \le 0$ and so in the good event we have
\[
 \Ext(z, \phi') \le n^{\maxpow(\bS, \bD', z) + 10 (m + s-d-1) \eps} = n^{10 (m+s-d - 1) \eps}
\]
and so 
\begin{equation}\label{eqn:crudeMax}
  \D \Ext(z, \phi)  \le n^{10 (\ell+s-d - 1) \eps}.
\end{equation}
Now we bound $\E[\Delta \Ext(z, \phi)|\mc{H}_i]$. We need only an upper bound and so we will consider only positive contributions to $\Ext(z, \phi)$ at step $i$. Such a positive contribution must arise from the newly colored edge $e_i$ and some extension of $\phi$ to a clique of type $z'$ where $z' = \pre(z, \be_i)$ is some predecessor of $z$ formed by uncoloring the Platonic edge $\be_i$ corresponding to $e_i$. 

If the color of $\be_i$ under $z$ is Platonic and no other colored edge has the same color, then the edge could be recolored with almost any color and still fit the type $z$. In that case, for any $\bS'$ with $\bD \subseteq \bS' \subseteq \bS$ we have $\pow(\bS, \bS', z') = \pow(\bS, \bS', z)$. Thus in the good event we have
\[
\Ext(z', \phi) \le n^{\maxpow(\bS, \bD, z')+ 10 (m + s-d-1) \eps} = n^{\maxpow(\bS, \bD, z)+ 10 (m + s-d-1) \eps}.
\]
Now, in order for any particular extension in $\Ext(z', \phi)$ to become an extension in $\Ext(z, \phi)$, the colored edge $e_i$ must be a particular edge. Thus the expected positive contribution to $\Ext(z, \phi)$ is at most  (recalling that $1-t \ge n^{-\eps}$)
\[
 n^{\maxpow(\bS, \bD, z)+ 10 (m + s-d-1) \eps} \cdot \frac{1}{(1-t) \binom{n}{2}} \le n^{\maxpow(\bS, \bD, z) -2 + [10 (m + s-d-1) +1]\eps}.
\]

The other case we need to consider is where the color of $\be_i$ under $z$ is a real color or a repeated Platonic color. In this case,  for any $\bS'$ with $\bD \subseteq \bS' \subseteq \bS$ we have $\pow(\bS, \bS', z') \le \pow(\bS, \bS', z)+ \f $. Thus in the good event we have
\[
\Ext(z', \phi) \le n^{\maxpow(\bS, \bD, z')+ 10 (m + s-d-1)\eps} \le n^{\maxpow(\bS, \bD, z)+ \f  + 10 (m + s-d-1)\eps}.
\]
Now, in order for any particular extension in $\Ext(z', \phi)$ to become an extension in $\Ext(z, \phi)$, the colored edge $e_i$ must be a particular edge and it must receive a particular color. Thus the expected positive contribution to $\Ext(z, \phi)$ is at most (recalling that $a(t)\geq n^{-\eps}$ and $1-t \ge n^{-\eps}$)
\[
n^{\maxpow(\bS, \bD, z)+ \f  + 10 (m + s-d-1)\eps} \cdot \frac{1}{(1-t) \binom{n}{2}|C|a} \le n^{\maxpow(\bS, \bD, z) -2 + [10 (m + s-d-1) + 3]\eps}.
\]
Thus in both cases we have 
\begin{equation}\label{eqn:crudeE}
    \E[\Delta \Ext(z, \phi)|\mc{H}_i] \le n^{\maxpow(\bS, \bD, z) -2 + [10 (m + s-d-1) + 3]\eps}.
\end{equation}

By line \eqref{eqn:crudeE} we have that $\E[\Delta \Ext^+(z, \phi)|\mc{H}_i] \le 0$, i.e.\ $\Ext^+(z, \phi)$ is a supermartingale. We also have 
\[
\D \Ext^+(z, \phi)  \le \D \Ext(z, \phi) \le n^{10 (m+s-d - 1) \eps}
\]
by line \eqref{eqn:crudeMax}. Thus for our application of Freedman's inequality we can use $D = n^{10 (m+s-d - 1) \eps}$ Also we have
\begin{align*}
   \Var[ \Ext^+(z, \phi, k)| \mc{H}_{k}] &= \Var[ \Ext(z, \phi, k)| \mc{H}_{k}] \\
   &\le \E[ \rbrac{\D \Ext(z, \phi, k)}^2| \mc{H}_{k}] \\
   &\le \E[ n^{10 (m+s-d - 1) \eps} \cdot |\D \Ext(z, \phi, k)| | \mc{H}_{k}] \\
   &\le n^{10 (m+s-d - 1) \eps} \cdot n^{\maxpow(\bS, \bD, z) -2 + [10 (m + s-d-1) + 3]\eps}\\
   & = n^{\maxpow(\bS, \bD, z) -2 + [20 (m + s-d-1) + 3]\eps}
\end{align*}

and so for all $i \le i_{max} < n^2$ we have
\[
V(i) = \sum_{0 \le k \le i} \Var[ \Ext^+(z, \phi, k)| \mc{H}_{k}] \le n^{\maxpow(\bS, \bD, z) + [20 (m + s-d-1) + 3]\eps},
\]
so for our application of Freedman's inequality we can use $b = n^{\maxpow(\bS, \bD, z) + [20 (m + s-d-1) + 3]\eps}$. We will choose $\lambda = n^{\frac 12 \maxpow(\bS, \bD, z) + [10 (m + s-d-1) + 2]\eps}$. Note that since $z$ colors at least one edge outside $\bD$, $\Ext^+(z, \phi, 0)=\Ext(z, \phi, 0)=0$. Then, Freedmann's inequality (Lemma \ref{lem:Freedman}) tells us that
\begin{align*}
  &\P\left[\exists i: \Ext^+(z, \phi, i) \geq \lambda \right]\\
  &\leq  \exp\left(-\frac{\lambda^2}{2(b+D\lambda) }\right)  \\
  & = \exp\left(-\frac{n^{\maxpow(\bS, \bD, z) + [20 (m + s-d-1) + 4]\eps}}{2\rbrac{n^{\maxpow(\bS, \bD, z) + [20 (m + s-d-1) + 3]\eps}+n^{\frac 12 \maxpow(\bS, \bD, z) + [20 (m + s-d-1) + 2]\eps}} } \right) \\
  &= \exp \rbrac{-\Omega\rbrac{n^{\eps}}}.
\end{align*}
But if the good event holds and $\Ext^+(z, \phi, i) < \lambda$ then 
\begin{align*}
   \Ext(z, \phi, i) &< n^{\maxpow(\bS, \bD, z)-2 + 10 (m + s-d) \eps} \cdot i + n^{\frac 12 \maxpow(\bS, \bD, z) + [10 (m + s-d-1) + 2]\eps}\\
   &< n^{\maxpow(\bS, \bD, z) + 10 (m + s-d) \eps}.
\end{align*}

\end{proof}

\section{Dynamic concentration of \texorpdfstring{$\Ava(e)$}{2}}\label{sec:A}

In this section we bound the probability that $\mc{E}_{i_{max}}$ fails due to Condition \eqref{cond:atraj}. We define a family of random variables that we will show are supermartingales, to which we will then apply Freedman's inequality. For each edge $e$ we define variables
\[
\Ava^\pm(e)=\Ava^\pm(e, i):=\begin{cases} 
& \Ava(e, i) -  \ava(t) \mp |C|f_\Ava(t) \;\;\; \mbox{ if $\mc{E}_{i-1}$ holds and $e \in U_i$},\\
& \Ava^\pm(e,i-1) \;\;\; \mbox{\hskip0.28\textwidth\relax otherwise}.
\end{cases}
\]
Note that in the above definition we are using ``$\pm$'' to simultaneously define $\Ava^+(e)$ and $\Ava^-(e)$. We will establish that the variables $\Ava^+(e)$ are supermartingales. We claim that $\Ava^-(e)$ are submartingales (or equivalently, $-\Ava^-(e)$ are supermartingales), but we will not separately justify it since it is very similar to showing that the $\Ava^+(e)$ are supermartingales.

\subsection{Establishing that \texorpdfstring{$\Ava^+(e)$}{2} is a supermartingale}

We need to show that $ \E[\D \Ava^+(e)| \mc{H}_i]\leq 0$. First we will estimate $ \E[\D \Ava(e)| \mc{H}_i]$. Let $\mc{F}(e, c, i)$ be the event that $c$ is forbidden at $e$ at step $i$ (this event excludes the possibility that $e_i=e$). Then we have 
\begin{equation}\label{eqn:EDA1}
   \E[\D \Ava(e)| \mc{H}_i] = -\sum_{c \in \Ava(e)} \P[\mc{F}(e,c, i)|\mc{H}_i].
\end{equation}
and
\begin{equation*}
   \P[\mc{F}(e,c, i)|\mc{H}_i] =  \frac{(1-t)\binom n2-1}{(1-t)\binom n2} \cdot \P[\mc{F}(e,c, i)|\mbox{$\mc{H}_i$ and $e_i \neq e$}]
\end{equation*}

Fix a color $c \in \Ava(e)$. Let $F(e, c)$ be the number of pairs $(e', c')$ such that $e' \neq e$ is uncolored, $c'$ is available at $e'$, and assigning $c'$ to $e'$ would forbid $c$ at $e$.
Then  
 \begin{align}
     \P[\mc{F}(e,c, i)|\mbox{$\mc{H}_i$ and $e_i \neq e$}] &= \sum_{e' \in U_i \sm \{e\}} \frac{1}{(1-t)\binom n2-1} \cdot \P[\mc{F}(e,c , i) \; | \; e_i = e'\mbox{ and } \mc{H}_i] \nonumber \\
     & =\frac{1}{(1-t)\binom n2-1 } \sum_{e' \in U_i \sm \{e\}}  \frac{\Big| \{ c': (e', c') \in F(e, c)  \} \Big|}{\Ava(e')} \nonumber\\
     & \ge \frac{\Big| F(e, c)   \Big|}{\sbrac{(1-t)\binom n2-1} \cdot (\ava + |C|f_\Ava)}  \nn
 \end{align}
 and so 
 \begin{align}\label{eqn:prF}
     \P[\mc{F}(e,c, i)|\mc{H}_i] &\ge \frac{\Big| F(e, c)   \Big|}{(1-t)\binom n2 \cdot (\ava + |C|f_\Ava)} \nn\\
     & \ge \rbrac{1 - \frac{f_\Ava}{a}} \frac{\Big| F(e, c)   \Big|}{(1-t)\binom n2 \cdot \ava }
 \end{align}
where on the last line we used that  
\begin{equation*}
    \frac{1}{1 + \frac{f_\Ava}{a}} \ge 1  - \frac{f_\Ava}{a} .
\end{equation*}

Note that for each pair $(e',c') \in F(e, c)$ there exists at least one $S\subseteq V$ such that coloring $e'$ the color $c'$ at step $i$ would cause $c$ to be forbidden at $e$ through $S$. We define the following notation. 

\begin{defn}\label{def:Sec}
 We let $S_{e, c}(e', c') = S_{e, c}(e', c', i)$ be the family of minimal sets $S\subseteq V$ such that assigning $c'$ to $e'$ at step $i$ would forbid $c$ at $e$ through $S$.
\end{defn} 

We would like to claim (with some explicitly bounded error term) that
\[
|F(e, c)|  =  \Big| \Big\{ (e', c'): (e',c') \in F(e, c)\Big\}  \Big| \approx \Big| \Big\{ (S, e', c'):  (e',c') \in F(e, c), S \in S_{e, c}(e', c')\Big\} \Big|.
\]
In other words, we claim that for almost all pairs $(e', c')\in F(e, c)$,  $S_{e, c}(e', c')$ is just a singleton and that there are not too many pairs $(e', c')$ where $S_{e, c}(e', c')$ is too large.
We have the precise identity
\begin{equation}\label{eqn:Fbound}
     \Big| \Big\{ (S, e', c'):  (e',c') \in F(e, c), S \in S_{e, c}(e', c')\Big\} \Big| = |F(e, c)| + \sum_{\substack{(e',c') \in F(e, c) \\ |S_{e, c}(e', c')| \ge 2 }} \Big( |S_{e, c}(e', c')| -1 \Big)
\end{equation}
so we try to bound the sum on the right. First we bound the number of terms, which we denote by $B_1$.  Fix  an edge $e$ and a real color $c$. We want to bound the number of pairs $(e', c')$ such that there exist distinct sets $S_1, S_2 \in S_{e, c}(e', c')$. The number of such pairs $(e', c')$ is at most a constant times the number of possible sets $S = S_1 \cup S_2$ so we will bound the latter instead.
     
We will use Lemma \ref{lem:appendpreforb}. For $j=1, 2$ say $S_j$ fits the type $y_j$ on the set of Platonic vertices $\bS_j$, where $y_j$ is a $(\be, c, \be', c'')$-preforbidder (here $\be$ corresponds to $e$, $\be'$ corresponds to $e'$, and $c''$ corresponds to $c'$ (meaning that $c''=c$ if $c'=c$  or else $c''$ is Platonic if $c'\neq c$)). Assume we choose $y_1, y_2$ to be compatible, and so $S=S_1 \cup S_2$ fits the type $y_1 \cup y_2$.  To bound the number of possible sets $S = S_1 \cup S_2$, we bound $\Ext(y_1 \cup y_2, \phi)$ where $\phi$ is of order 2 (mapping $\be$ to $e$). Applying Lemma \ref{lem:appendpreforb}  with $\be_1=\be_2 = \be,\; \be_3= \be',\; c_2=c,\; c_3=c''$ we have
\[
\maxpow(\bS_1 \cup \bS_2, \be, y_1 \cup y_2) \le \pow(\bS_1, \be, y_1) - \m \le \frac{p-3}{\binom p2 - q + 1},
\]
where the last inequality follows from Observation \ref{obs:preforbpow}.
Thus by Corollary \ref{cor:crudepractical} we have in the good event that
\[
|\Ext(y_1 \cup y_2, \phi)|  = O\rbrac{n^{\frac{p-3}{\binom p2 - q + 1} + 20p^2 \eps}}.
\]
Of course $|\Ext(y_1 \cup y_2, \phi)|$ only counts the contribution to $B_1$ due to some fixed $y_1, y_2$. But since there are only a constant number of choices for the colored portions $\colport(y_1), \colport(y_2)$  there are only a constant number of relevant types $y_1 \cup y_2$. Thus the number of terms in the sum on line \eqref{eqn:Fbound} is
\[
B_1 = \Big| \Big\{ (e',c') \in F(e, c): |S_{e, c}(e', c')| \ge 2 \Big\} \Big|= O\rbrac{n^{\frac{p-3}{\binom p2 - q + 1} + 20p^2 \eps}}.
\]
Now we will bound the maximum possible size of any term in the sum on line \eqref{eqn:Fbound}. Suppose $e, e'$ are fixed real edges and $c, c'$ are fixed real colors. We bound $|S_{e, c}(e', c')|$ by the sum of a constant number of terms of the form $|\Ext(y, \phi)|$ where $\phi$ has domain $\be \cup \be'$ sending $\be$ to $e$ and $\be'$ to $e'$, and $y$ is a $(\be, c, \be', c')$-preforbidder containing no real colors outside of $\{c,c'\}$.
Consider any $\bS'$ with $\be \cup \be' \subseteq \bS' \subseteq \bS$. If $\bS' \neq \bS$ then $\pow(\bS, \bS', y)<0$ by Lemma~\ref{lem:preforbpowbound}. On the other hand if $\bS' = \bS$ then of course $\pow(\bS, \bS', y)=0$. Altogether we have $\maxpow(\bS, \be \cup \be', y)=0$ and so by Corollary \ref{cor:crudepractical} we have
\[
|\Ext(y, \phi)| \le n^{ 20p^2 \eps},
\]
and so the maximum possible size of any term in the sum on line \eqref{eqn:Fbound} is at most $B_2$ where
\[
B_2 = O\rbrac{n^{ 20p^2 \eps}}.
\]

 We return to our calculation of $\P[\mc{F}(e,c, i)|\mc{H}_i]$ on \eqref{eqn:prF}. Since the sum on line \eqref{eqn:Fbound} is at most $B_1 B_2 = O\rbrac{n^{\frac{p-3}{\binom p2 - q + 1} + 40p^2 \eps}}$ we have
  \begin{equation}\label{eqn:triples}
    |F(e, c)|  = \Big| \Big\{ (S, e', c'):  (e',c') \in F(e, c), S \in S_{e, c}(e', c')\Big\} \Big| +  O\rbrac{n^{\frac{p-3}{\binom p2 - q + 1} + 40p^2\eps}}. 
 \end{equation}
Now for each triple $(S, e', c')$ counted above we have that 
before we colored $e'$, there was an embedding $\phi$ with image $S$ fitting a type $z$ such that $(z,\be',c')\in\prefbd(\be, c)$, where $\phi(\be')=e'$. But each such $S$ is counted by at most $\equiv(z, \be') \cdot (\ava(t) + |C|f_\Ava(t))^{\binom{s}2 - \ell -2}$ types. So \eqref{eqn:triples} is (using $\ell=\ell(z)$, $r=r(z)$,  $k=k(z)$, $s=|\bS(z)|$) 

\begin{align}
 & \sum_{(z, \be', c^*) \in \prefbd(\be, c)  }  \frac{|\Ext(z, \phi)|}{\equiv(z, \be')(\ava + |C|f_\Ava)^{\binom{s}{2}-\ell-2}} +  O\rbrac{n^{\frac{p-3}{\binom p2 - q + 1} + 40p^2 \eps}}\nn\\
    \ge  &  \sum_{(z, \be', c^*) \in \prefbd(\be, c)  } \frac{\ext_z- n^{s-2}|C|^{-r-k}f_z}{\equiv(z, \be')(\ava + |C|f_\Ava)^{\binom{s}{2}-\ell-2}} +  O\rbrac{n^{\frac{p-3}{\binom p2 - q + 1} + 40p^2 \eps}}\label{eqn:triples2}
\end{align}
We estimate the sum in \eqref{eqn:triples2} by splitting the numerator. For the first term in the numerator we write
\begin{align}
&  \sum_{(z, \be', c^*) \in \prefbd(\be, c)  } \frac{\ext_z}{\equiv(z, \be')(\ava + |C|f_\Ava)^{\binom{s}{2}-\ell-2}}\nn\\
& \ge \rbrac{1 - p^2\frac{f_\Ava}{a}} \sum_{(z, \be', c^*) \in \prefbd(\be, c)  } \frac{\ext_z}{\equiv(z, \be')\ava^{\binom{s}{2}-\ell-2}}\nn\\
& = \rbrac{1 - p^2\frac{f_\Ava}{a}}|C| a(t) (1-t) h(t) + \TO\rbrac{ n^{-\mf } }\nn\\
& = |C| a(t) (1-t) h(t) - p^2 |C| (1-t)h(t)f_\Ava + \TO\rbrac{ n^{-\mf } }.\label{eqn:triples3}
\end{align}
Indeed, on the second line we used that since $f_\Ava/a = o(1)$ for all $t \le t_{max}$, we have that for any $\ell \ge 0$ 
\begin{equation*}
    \frac{1}{\rbrac{1 + \frac{f_\Ava}{a}}^{\binom{s}2 - 2-\ell}} \ge \frac{1 }{\rbrac{1 + \frac{f_\Ava}{a}}^{\binom{p}2 - 2}} \ge 1  - p^2\frac{f_\Ava}{a} .
\end{equation*} 
On the third line we have replaced the expression from \eqref{eqn:sumpreforb} with \eqref{eqn:h}. 
Now for the second term in the numerator of \eqref{eqn:triples2} we note that all terms have $r+k = R(s) \ge \frac{(s-2)\rbrac{\binom{p}{2}-q+1}}{p-2} -1$ so we have
\begin{align}
n^{s-2}|C|^{-r-k}f_z & \le n^{s-2}|C|^{-\frac{(s-2)\rbrac{\binom{p}{2}-q+1}}{p-2} +1}f_z \nn\\
&= n^{s-2} |C| \sbrac{\k n^{\f} \log^{-\m} n}^{-\frac{(s-2)\rbrac{\binom{p}{2}-q+1}}{p-2} }  \cdot a(t)^{\binom s2 - 1 - \ell}f_\Ext \nn\\
& \le \k^{-\frac{(s-2)\rbrac{\binom{p}{2}-q+1}}{p-2}} |C| \log n \cdot a(t)^{\binom s2 - 1 - \ell}f_\Ext\nn
\end{align}
and so
\begin{align}
  &\sum_{(z, \be', c^*) \in \prefbd(\be, c)  } \frac{ n^{s-2}|C|^{-r-k} f_z}{\equiv(z, \be')(\ava + |C|f_\Ava)^{\binom{s}{2}-\ell-2}}\nn\\
  & \le  \sum_{(z, \be', c^*) \in \prefbd(\be, c)  } \frac{ \k^{-\frac{(s-2)\rbrac{\binom{p}{2}-q+1}}{p-2}} |C| \log n \cdot a(t)^{\binom s2 - 1 - \ell}f_\Ext}{\equiv(z, \be')\ava ^{\binom{s}{2}-\ell-2}}\nn\\
  & =  \sum_{\substack{ 3 \le s \le p \\ 1\leq \ell\leq \binom{s}2 -2\\ 1\leq \eta\leq s! |\bC|!}} \mu(\ell, \eta, s) |C|^{\binom{s}{2}-\ell-2} \cdot \frac{ \k^{-\frac{(s-2)\rbrac{\binom{p}{2}-q+1}}{p-2}} |C| \log n \cdot a(t)^{\binom s2 - 1 - \ell}f_\Ext}{\eta \ava ^{\binom{s}{2}-\ell-2}}\nn\\
  & \le |C| \log n \cdot a(t) f_\Ext \cdot \k^{-\frac{\binom{p}{2}-q+1}{p-2}} \sum_{\substack{ 3 \le s \le p \\ 1\leq \ell\leq \binom{s}2 -2\\ 1\leq \eta\leq s! |\bC|!}}   \frac{ \mu(\ell, \eta, s)}{\eta }\nn\\
  & \le \eps|C| \log n \cdot a(t) f_\Ext.\label{eqn:tripleserror}
\end{align}
where we have used \eqref{eqn:pickkappa}.
Putting \eqref{eqn:triples3} and \eqref{eqn:tripleserror} together, line \eqref{eqn:triples2} is at least 
\[
|C| a(t) (1-t)h(t) - p^2|C|(1-t)h(t)f_\Ava- \eps |C|\log n \cdot a(t) f_\Ext  + O\rbrac{n^{\frac{p-3}{\binom p2 - q + 1} + 40p^2 \eps}}.
\]

Returning to line \eqref{eqn:prF}, we see that 
\begin{align}
    &\P[\mc{F}(e,c, i)|\mc{H}_i] \nonumber\\
    & \ge \frac{1 - \frac{f_\Ava}{a} }{(1-t)\binom n2 \cdot \ava } \sbrac{|C| a(t) (1-t)h(t) - p^2|C|(1-t)h(t)f_\Ava- \eps |C|\log n \cdot a(t) f_\Ext  + O\rbrac{n^{\frac{p-3}{\binom p2 - q + 1} + 40p^2 \eps}}} \nonumber\\
    & \ge \frac{1 - \frac{f_\Ava}{a} }{\binom n2} \sbrac{h(t) - p^2 \eps \log n \cdot \frac{f_\Ava}{a}- \eps \log n \cdot \frac{f_\Ext}{1-t}  + O\rbrac{n^{-\m + 50p^2 \eps}}} \nonumber\\
    & \ge \frac{1}{\binom n2 } \sbrac{ h(t) - 2p^2 \eps \log n \cdot \frac{f_\Ava}{a}- \eps \log n \cdot \frac{f_\Ext}{1-t}} + O\rbrac{n^{-2-\m + 50p^2\eps}}.\label{equation for PF(e,c)}
\end{align}
Similarly,
\begin{align}\label{eqn:Fupper}
    \P[\mc{F}(e,c, i)|\mc{H}_i] &\le  \frac{1}{\binom n2 } \sbrac{ h(t) + 2p^2 \eps \log n \cdot \frac{f_\Ava}{a}+ \eps \log n \cdot \frac{f_\Ext}{1-t}} + O\rbrac{n^{-2-\m + 50p^2\eps}}.\notag \\
    &\leq \frac{2\epsilon \log n}{\binom{n}2},
\end{align}
which we will need later. Note that for the last line we have used the fact that $h(t) \le \eps \log n$, $f_\Ava(t)/a(t) = o(1)$ and $f_\Ext(t)/(1-t)=o(1)$. 

Returning to line \eqref{eqn:EDA1} we see that
\begin{align}
   & \E[\D \Ava(e)| \mc{H}_i]= -\sum_{c \in \Ava(e)} \P[\mc{F}(e, c, i)|\mc{H}_i]\nn\\
    & \le -\Big(\ava(t)-|C|f_\Ava(t)\Big) \cbrac{ \frac{1}{\binom n2 } \sbrac{  h(t) - 2p^2 \eps \log n \cdot \frac{f_\Ava}{a}- \eps \log n \cdot \frac{f_\Ext}{1-t}} + O\rbrac{n^{-2-\m + 50p^2\eps}}} \nn\\
    & \le -\frac{|C|}{\binom n2 }   \sbrac{ a(t)h(t) - 3p^2\eps \log n \cdot f_\Ava- \eps \log n \cdot \frac{a(t)f_\Ext}{1-t}} + O\rbrac{n^{-2 + \frac{p-3}{\binom p2 - q + 1}+ 50p^2\eps}} \label{eqn:EDava}
\end{align}

Recall that by Taylor's theorem, we have that for any twice differentiable function $g(t)$ that for some $\t \in [t, t+ \D t]$ that 
\[
g(t+\D t)- g(t) =  g'(t)\D t + \frac 12 g''(\t)(\D t)^2.
\]
We will apply the above to the function $\ava(t) + |C|f_\Ava(t) = |C|(a(t) + f_\Ava(t))$ with $\D t = 1/\binom n2$ (the change in $t$ from step $i$ to $i+1$). Note that 
\[
a'(t) = -a(t) h(t), \qquad a''(t) = a(t) h(t)^2 - a(t) h'(t) 
\]
so $a''(t) = O(\log^2 n)$ for all $0 \le t \le t_{\max}$. Using \eqref{eqn:pickfava}  we have
\begin{align}
    f_\Ava'(t) &= n^{-\hm + 20p^2\eps} \rbrac{5p^2(1-t)^{-5p^2-1} e^{10p^4 \eps t \log n} + 10p^4 \eps \log n \cdot (1-t)^{-5p^2} e^{10p^4 \eps t \log n}}\nn\\
    & \ge n^{-\hm + 20p^2\eps} \cdot 10p^4 \eps \log n \cdot (1-t)^{-5p^2} e^{10p^4 \eps t \log n}\label{eqn:favaprime}
\end{align}
    Also, for all $1 \le t \le t_{max}$ we have
   \begin{align*}
   f_\Ava''(t) &= n^{-\hm + 20p^2\eps} \Bigg[5p^2(5p^2+1)(1-t)^{-5p^2-2} e^{10p^4 \eps t \log n} +100p^6 \eps \log n (1-t)^{-5p^2-1} e^{10p^4 \eps t \log n}\nn\\
   & \qquad \qquad \qquad \qquad + 100p^8 \eps^2 \log^2 n(1-t)^{-5p^2} e^{10p^4 \eps t \log n}\Bigg]\nn\\
   & \le n^{-\hm + 20p^2\eps}O \rbrac{n^{(10p^4+5p^2+2 )\eps} \log^2 n } \nn\\
   &=o(1),
\end{align*}
where on the last line we have used that since $\eps = 10^{-3}p^{-6}$ and $\hm \ge \frac12 p^{-2}$ we have that the power of $n$ is
\begin{equation*}
    -\hm + 20p^2 \eps  +(10p^4+5p^2+2 )\eps \le \frac{-490p^4+25p^2+2 }{2000 p^6} <0.
\end{equation*}
Thus Taylor's theorem gives us that 
\[
\D\sbrac{|C|(a(t) + f_\Ava(t))} = \frac{|C|}{\binom n2}\sbrac{-a(t)h(t) + f'_\Ava(t)} + O\rbrac{\frac{|C|\log^2 n}{n^4}}.
\]
Using the above and \eqref{eqn:EDava} we get 
\begin{align}
    \E[\D &\Ava(e)^+| \mc{H}_i]\nn\\
    &\le -\frac{|C|}{\binom n2 } \sbrac{ a(t)h(t) -3p^2 \eps \log n \cdot f_\Ava - \eps \log n \cdot \frac{a(t) f_\Ext}{1-t}} - \frac{|C|}{\binom n2}\sbrac{-a(t)h(t) + f'_\Ava(t)}\nn\\
    & \qquad + O\rbrac{n^{-2 +\mf + 50p^2\eps}}\nn\\
   & = \frac{|C|}{\binom n2}\sbrac{ 3p^2 \eps \log n \cdot f_\Ava + \eps \log n \cdot \frac{a(t) }{1-t} f_\Ext - f_\Ava'(t)} + O\rbrac{n^{-2+ \mf + 50p^2\eps}}.\label{eqn:EDavaplus}
\end{align}
Finally, we verify that we have chosen $f_\Ava, f_\Ext$ that make \eqref{eqn:EDavaplus} negative. Indeed, using \eqref{eqn:pickfava}, \eqref{eqn:pickfext} and \eqref{eqn:favaprime} we have for all $0 \le t \le t_{max}$,
\begin{align*}
  & 3p^2 \eps \log n \cdot f_\Ava + \eps \log n \cdot \frac{a(t) }{1-t} f_\Ext - f_\Ava'(t) \\
  &\le n^{-\hm + 20p^2\eps}  \Bigg( 3p^2 \eps \log n \cdot  (1-t)^{-5p^2} e^{10p^4 \eps t \log n} + \eps \log n \cdot (1-t)^{-5p^2} e^{10p^4 \eps t \log n}\\
  & \qquad \qquad - 10p^4 \eps \log n (1-t)^{-5p^2 }e^{10p^4 \eps t \log n}  \Bigg)\\
  &= \rbrac{ - 10p^4 + 3p^2   + 1    }\eps n^{-\hm + 20p^2\eps} \log n \cdot  (1-t)^{-5p^2} e^{10p^4 \eps t \log n}\\
  & = -\TOM \rbrac{n^{-\hm + (10p^4+25p^2) \eps }}.
\end{align*}
Therefore \eqref{eqn:EDavaplus} is at most 
\begin{align*}
    &-\TOM \rbrac{n^{-2+\frac{p-5/2}{\binom p2 - q + 1} + (10p^4+25p^2) \eps }}+ O\rbrac{n^{-2+ \mf + 50p^2\eps}} < 0.
\end{align*}
The last inequality follows because the first power of $n$ is larger than the second. Indeed, if we subtract the second power from the first, and use the fact that $\binom p2 - q + 1 \le p^2$, we have
\[
\hm + (10p^4-25p^2)\eps \geq \frac 1{2p^2} + \frac{10p^4-25p^2}{1000p^6} = \frac{102p^2-5}{200p^4} >0.
\]

\subsection{Applying Freedman's inequality to \texorpdfstring{$\Ava^+(e)$}{2}}

First we bound $|\D \Ava(e)|$. Suppose at step $i$ that the edge $e_i$ receives the color $c_i$ causing some color to be forbidden at $e$ through $S$. Then $S$ fits some colors-only type $z$ that is a $(\be, \bc, \be', c_i)$-preforbidder (where $\be, \be'$ correspond to $e, e_i$, and $\bc$ corresponds to the color being forbidden). Then by Lemma~\ref{lem:preforbpowbound} we have that $\maxpow(\bS(z), \be \cup \be', z) \le \frac{p-3}{\binom p2 - q + 1}$ and so by Corollary \ref{cor:crudepractical}, if $\phi$ is the partial embedding on $\be\cup\be'$ mapping $\be$ to $e$ and $\be'$ to $e_i$, we have $\Ext(z, \phi) = O\rbrac{n^{\frac{p-3}{\binom p2 - q + 1}+20p^2 \eps}}.$ Since each extension in $\Ext(z, \phi)$ gives a constant number of colors forbidden at $e$, and since there are a constant number of choices for the preforbidder $z$, we have $|\D \Ava(e)|= O\rbrac{n^{\frac{p-3}{\binom p2 - q + 1}+20p^2 \eps}}$.
Also the one-step change in the deterministic part of $\Ava^+(e)$ is $|\D (\ava + |C|f_\Ava)| = o(1)$ so we have $|\D \Ava(e)^+| \le O\rbrac{n^{\frac{p-3}{\binom p2 - q + 1}+20p^2 \eps}}$. Thus  in our application of Lemma~\ref{lem:Freedman}, we will use $D=O\rbrac{n^{\frac{p-3}{\binom p2 - q + 1}+20p^2 \eps}}$. Also we have
\begin{align*}
  \Var[ \Delta \Ava^+(e,k)| \mc{H}_{k}] = \Var[ \Delta \Ava(e,k)| \mc{H}_{k}] 
  &\le \E[ \rbrac{\Delta \Ava(e,k)}^2| \mc{H}_{k}]\\
 &\le \E\sbrac{ n^{\mf+ 20p^2 \eps} \cdot |\Delta \Ava(e,k)|| \mc{H}_{k}}\\
 &\le n^{\mf+20p^2 \eps} \cdot \frac{2|C|\epsilon \log n}{\binom n2}\le \frac{n^{\frac{2p-5}{\binom p2 - q + 1} + 30p^2 \eps}}{\binom n2}
\end{align*}
where we have used lines \eqref{eqn:EDA1} and \eqref{eqn:Fupper}.
and so for our application of Freedman we have for all $i \le i_{max} < \binom n2$ that 
\[
V(i) = \sum_{0 \le k \le i} \Var[ \D \Ava^+(e,k)| \mc{H}_{k}] \le n^{\frac{2p-5}{\binom p2 - q + 1} + 30p^2 \eps},
\]
 so we will use $b = n^{\frac{2p-5}{\binom p2 - q + 1} + 30p^2 \eps}.$ Note that at the beginning of the process we have $\Ava^+(e,0) = \Ava(e,0) - \ava(0) - |C|f_\Ava(0) = -|C|f_\Ava(0).$ For $\Ava^+(e,i)$ to become positive would therefore be a positive change of $\lambda := |C|f_\Ava(0)=\TT\rbrac{n^{\frac{p-5/2}{\binom p2 - q + 1} + 20p^2 \eps}}$. Freedman's inequality gives us a failure probability of at most
\[
\exp\left(-\frac{\lambda^2}{2(b+D\lambda) }\right) \le \exp\rbrac{-n^{\Omega(\eps)}}
\]
which beats any polynomial union bound.

\section{Dynamic concentration of  \texorpdfstring{$\Ext(z, \phi)$}{2}}\label{sec:G}

In this section we bound the probability that $\mc{E}_{i_{max}}$ fails due to Condition \ref{cond:exttraj}. Let $(z, \be)$ be a rooted trackable type on $\bS = \bS(z)$, let $s:=|\bS|$ and let $\phi$ have domain $\be$ (so $\phi$ has order 2), where $\phi(\be)=e$. Let 
\[
\Ext^{\pm}(z, \phi)=\Ext^{\pm}(z, \phi, i):=\begin{cases} 
& \Ext(z, \phi) -  \ext_z \mp n^{s-2}|C|^{-r-k} f_z \;\;\; \mbox{ if $\mc{E}_{i-1}$ holds and $e \in U_i$},\\
& \Ext^{\pm}(z, \phi, i-1) \;\;\; \mbox{\hskip0.2\textwidth\relax otherwise}.
\end{cases}
\]
We will show that $\Ext^{+}(z, \phi)$ is a supermartingale, i.e.\ that $ \E[\D \Ext^+(z, \phi)| \mc{H}_i]\le 0$. Thus we will need to estimate  $ \E[\D \Ext(z, \phi)| \mc{H}_i]$. To apply Freedman's inequality we also need to bound $|\D \Ext(z, \phi)|$ in the good event.  We have
\[
\D \Ext(z, \phi) = D_1 - D_2 - D_3,
\]
where $D_1$ is the number of new extensions $\phi'$ that come into $\Ext(z,  \phi)$, $D_2$ is the number of extensions $\phi'$ that leave $\Ext(z, \phi)$ due to edges in $\phi'(\bS)$ getting colored, and $D_3$ is the number of extensions $\phi'$ that leave $\Ext(z, \phi)$ due to colors being forbidden on edges in $\phi'(\bS)$ (except for the extensions already counted by $D_2)$. We will handle $D_1, D_2, D_3$ separately, finding the expected change and maximum possible change for each one. 
Estimating $ \E[\D \Ext^+(z, \phi)| \mc{H}_i]$ will require us to estimate the one-step change in the deterministic function $n^{s-2}|C|^{-r-k}\ext_z(t)$. Recall from \eqref{eqn:exttraj} that 
\[
\ext_z(t) = n^{s-2}|C|^{-r-k}t^\ell (1-t)^{\binom{s}{2}-\ell-1}a(t) ^{\binom{s}{2}-\ell-1} 
\]
and its derivative with respect to $t$ is (using the product rule and $a'=-ah$)
\[
\ext_z'(t) = n^{s-2}|C|^{-r-k} ( g_1(t) - g_2(t) - g_3(t))
\]
where
\begin{align}
    & g_1(t) := \ell t^{\ell-1} (1-t)^{\binom{s}{2}-\ell-1} a(t) ^{\binom{s}{2}-\ell-1} \label{def:g1}\\
   &g_2(t):=  \rbrac{\binom{s}{2}-\ell-1} t^\ell (1-t)^{\binom{s}{2}-\ell-2} a(t) ^{\binom{s}{2}-\ell-1} \label{def:g2}\\
   &g_3(t):=\rbrac{\binom{s}{2}-\ell-1} t^\ell (1-t)^{\binom{s}{2}-\ell-1} a(t) ^{\binom{s}{2}-\ell-1}h(t).\label{def:g3}
\end{align}
We discuss a little motivation for the coming calculations. We will show that the expected one-step change in $\Ext(z, \phi)$ is approximately equal to the one-step change in $\ext_z(t)$. Note that 
\[
  \E[\D \Ext(z, \phi)| \mc{H}_i] =   \E[D_1| \mc{H}_i]-\E[D_2| \mc{H}_i]-\E[D_3| \mc{H}_i],
\]
and the three terms above naturally correspond to the three terms in $\D \ext_z(t)$. In particular, we will show that for $1\leq j\leq 3$,
\begin{equation}\label{eqn:Dj}
\E[D_j\mid \mc{H}_i]\approx \frac{n^{s-2}|C|^{-r-k} }{\binom{n}2} g_j(t),
\end{equation}
where $\approx$ will be made rigorous. Note that we have
\begin{equation}\label{eqn:g estimate}
g_1(t), g_2(t) \le p^2 a(t)^{\binom{s}{2}-\ell-1}, \quad g_3(t) \le p^2  a(t)^{\binom{s}{2}-\ell-1} h(t)
\end{equation}

We will now use Taylor's theorem to estimate the one-step change in $\ext_z(t)$. First let us give a crude bound on the second derivative of $\ext_z(t)$. Recall from \eqref{eqn:hdef} that $h(t)$ is a polynomial in $t$, whose coefficients only depend on $n$ in that they have factors of the form $\log^\frac{s-2}{p-2} n$ for $s \le p$. Then,
\[
\frac{d}{dt}\sbrac{\frac{\ext_z(t)}{n^{s-2}|C|^{-r-k}}}=\frac{d}{dt} \sbrac{t^\ell (1-t)^{\binom{s}{2}-\ell-1}a(t) ^{\binom{s}{2}-\ell-1}} = P_1(t) a(t) ^{\binom{s}{2}-\ell-1}
\]
where $P_1(t)$ is a polynomial in $t$ with coefficients that may have factors $\log^\frac{s-2}{p-2} n \le \log n$. Taking another derivative with respect to $t$, we get
\[
\frac{d^2}{dt^2}\sbrac{\frac{\ext_z(t)}{n^{s-2}|C|^{-r-k}}}= P_2(t) a(t) ^{\binom{s}{2}-\ell-1},
\]
where $P_2$ is again a polynomial in $t$, and the largest power of $\log n$ that shows up as a coefficient is $\log^2 n$. Therefore, since $P_2/\log^2 n$ and $a$ are bounded when $0\leq t\leq 1$,
\[
\ext_z''(t) = \TO\rbrac{n^{s-2}|C|^{-r-k}}.
\]
Taylor's theorem then gives us that from step $i$ to $i+1$, the change in $\ext_z(t)$ is 
\begin{align}
   \D [\ext_z(t)] &=  \ext'(t) \frac{1}{\binom n2} + \TO\rbrac{\frac{n^{s-2}|C|^{-r-k} }{\binom n2 ^2}}\nn\\
   &= \frac{n^{s-2}|C|^{-r-k}}{\binom n2} ( g_1(t) - g_2(t) - g_3(t))+ \TO\rbrac{n^{\pow(\bS, \be, z) - 4 }} \label{eqn:Dext}
\end{align}
Now we will similarly use Taylor's theorem to estimate $\D[f_z(t)]$. Recall from \eqref{eqn:pickfext} and \eqref{eqn:pickfz} that 
\[
f_z(t) = n^{-\hm + 20p^2\eps} (1-t)^{-5p^2+1} a(t)^{\binom s2 - 2 - \ell}e^{10p^4 \eps t \log n}.
\]
Then we can see that 
\[
f_z''(t) = n^{-\hm + 20p^2\eps} (1-t)^{-5p^2-1} Q(t) a(t)^{\binom s2 - 2 - \ell}e^{10p^4 \eps t \log n} 
\]
where $Q(t)$ is some polynomial in $t$ whose only dependence on $n$ is in the form of powers of $\log n$ up to $\log^2 n$. Thus for all $t \le t_{\max} = 1-n^{-\eps}$ we have 
\[
f_z''(t) = \TO\rbrac{n^{-\hm + 20p^2\eps + (5p^2+1)\eps  + 10p^4 \eps} } = o(1).
\]
Taylor's theorem then gives us that  the change in $f_z(t)$ is 
\begin{align}
   \D [f_z(t)] &=  f_z'(t) \frac{1}{\binom n2} + o\rbrac{\frac{1}{\binom n2 ^2}}= \frac{1}{\binom n2} f_z'(t)  + o\rbrac{\frac{1}{n^4}} \label{eqn:Dfz}
\end{align}
Now putting \eqref{eqn:Dext} and \eqref{eqn:Dfz} together we see that 
\begin{align*}
    \D \sbrac{ \ext_z + n^{s-2}|C|^{-r-k} f_z} & = \frac{n^{s-2}|C|^{-r-k}}{\binom n2} \sbrac{ g_1(t) - g_2(t) - g_3(t)+ f_z'(t)} + \TO\rbrac{n^{\pow(\bS, \be, z) - 4 }}
\end{align*}

Having just estimated the one-step change in the deterministic terms in $\Ext^+(z, \phi)$, in the coming subsections we will estimate the expectations of $D_1, D_2, D_3$, formalizing \eqref{eqn:Dj}. We will also find absolute bounds on those variables which hold in the good event. 

\subsection{\texorpdfstring{$D_1$}{1}}

\subsubsection{Estimating  \texorpdfstring{$ \E[D_1| \mc{H}_i]$}{1}}

For an extension $\phi'$ counted by $D_1$, in the previous step of the process the image $\phi'(\bS)$ must have fit some type of the form $\pre(z, \be')$. Each such extension $\phi'$ counts exactly one way to create a new extension of type $z$ in this step (the edge $\phi'(\be')$  corresponding to the edge $\be'$ in $z$ needs a certain color). Let $E_C$ be the set of edges that are colored under $z$. 

We estimate $\E[D_1|\mc{H}_i]$ as follows.
For each $\be'\in E_C$, each $\phi'$ that is of a type $z'\in \pre(z,\be')$ has a chance of becoming an extension counted by $D_1$ if on the $i$th step, $\phi'(\be')$ is the edge colored, and the color chosen for this edge is compatible with $z$. Thus, we will write $D_1\leq \sum X(e',c')$, where the sum first goes over all choices of $\be'$, then over all $z'\in \pre (z,\be')$, and finally over all $\phi'\in \Ext(z',\phi)$, and where, given one such fixed triple $(\be',z',\phi')$, the edge $e'=\phi'(\be')$ and we have that $z'(\be')=\unc{c'}$, and $X(e',c')$ is the indicator random variable that is $1$ if $e'$ gets colored $c'$ on step $i$ (note that we have $\leq$ for this expression for $D_1$ here because even if $\phi'(\be')$ is assigned a color compatible with $z$, this assignment could cause a color to be forbidden on a uncolored edge of $\phi'$, which may result in the extension not being of type $z$. We will account for this a little later using an error term $D_1^{-}$). Then,

\begin{align} \label{eqn:D11}
     \E[ D_1| \mc{H}_i]  & \le \sum_{ \substack{\be'\in E_C \\ z' \in \pre(z, \be') \\ \phi' \in \Ext(z', \phi)}}\P[X(e',c')\mid\mc{H}_i]=\sum_{ \substack{\be'\in E_C \\ z' \in \pre(z, \be') \\ \phi' \in \Ext(z', \phi)}} \frac{1}{\binom n2 (1-t) |\Ava(\phi'(\be'))|} 
\end{align}
Now let us estimate the above. For a fixed $\be'$ colored by $z$, if $z(\be')= \col{c}$ for a real or repeated Platonic color $c$, then $|\pre(z,\be')|=1$ and $z' \in \pre(z, \be')$ has $\ell'=\ell-1$ colored edges, and it has $r'$ repeats and $k'$ real colors used on colored edges where $r'+k' = r+k-1$ (we have either lost a repeat or a real color). Therefore the number of terms in \eqref{eqn:D11} corresponding to this fixed $\be'$ is 
\begin{align}
\Ext(z', \phi) &\le  \ext_{z'}(t) +n^{s-2}|C|^{-r-k+1}f_{z'}(t)\nn\\
&= n^{s-2}|C|^{-r-k+1} \sbrac{ t^{\ell-1} (1-t)^{\binom{s}{2}-\ell}a(t) ^{\binom{s}{2}-\ell} +f_{z'}(t)}.\label{eqn:zprime}
\end{align}
Meanwhile, if $z(\be') = \col{c}$ for a nonrepeated Platonic color $c$, then $|\pre(z, \be')|=|C|$ and each $z'' \in \pre(z, \be')$ has $\ell''=\ell-1$ colored edges,  $r''=r$ repeats and $k''=k$ real colors used on colored edges. Therefore the number of terms in \eqref{eqn:D11} corresponding to this fixed $\be'$ is
\begin{align}
   \sum_{z''\in\pre(z,\be')}\Ext(z'', \phi) &\le |C|(\ext_{z'} + n^{s-2}|C|^{-r-k}f_{z'}(t))\nn\\
  &=  n^{s-2}|C|^{-r-k+1} \sbrac{ t^{\ell-1} (1-t)^{\binom{s}{2}-\ell}a(t) ^{\binom{s}{2}-\ell} +f_{z'}(t)}\label{eqn:zdubprime}
\end{align}
where $z'$ is any element of $\pre(z,\be')$, noting that the values of $\ext_{z'}$ and $f_{z'}$ is the same regardless of which one we choose. Note that the value of $f_{z'}$ in line \eqref{eqn:zprime} is the same as $f_{z''}$ in line \eqref{eqn:zdubprime} (see \eqref{eqn:pickfz}), and so in either case \eqref{eqn:zprime} and \eqref{eqn:zdubprime} are equal, so the number of terms in \eqref{eqn:D11} is at most $|E_C|=\ell$ times the expression in \eqref{eqn:zprime}. Therefore line \eqref{eqn:D11} tells us that 
\begin{align}
    \E[ D_1 | \mc{H}_i]  & \le  \ell \cdot  \frac{n^{s-2}|C|^{-r-k+1} \sbrac{ t^{\ell-1} (1-t)^{\binom{s}{2}-\ell}a(t) ^{\binom{s}{2}-\ell} +f_{z'}(t)}}{\binom n2 (1-t)(\ava(t)-|C|f_\Ava(t))} \nn \\
     & = \frac{n^{s-2} |C|^{-r-k}}{\binom n2}\sbrac{\frac{\ell  t^{\ell-1} (1-t)^{\binom{s}{2}-\ell-1} a(t)^{\binom{s}{2}-\ell-1}+\frac{\ell f_{z'}(t)}{(1-t)a(t)} }{1-\frac{f_\Ava(t)}{a(t)}}}\nn\\
     & \le \frac{n^{s-2} |C|^{-r-k}}{\binom n2} \sbrac{\rbrac{1 + 2 \frac{f_\Ava(t)}{a(t)} }  \ell  t^{\ell-1} (1-t)^{\binom{s}{2}-\ell-1} a(t) ^{\binom{s}{2}-\ell-1} +  \frac{2\ell f_{z'}(t)}{(1-t)a(t)} }\label{eqn:ed1est}
     \end{align}
     where the last line follows from the inequalities 
     \[
     \frac{1}{1-\frac{f_\Ava(t)}{a(t)}} \le 1 + \frac{2f_\Ava(t)}{a(t)} \le 2.
     \]
     Now, using \eqref{def:g1}, we can write \eqref{eqn:ed1est} as
     \begin{align}
     & \frac{n^{s-2} |C|^{-r-k}}{\binom n2} \sbrac{ \rbrac{1 + 2 \frac{f_\Ava(t)}{a(t)} }  g_1(t) +  \frac{2\ell f_{z'}(t)}{(1-t)a(t)} }\nn\\
     & \le  \frac{n^{s-2} |C|^{-r-k}}{\binom n2}  \sbrac{  g_1(t) + 2p^2  a(t)^{\binom s2 - \ell - 2} f_\Ava(t) +   \frac{2\ell a(t)^{\binom s2 - \ell - 1}f_{\Ext}(t)}{(1-t)}}, \label{eqn:ED1}
\end{align}
where the above inequality follows from \eqref{eqn:g estimate} and \eqref{eqn:pickfz}.

While the above work will be enough for us to establish that $\Ext^+(z, \phi)$ is a supermartingale, there is some more work required to show that $\Ext^-(z, \phi)$ is a submartingale. In particular we would also need a lower bound on $\E[ D_1 | \mc{H}_i]$. The bound given above overcounts situations in which an edge is colored in such a way that a clique would become the right type if it was not for the fact that the newly colored edge also caused a color on one of the uncolored edges of that clique to become forbidden. To deal with this, we will let $D_1^-$ be the number of such extensions that do not become type $z$ for that reason. Note that we have
\begin{align}
     \E[ D_1+D_1^-| \mc{H}_i]  &= \sum_{ \substack{\be'\in E_C \\ z' \in \pre(z, \be') \\ \phi' \in \Ext(z', \phi)}} \frac{1}{\binom n2 (1-t) |A_{\phi'(\be')}|}\nonumber\\
     &\geq  \frac{n^{s-2} |C|^{-r-k}}{\binom n2}  \sbrac{  g_1(t) - 2p^2  a(t)^{\binom s2 - \ell - 2} f_\Ava(t) -   \frac{2 \ell a(t)^{\binom s2 - \ell - 1}f_{\Ext}(t)}{(1-t)}}
     \label{eqn:D1submtg},
\end{align}
where the inequality follows from essentially the same work used to derive \eqref{eqn:ED1}, but applied to a lower bound. Now, to get a lower bound on $E[D_1|\mc{H}_i]$, it will suffice to subtract an upper bound on $E[D_1^-|\mc{H}_i]$ from \eqref{eqn:D1submtg}.

\begin{claim}\label{clm:D1-}
In the good event we have
   \begin{equation}\label{eqn:D1-bound}
     \E[ D_1^-| \mc{H}_i]  \le n^{\pow(\bS, \be, z) -2 - \m + 30p^2\eps}.
\end{equation} 
\end{claim}

\begin{proof}
    We intend to use Lemma \ref{lem:appendpreforb}. In order for an extension $\phi'$ to be counted in $D_1^-$ the set $S_1 = \phi'(\bS)$  fits some type $z_1 \in \pre(z, \be_2)$ which is a predecessor of $z$ rooted at $\be_1:=\be$ and $\be_2$ is some edge colored by $z$. 
    
    Say $z(\be_2)=\col{c_2}$, and first consider the case when $|\pre(z, \be_2)|=1$, i.e.\ the case where either $c_2$ is either real or a Platonic repeated color in $z$. Then we have only one choice for $z_1$ and $z_1(\be_2)=\unc{c_2}$. Note that in this case we have $\pow(\bS, \be, z_1)= \pow(\bS, \be, z) + \f$ since $z_1$ has one fewer coincidence than $z$. There must be some $(\be_2, c_2, \be_3, c_3)$-preforbidder $z_2$ on say $\bS_2$, compatible with $z_1$, where $z_1(\be_3)=\unc{c_3}$. Indeed, the only real colors that would potentially need to appear on $z_2$ are $c_2, c_3$ and possibly one real color that appears on $z_1$. Note that in this case we must have $\bS_2\setminus \bS_1\neq \emptyset$ since if it were the case that $\bS_2\subseteq \bS_1$, then the fact that $(z_2)_{\be_2,\be_3}$ is not legal (see Definition~\ref{def:ecpreforb}~\eqref{def:ecpreforb3}) would contradict Observation~\ref{obs:trackable}~\eqref{obs:trackable1}.
    
    We bound the contribution to $\E[ D_1^-| \mc{H}_i]$ from this case as follows. There are a constant number of choices for the predecessor $z_1$ and the preforbidder $z_2$. For each fixed choice $z_1, z_2$ we have by Lemma \ref{lem:appendpreforb} and Corollary \ref{cor:crudepractical} that 
    \[
    |\Ext(z_1 \cup z_2, \phi)| \le n^{\pow(S,  \be, z_1) - \m + 20p^2\eps} =  n^{\pow(S, \be, z) + \frac{p-3}{\binom p2 - q + 1}+ 20p^2\eps}.
    \]
    Now for a fixed extension $\phi' \in \Ext(z_1 \cup z_2, \phi)$ to contribute to $D_1^-$ at step $i$, note that we need to color a fixed edge a fixed color. Indeed, we are in the case where the predecessor $z_1$ was formed by uncoloring an edge whose color $c_2$ was either real or a Platonic repeat. If $c_2$ is real then the color chosen at step $i$ must be $c_2$. If $c_2$ was a Platonic repeat in $z$ then $z_1$ still has an edge colored $c_2$, and in this extension $\phi'$ we must have that $c_2$ is represented by some fixed real color which must then be the color chosen at step $i$. Coloring a fixed edge a fixed color  has probability at most 
    \[
  O\rbrac{ \frac{1}{(1-t)\binom n2 |C| a(t)}} \le n^{-2-\f + 3\eps},
  \]
    and so the contribution to $\E[ D_1^-| \mc{H}_i]$ from this case is at most
    \[
    O\rbrac{n^{\pow(S, \be, z) + \frac{p-3}{\binom p2 - q + 1}+ 20p^2\eps} \cdot n^{-2-\f + 3\eps}} =  O\rbrac{n^{\pow(\bS, \be, z) -2 - \m + 30p^2\eps}}.
    \]
    
    Now consider the case where $|\pre(z, \be_2)|>1$, i.e.\ $c_2$ is a Platonic color that is nonrepeated in $z$. In this case there are $|C| = O\rbrac{n^{\f}}$ choices for the predecessor $z_1$, and for any $z_1$ we have  $\pow(\bS, \be, z_1)= \pow(\bS, \be, z)$. We still have a constant number of choices for the preforbidder $z_2$, and for each fixed choice $z_1, z_2$ we have by Lemma \ref{lem:appendpreforb} and Corollary \ref{cor:crudepractical} that 
    \[
    |\Ext(z_1 \cup z_2, \phi)| \le n^{\pow(S,  \be, z_1) - \m+ 20p^2\eps} =  n^{\pow(S, \be, z) -\m + 20p^2\eps}.
    \]
    Here each extension $\phi' \in \Ext(z_1 \cup z_2, \phi)$ represents a potential contribution to $D_1^-$ at step $i$, and for this potential contribution to be realised we need to color a fixed edge a fixed color.
    Altogether, the contribution to $\E[ D_1^-| \mc{H}_i]$ from this case is at most
    \[
    O\rbrac{n^{\f} \cdot n^{\pow(S, \be, z) - \m + 20p^2\eps} \cdot n^{-2-\f + 3\eps}} =  O\rbrac{n^{\pow(\bS, \be, z) -2 - \m + 30p^2\eps}}.
    \]
    Summing the contributions to $\E[ D_1^-| \mc{H}_i]$ gives the bound \eqref{eqn:D1-bound}.
\end{proof}

Subtracting \eqref{eqn:D1-bound} from \eqref{eqn:D1submtg} yields
\begin{align*}
    \E[ D_1 | \mc{H}_i] & \ge \frac{n^{s-2} |C|^{-r-k}}{\binom n2}  \sbrac{  g_1(t) - 2p^2  a(t)^{\binom s2 - \ell - 2} f_\Ava(t) -   \frac{2a(t)^{\binom s2 - \ell - 1}f_{\Ext}(t)}{(1-t)}}\nn\\
    & \qquad \qquad+ O\rbrac{ n^{\pow(\bS, \be, z) -2 - \m + 30p^2 \eps}}.
\end{align*}

Via a similar calculation, we can show that
\begin{align}
    \E[ D_1 | \mc{H}_i] &\leq \frac{n^{s-2} |C|^{-r-k}}{\binom n2}  \sbrac{  g_1(t) + 2p^2  a(t)^{\binom s2 - \ell - 2} f_\Ava(t) +   \frac{2a(t)^{\binom s2 - \ell - 1}f_{\Ext}(t)}{(1-t)}}\nn\\
    & \qquad \qquad + O\rbrac{ n^{\pow(\bS, \be, z) -2 - \m + 30p^2 \eps}}.\label{equation ED1}
\end{align}

It will be useful to have a crude upper bound on $\E[ D_1 | \mc{H}_i]$. To that end, using the simple bound $a(t)\leq 1$ with \eqref{eqn:g estimate}, we can see that $g_1(t)\leq p^2$. Furthermore, using $1-t\geq n^{-\epsilon}$, we can see from \eqref{eqn:pickfava} that $f_\Ava\leq n^{-\hm + (25p^2+10p^4)\eps}=o(1)$, and from \eqref{eqn:pickfext} that $\frac{a(t)\cdot f_\Ext}{1-t}\leq n^{-\hm + (25p^2+10p^4)\eps}=o(1)$. Then we can rewrite \eqref{equation ED1} as
\begin{align}
    \E[ D_1 | \mc{H}_i] \leq \frac{n^{s-2} |C|^{-r-k}}{\binom n2}  \sbrac{p^2 + o(1)} + O\rbrac{ n^{\pow(\bS, \be, z) -2 - \m + 30p^2 \eps}}\leq \frac{n^{s-2} |C|^{-r-k}}{\binom n2}(2p^2).\label{equation ED1 2}
\end{align}

\subsubsection{Bounding \texorpdfstring{$D_1$}{1}}

Here we give an absolute bound on how large $D_1$ can be in the good event. Let $e_i\in \binom{[n]}{2}$. We want to bound $D_1$ when we color the edge $e_i$ the color $c_i$. We can bound this by considering predecessors of $z$. 

First let us consider the case where $c_i$ is represented by itself in $z$ (i.e.\ when $z(\be')=\col{c_i}$ where $\be'$ is the Platonic edge corresponding to $e_i$). Let $z'$ be the unique element of $ \pre(z,\be')$. Now, $D_1$ is bounded by $\Ext(z', \phi')$ where $\bD(\phi') = \be \cup \be'$. Towards using Corollary \ref{cor:crudepractical}, let $\bS'$ be such that $\be\cup \be'\subseteq \bS' \subseteq \bS$. Then,
\begin{align*}
\pow(\bS,\bS',z')&=\pow(\bS,\be,z')-\pow(\bS',\be,z')\\
&=\pow(\bS,\be,z)+\f-\left(\pow(\bS', \be,z)+\f \right)\\
&=\pow(\bS,\be,z)-\pow(\bS', \be,z)\le \pow(\bS,\be,z)-\m,
\end{align*}
where the first equality follows from additivity, then the second follows from the fact that $z'$ is a predecessor of $z$ in which an edge with a real color $\be'\in \bS' \subseteq \bS$ was uncolored, and then the final inequality since $|\bS'|\geq 3$ and $z$ is a trackable type, so $\pow(\bS',e,z)\geq \m$. Thus by Corollary \ref{cor:crudepractical}, we have that
\[
D_1\leq |\Ext(z',\phi')|\leq n^{\pow(\bS,\be,z)-\m + 20p^2\eps}.
\]

Now consider the case where  $z(\be')=\col{\bc'}$, where $\bc'$ is a Platonic color corresponding to the real color $c_i$. Let $z''$ be the (not necessarily trackable) type that is identical to $z$ except  every label of the form $\unc{\bc'}$ or $\col{\bc'}$ is replaced with $\unc{c_i}$ or $\col{c_i}$ respectively. We let $z'''$ be the unique element of $\pre(z'', \be')$.  Any extension that enters $\Ext(z, \phi)$ at step $i$ must have previously been in $\Ext(z''', \phi')$ where $\bD(\phi') = \be \cup \be'$. Thus we want to bound the size of $\Ext(z''', \phi')$, so again, towards using Corollary~\ref{cor:crudepractical}, we let $\bS'$ be such that $\be\cup \be'\subseteq \bS' \subseteq \bS$. Then by the construction of $z''$ and $z'''$ we have \[
\pow(\bS ,\be,z) = \pow(\bS ,\be,z'') - \f  = \pow(\bS ,\be,z'''),
\]
since going from $z$ to $z''$ replaces all instances of a Platonic color with a real color that did not appear on $z$, and that does appear on at least one colored edge of $z''$ (namely $\be'$), thus giving us exactly one new coincidence, and then going from $z''$ to $z'''$ uncolors an edge that was colored a real color, causing exactly one less coincidence to occur. Similarly $\pow(\bS' ,\be,z''') =  \pow(\bS' ,\be,z)$, so
\begin{align*}
\pow(\bS,\bS',z''')&=\pow(\bS ,\be,z''')-\pow(\bS',\be,z''')\\
&=\pow(\bS,\be,z)-\pow(\bS',\be,z)\\
&\le \pow(\bS,\be,z)-\m,
\end{align*}
where the last inequality follows from the fact that $\pow(\bS',\be,z)\geq \m$ since $z$ is a trackable type and $|\bS'|\ge 3$. Thus we get
\[
D_1\leq |\Ext(z''',\phi')|\leq n^{\pow(\bS,\be,z)-\m + 20p^2 \eps}.
\]

Thus in either case, 
\begin{equation}\label{eqn:D1bound}
    D_1 \leq  n^{\pow(\bS,\be,z)-\m + 20p^2 \eps}.
\end{equation}

\subsection{\texorpdfstring{$D_2$}{2}}

\subsubsection{Estimating  \texorpdfstring{$ \E[D_2| \mc{H}_i]$}{2}}

It is easy to see that
\begin{align}\label{eqn:ED2}
\E[ D_2 | \mc{H}_i] &= \frac{\binom{s}2 - \ell - 1}{\binom n2 (1-t)}\Ext(z, \phi)  \nn\\
&\ge  \frac{\binom{s}2 -\ell - 1}{\binom n2 (1-t)} \rbrac{n^{s-2}|C|^{-r-k}t^\ell (1-t)^{\binom{s}{2}-\ell-1}a(t) ^{\binom{s}{2}-\ell-1} -n^{s-2}|C|^{-r-k}f_z }\nn\\
& \ge  \frac{n^{s-2}|C|^{-r-k}}{\binom n2 } \sbrac{ g_2(t) - \frac{p^2 a(t) ^{\binom{s}{2}-\ell-1}f_\Ext}{1-t} },
\end{align}
where the last line follows from \eqref{eqn:pickfz}, \eqref{def:g2} and $\binom{s}{2}-\ell-1 \le p^2$.
Similarly, 
\begin{align}
\E[ D_2 | \mc{H}_i] &\le  \frac{n^{s-2}|C|^{-r-k}}{\binom n2 } \sbrac{ g_2(t) + \frac{p^2 a(t) ^{\binom{s}{2}-\ell-1}f_\Ext}{1-t} }.\label{equation ED2}
\end{align}

We will need a crude upper bound on $\E[ D_2 | \mc{H}_i] $, so we note that using the simple bound $a(t)\leq 1$ and \eqref{eqn:g estimate}, we have that $g_2(t)\leq p^2$, and then using $1-t\geq n^{-\epsilon}$, we can see from \eqref{eqn:pickfext} that $\frac{a(t)\cdot f_\Ext}{1-t}\leq n^{-\hm + (25p^2+10p^4)\eps}=o(1)$. Then, using this we can rewrite \eqref{equation ED2} as
\begin{align}
\E[ D_2 | \mc{H}_i] &\le  \frac{n^{s-2}|C|^{-r-k}}{\binom n2 } \sbrac{ p^2 + o(1)}\leq \frac{n^{s-2}|C|^{-r-k}}{\binom n2}(2p^2) .\label{equation D2 2}
\end{align}

\subsubsection{Bounding \texorpdfstring{$D_2$}{2}}

Any extension $\phi'\in \Ext(z,\phi)$ that is counted by $D_2$ must contain both $\phi(\be)$ and the edge $e_i$ that is being colored on this step. Thus, $D_2$ is at most $\Ext( z , \phi'')$, where $\bD(\phi'')=\be\cup \be'$ (where $\phi''(\be')=e_i$). For any such $\phi'$ and any $\bS'$ with $\bD(\phi'') \subseteq \bS' \subseteq \bS$, we have
\[
\pow(\bS, \bS', z) = \pow(\bS, \be, z) - \pow(\bS', \be, z) \le  \pow(\bS, \be, z) - \m,
\]
where the last inequality follows from Observation~\ref{obs:trackable}~\eqref{obs:trackable2} since $\be\subsetneq\bS'$ and $z$ is a trackable type. Thus, by Corollary \ref{cor:crudepractical} we have that 
\begin{equation}\label{eqn:D2bound}
  D_2 \leq |\Ext(z,\phi'')|\leq  n^{\pow(\bS,\be,z)-\m + 20p^2 \eps}.  
\end{equation}

\subsection{\texorpdfstring{\href{https://www.imdb.com/title/tt0116000/}{$D_3$}}{2}}

\subsubsection{Estimating  \texorpdfstring{$ \E[D_3| \mc{H}_i]$}{2}}
Note first that if $\ell = \binom s2 - 1$ then $\be$ is the only edge in $\bS$ that is uncolored, and therefore it is impossible for $\Ext(z, \phi)$ to decrease due to a color being forbidden at an edge. In this case we deterministically have $D_3=0 = \frac{n^{s-2}|C|^{-r-k}}{\binom n2}g_3(t)$. For the rest of this subsection assume that $\ell \le \binom s2 - 2$.

 Recall that $\mc{F}(e', c')=\mc{F}(e', c', i)$ is the event that the color $c'$ is forbidden at $e'$ at step $i$. Consider some $\phi' \in \Ext(z, \phi)$ and some uncolored edge $e'\neq \phi(\be)$ in $\phi'(\bS)$ (where $\bS = \bS(z)$), and let $\psi$ be the color map of $\phi'$ (see Definition~\ref{def:ext}). The type $z$ indicates that some color $c'$ is available at $\be'$ (where $\be'$ is the Platonic edge such that $\phi'(\be')=e'$). The color $c'$ might be Platonic, but in that case $\phi'$ must use some real color $\psi(c')$ in place of $c'$. To ease the notation,  we define the event $\mc{F}_{\phi'}(\be',c'):=\mc{F}(\phi'(\be'),\psi(c'))$. Thus, we want to know for which embeddings $\phi'$ does the event $\mc{F}_{\phi'}(\be',c')$ happen on at least one edge $\be'$ in $\bS$, but also want to avoid counting any embeddings which end up having an uncolored edge become colored on this step (as those are counted in $D_2$). We will let $D_3^-$ denote the embeddings $\phi'$ for which both $\mc{F}_{\phi'}(\be',c')$ happens to at least one edge $\be'\subseteq \bS$, and an edge in $\phi'(\bS)$ is assigned a color on this step. Then $D_3+D_3^-$ simply counts embeddings removed from $\Ext(z,\phi)$ because a color that was originally prescribed as available at an edge was forbidden at an edge (regardless of if an edge in the embedding was also colored). Thus we have
\begin{align*}
 \E[ D_3 +D_3^-| \mc{H}_i] =  \sum_{\phi' \in \Ext(z, \phi)} \P\sbrac{\bigcup_{(\be', c'): z(\be') = \unc{c'}} \mc{F}_{\phi'}(\be',c')}.
\end{align*}
 Now we use the simple bounds for any events $E_j$ with $j\in J$ for an index set $J$,
 \[
 \sum_{j \in J} \P[E_j] \ge \P\sbrac{\bigcup_{j \in J} E_j} \ge \sum_{j \in J} \P[E_j] - \sum_{\substack{j, j' \in J \\ j \neq j'}} \P[E_j \cap E_{j'}],
 \]
 which imply
 \[
 \P\sbrac{\bigcup_{j \in J} E_j} = \sum_{j \in J} \P[E_j] + O\rbrac{ \sum_{\substack{j, j' \in J \\ j \neq j'}} \P[E_j \cap E_{j'}]}.
 \]
Thus we have 
\begin{align*}
 \E[ D_3 +D_3^-| \mc{H}_i] =  \sum_{\substack{ \phi' \in \Ext(z, \phi) \\ (\be', c'): z(\be') = \unc{c'}}} \P\sbrac{\mc{F}_{\phi'}(\be',c')} + O(B_3),
\end{align*}
where
\[
B_3 := \sum_{\substack{\phi' \in \Ext(z, \phi)\\ (\be', c'): z(\be') = \unc{c'} \\ (\be'', c''): z(\be'') = \unc{c''}  \\ \be' \neq \be''}}  \P\sbrac{\mc{F}_{\phi'}(\be', c') \cap \mc{F}_{\phi'}(\be'', c'')}.
\]
Now we will bound $B_3$. 
\begin{claim}
\[
B_3 =  O\rbrac{n^{\pow(\bS, \be, z) - 2 - \m + 40p^2 \eps}}. 
\]
\end{claim}
\begin{proof}
First we claim that for any $e', c', e'', c''$ with $e' \neq e''$ there are at most $n^{\frac{p-3}{\binom p2 - q + 1} + 20p^2\eps}$ pairs  $e''', c'''$ such that assigning $c'''$ to $e'''$ would simultaneously forbid $c'$ at $e'$ and $c''$ at $e''$. Indeed, let $S_1 \in S_{e', c'}(e''', c''')$ and consider the following two cases (handled in the following paragraphs): either $e'' \subseteq S_1$, or else $e'' \not\subseteq S_1$ in which case there must be some $S_2 \in S_{e'', c''}(e''', c''')$ with $S_1 \neq S_2$.

Consider first the case where $e'' \subseteq S_1$. The number of possible sets $S_1$ is bounded by the sum of a constant number of terms of the form $O(\Ext(y, \rho))$ where $y$ is a colors-only $(\be', c', \be''', c^*)$-preforbidder on $\bS_1$ (where $c^*=c'$ if $c'''=c'$ and otherwise $c^*$ is a Platonic color representing $c'''$) and $\rho$ has domain  $\be' \cup \be''$. Assume that $y$ is chosen so that no real colors other than possibly $c'$ appear on colored edges. We have by Observation \ref{obs:preforbpow} that $\pow(\bS_1, \be', y) \le \f.$ Also since $y$ is the colored portion of some trackable type rooted at $\be'$, by Observation~\ref{obs:trackable}~\eqref{obs:trackable2} we have  $\pow(\bS', \be', y) \ge \m$ whenever $\be' \subsetneq \bS' \subseteq \bS_1.$ Therefore
\[
\maxpow(\bS_1, \be' \cup \be'', y) \le \mf
\]
and so $\Ext(y, \rho)= O\rbrac{n^{\mf + 20p^2\eps}}$ by  Corollary~\ref{cor:crudepractical}. Note that there is a constant number of choices for $y$, and each element of $\Ext(y, \rho)$ represents a constant number of choices for $(e''', c''')$. Thus the number of choices for $(e''', c''')$ is $O\rbrac{n^{\mf+20p^2\eps}}$.

Now consider the second case, where $e'' \not \subseteq S_1$. We will use Lemma~\ref{lem:union}. $S_1$ must fit some colors-only  $(\be', c', \be''', c^*)$-preforbidder $y_1$ on $\bS_1$ (where as before, $c^*=c'$ if $c'''=c'$ and otherwise $c^*$ is a Platonic color representing $c'''$). Likewise $S_2$ must fit some colors-only  $(\be'', c'', \be''', c^*)$-preforbidder $y_2$ on $\bS_2$. Clearly we can choose $y_1, y_2$ to be compatible. Then the number of possible sets $S_1 \cup S_2$ is bounded by a constant number of terms of the form $\Ext(y_1 \cup y_2, \rho)$ where $\bD(\rho) = \be' \cup \be''.$  Since $e'' \not\subseteq S_1$ we have $\be'' \not \subseteq \bS_1$. Thus by Lemma~\ref{lem:union} and Observation \ref{obs:preforbpow} we have
\[
\maxpow(\bS_1 \cup \bS_2, \be' \cup \be'', y_1 \cup y_2) \le \pow(\bS_1 , \be' , y_1 ) - \m \le \mf.
\]
Thus by Corollary \ref{cor:crudepractical} we conclude that $\Ext(y_1\cup y_2, \rho)=O\rbrac{n^{\mf+ 20p^2\eps}}.$ Since each possible set $S_1 \cup S_2$ counted by $\Ext(y_1\cup y_2, \rho)$ corresponds to a constant number of choices for $(e''', c''')$ we have that the number of such choices is also $O\rbrac{n^{\mf+20p^2\eps}}.$

Adding together the bounds from the two cases, we get that the total number of choices for pairs $(e''',c''')$ is $O\rbrac{n^{\mf+20p^2\eps}}$, so we have
\[
\P\sbrac{\mc{F}(e', c') \cap \mc{F}(e'', c'')} \le O\rbrac{\frac{n^{\mf+20p^2\eps}}{\binom n2 (1-t) |C|a(t) }} \le n^{-2 - \m + 30p^2 \eps}.
\]
Now since the number of terms in $B_3$ is at most $O(\Ext(z, \phi)) = \TO\rbrac{n^{\pow(\bS, \be, z)} }$ (by Condition~\eqref{cond:exttraj} from Definition~\ref{definition of the good event} and the fact that $z$ is trackable), we have that 
\[
B_3 =  O\rbrac{n^{\pow(\bS, \be, z) - 2 - \m + 40p^2 \eps}}. 
\]
\end{proof}

Now we bound $\E[D_3^-| \mc{H}_i]$.  
\begin{claim}\label{clm:D3-}
In the good event we have
   \begin{equation*}
     \E[ D_3^-| \mc{H}_i]  = O\rbrac{n^{\pow(\bS, \be, z) -2 - \f + 30p^2 \eps}}.
\end{equation*} 
\end{claim}

The proof is very similar to that of Claim \ref{clm:D1-}. 

\begin{proof}
     In order for an extension $\phi'$ to be counted in $D_3^-$ there must be two sets $S_1 = \phi'(\bS)$, which fits the type $z_1=z$ rooted at $\be_1:=\be$, and $S_2$, which fits a preforbidder type. Say $S_2$ fits the $(\be_2, c_2, \be_3, c_3)$-preforbidder type $z_2$ on $\bS_2$ and compatible with $z_1$, where $\be_2$ represents the edge (say $e'$) getting a color forbidden causing $\phi'$ to be counted in $D_1^-$ and $\be_3$ represents $e_i$. First consider the case where $\bS_2 \not\subseteq \bS_1$. For this case we will use Lemma \ref{lem:appendpreforb}. Since $e' \cup e_i \subseteq S_1 \cap S_2$ we have $\be_2 \cup \be_3 \subseteq \bS_1 \cap \bS_2$. Thus, if we let $y:=z_1\cup z_2$, applying Lemma \ref{lem:appendpreforb} we have 
    \[
   \maxpow(\bS_1 \cup \bS_2, \be, y) \le \pow(\bS, \be, z_1) - \m = \pow(\bS, \be, z) - \m.
   \]
 Thus by Corollary \ref{cor:crudepractical} we have that
    \[
    |\Ext(y, \phi)| \le  O\rbrac{n^{\pow(\bS, \be, z)- \m+ 20p^2\eps}}.
    \]
    Note there is a constant number of choices for the colored portion of $y_2$. There is a constant number of relevant types $y$, and each element of $\Ext(y, \phi)$ represents some constant number of potential ways to get a contribution to $D_3^-$. Each such way has probability at most $n^{-2-\f + 3\eps}$.
    Thus the expected contribution to $D_3^-$ from this case is at most
    \[
     O\rbrac{n^{\pow(\bS, e, z)-\m+ 20p^2\eps}} \cdot n^{-2-\f + 3\eps} \le n^{\pow(\bS, \be, z) -2 - \frac{p-1}{\binom p2 - q + 1} + 30p^2 \eps}
    \]
    which is even smaller than we need. 
    
    Now consider the case where $\bS_2 \subseteq \bS_1$. In this case each element of $\Ext(z, \phi)$ represents some constant number of potential ways to get a contribution to $D_3^-$. Thus the expected contribution to $D_3^-$ from this case is at most
    \[
     O\rbrac{|\Ext(z, \phi)|} \cdot n^{-2-\f + 3\eps} \le n^{\pow(\bS, \be, z) -2 - \f + 30p^2 \eps},
    \]
    where we use Condition~\eqref{cond:exttraj} from Definition~\ref{definition of the good event} to bound $|\Ext(z, \phi)|$.
\end{proof}

With the claim in hand, we have
\[
 \E[D_3 | \mc{H}_i] \ge \sum_{\substack{ \phi' \in \Ext(z, \phi) \\ (\be', c'): z(\be') = \unc{c'}}} \P\sbrac{\mc{F}_{\phi'}(\be',c')} + O\rbrac{n^{\pow(\bS, \be, z) -2 - \m + 40p^2 \eps}}.
\]
We then use Condition~\eqref{cond:exttraj} from Definition~\ref{definition of the good event} to lower bound the number of terms in the above sum, and \eqref{equation for PF(e,c)} to write
\begin{align}
   & \E[D_3 | \mc{H}_i]\nn\\
    & \ge \rbrac{\ext_z - n^{s-2}|C|^{-r-k}f_z} \cdot \frac{\binom s2 -1- \ell }{\binom n2 } \sbrac{ h(t) - 2p^2 \eps\log n \cdot \frac{f_\Ava}{a(t)}- \eps \log n \cdot \frac{f_\Ext}{1-t}} \nn\\
    & \qquad \qquad +O\rbrac{n^{\pow(\bS, \be, z) -2 - \m + 40p^2 \eps}}\nn\\
    & = \rbrac{\binom s2 -1- \ell }\frac{n^{s-2}|C|^{-r-k}}{\binom n2} a(t) ^{\binom{s}{2}-\ell-1}\sbrac{ t^\ell (1-t)^{\binom{s}{2}-\ell-1} - f_\Ext}\nn \\
    & \qquad \qquad \cdot \sbrac{ h(t) - 2p^2 \eps a(t)^{-1} \log n \cdot f_\Ava- \eps \log n \cdot \frac{f_\Ext}{1-t}} +O\rbrac{n^{\pow(\bS, \be, z) -2 - \m + 40p^2 \eps}}\nn\\
    & \ge \frac{n^{s-2}|C|^{-r-k}}{\binom n2} \Bigg[ g_3(t)  - 2p^4\eps a(t) ^{\binom{s}{2}-\ell-2} (1-t)\log n \cdot  f_\Ava - 2p^2\eps  a(t) ^{\binom{s}{2}-\ell-1} \log n \cdot f_\Ext\nn\\
    & \qquad \qquad + \TO\rbrac{f_\Ext f_\Ava + \frac{f_\Ext^2}{1-t}} \Bigg] + O\rbrac{n^{\pow(\bS, \be, z) -2 - \m + 40p^2 \eps}}\nn \\
    & \ge \frac{n^{s-2}|C|^{-r-k}}{\binom n2} \sbrac{g_3(t)  - 2p^4\eps a(t) ^{\binom{s}{2}-\ell-2} (1-t)\log n \cdot  f_\Ava - 2p^2\eps  a(t) ^{\binom{s}{2}-\ell-1} \log n \cdot f_\Ext }\nn\\
    & \qquad \qquad + O\rbrac{n^{\pow(\bS, \be, z) -2 - \m + (20p^4 + 60p^2 +2)\eps}}\label{eqn:ED3}
\end{align}
On the second-to-last line we used \eqref{def:g3}, the fact that $h(t) \le \eps \log n$, $\binom s2 - 1 - \ell \le p^2$, and that $\ell \le \binom s2 - 2$ so $t^\ell (1-t)^{\binom{s}{2}-\ell-1} \le 1-t.$ On the last line we have used that 
\[
f_\Ext f_\Ava + \frac{f_\Ext^2}{1-t} = \TO\rbrac{n^{ - \m + (20p^4 + 60p^2 +1)\eps}}.
\]

Similarly, we can write
\begin{align}
   \E[D_3 | \mc{H}_i]& \leq \frac{n^{s-2}|C|^{-r-k}}{\binom n2} \sbrac{g_3(t)  + 2p^4\eps a(t) ^{\binom{s}{2}-\ell-2} (1-t)\log n \cdot  f_\Ava + 2p^2\eps  a(t) ^{\binom{s}{2}-\ell-1} \log n \cdot f_\Ext }\nn\\
    & \qquad \qquad + O\rbrac{n^{\pow(\bS, \be, z) -2 - \m + (20p^4 + 60p^2 +2)\eps}}\label{eqn:ED3part2}
\end{align}

We will need a crude upper bound on $ \E[D_3 | \mc{H}_i]$, so using the simple bound $a(t)\leq 1$, along with \eqref{eqn:g estimate} and \eqref{eqn:hbound}, we have that $g_3(t)\leq p^2h(t)\leq p^2\epsilon\log n$. Furthermore, using $1-t\geq n^{-\epsilon}$, we can see from \eqref{eqn:pickfava} that $f_\Ava\leq n^{-\hm + (25p^2+10p^4)\eps}=o(1)$, and from \eqref{eqn:pickfext} that $a(t)\cdot f_\Ext\leq n^{-\hm + (25p^2+10p^4)\eps}=o(1)$. Using the above bounds, we can rewrite \eqref{eqn:ED3part2} as
\begin{align}
   \E[D_3 | \mc{H}_i]& \leq \frac{n^{s-2}|C|^{-r-k}}{\binom n2} \sbrac{p^2\epsilon\log n + o(\log n) }+ O\rbrac{n^{\pow(\bS, \be, z) -2 - \m + (20p^4 + 60p^2 +2)\eps}}\nn\\
   &\leq \frac{n^{s-2}|C|^{-r-k}}{\binom n2}(2p^2\epsilon\log n).\label{eqn:ED3part3}
\end{align}

\subsubsection{Bounding \texorpdfstring{$D_3$}{1}}

We describe how to get a contribution to $D_3$. To get a contribution we must have some extension $\phi' \in \Ext(z, \phi)$ whose image is some set $S_1 \supseteq e$ fitting $z$. Assigning $c_i$ to $e_i$ forbids some color $c'$ at $e' \subseteq S_1$, causing $\phi'$ to no longer be in $\Ext(z, \phi)$. So there is some $S_2 \in S_{e_i, c_i}(e', c')$ fitting some colors-only $(\be_2, c_2, \be_3, c_3)$-preforbidder type $y_2$ on say $\bS_2$, where $\be_2, \be_3$ correspond to $e_i$ and $e'$ respectively, $c_2=c_i$ and $c_3$ is either equal to $c'$ if $c'=c_i$ or else $c_3$ is Platonic representing $c'$. Assume we choose $y_2$ compatible with $z$.

We will use Lemma \ref{lem:union}. We use $y_1=z$, $\bS_1=\bS$ and $\be_1=\be$ (we have already defined $\bS_2, y_2, \be_2, \be_3$). Note that $e_i$ cannot be in $S_1$ (otherwise we would actually get a contribution to $D_3^-$ instead of $D_3$), and so $\be_2 \not \subseteq \bS_1$ as required. Lemma \ref{lem:union} then gives us that 
    \[
    \maxpow(\bS_1 \cup \bS_2, \be_1 \cup \be_2, y_1 \cup y_2) \le  \pow(\bS_1, \be_1, y_1) - \m =  \pow(\bS, \be, z) - \m.
    \]
Thus, by Corollary \ref{cor:crudepractical} we have for $\phi''$ with domain $\be_1 \cup \be_2$ that
\[
 |\Ext(y_1 \cup y_2,\phi'')|\leq  n^{\pow(\bS,\be,z)-\m + 20p^2 \eps}.
\]
There are a constant number of choices for $y_2$ and each element of $\Ext(y_1 \cup y_2,\phi'')$ gives a constant sized contribution to $D_3$. Thus
\begin{equation}\label{eqn:D3bound}
  D_3 =O\rbrac{  n^{\pow(\bS,\be,z)-\m + 20p^2 \eps}}.  
\end{equation}

\subsection{Establishing that \texorpdfstring{$\Ext^+(z, \phi)$}{1} is a supermartingale}

Note that since $f_z(t) = a(t)^{\binom{s}{2}-\ell-1}f_\Ext(t)$ we have
\begin{align*}
f_z'(t) &= a(t)^{\binom{s}{2}-\ell-1}f_\Ext'(t) - \rbrac{\binom{s}{2}-\ell-1}a(t)^{\binom{s}{2}-\ell-1}h(t) f_\Ext(t)\\
&\ge a(t)^{\binom{s}{2}-\ell-1} \sbrac{f_\Ext'(t) - p^2 \eps \log n \cdot f_\Ext(t)}.  
\end{align*}

Therefore
\begin{align}
   & \E[\D \Ext^+(z, \phi) | \mc{H}_i]\nn \\
   &= \E[\D \Ext(z, \phi) | \mc{H}_i] - \frac{n^{s-2}|C|^{-r-k}}{\binom n2} \sbrac{ g_1(t) - g_2(t) - g_3(t)+ f_z'(t)} + O\rbrac{n^{\pow(\bS, \be, z) - 4 + \eps}}\nn\\
    & \le \E[D_1 - D_2 - D_3 | \mc{H}_i] - \frac{n^{s-2}|C|^{-r-k}}{\binom n2}\sbrac{g_1-g_2-g_3+ a(t)^{\binom{s}{2}-\ell-1} \sbrac{f_\Ext'(t) - p^2 \eps \log n \cdot f_\Ext(t)}}\nn\\
    & \qquad \qquad + O\rbrac{n^{\pow(\bS, \be, z) - 4 + \eps}} .\label{eqn:extsupermtg1}
\end{align}
Now using \eqref{eqn:ED1} we have
\begin{align}
    \E[D_1 |\mc{H}_i] -\frac{n^{s-2}|C|^{-r-k}}{\binom n2}g_1 &\le  \frac{n^{s-2} |C|^{-r-k}}{\binom n2}  \sbrac{   2p^2  a(t)^{\binom s2 - \ell - 2} f_\Ava(t) +   \frac{2a(t)^{\binom s2 - \ell - 1}f_{\Ext}(t)}{1-t}}\label{eqn:d1g1}
\end{align}
Using \eqref{eqn:ED2} we have
\begin{align}
    \E[D_2 |\mc{H}_i] -\frac{n^{s-2} |C|^{-r-k}}{\binom n2} g_2 &\ge -\frac{n^{s-2} |C|^{-r-k}}{\binom n2} \cdot \frac{p^2 a(t) ^{\binom{s}{2}-\ell-1}f_\Ext}{1-t}.\label{eqn:d2g2}
\end{align}
Using \eqref{eqn:ED3} we have
\begin{align}
    &\E[D_3 |\mc{H}_i] -\frac{n^{s-2}|C|^{-r-k}}{\binom n2} g_3 \nn\\
    &\ge -\frac{n^{s-2}|C|^{-r-k}}{\binom n2} \sbrac{  2p^4 \eps a(t) ^{\binom{s}{2}-\ell-2}(1-t)\log n \cdot f_\Ava + 2p^2 \eps  a(t) ^{\binom{s}{2}-\ell-1} \log n \cdot f_\Ext}\nn\\
    & \qquad \qquad +O\rbrac{n^{\pow(\bS, \be, z) -2 - \m + (20p^4 + 60p^2 +2)\eps}}.\label{eqn:d3g3}
\end{align}
Now to bound line \eqref{eqn:extsupermtg1} we take \eqref{eqn:d1g1} minus \eqref{eqn:d2g2} minus \eqref{eqn:d3g3} minus 
\[
\frac{n^{s-2}|C|^{-r-k}}{\binom n2}a(t)^{\binom{s}{2}-\ell-1} \sbrac{f_\Ext'(t) - p^2 \eps \log n \cdot f_\Ext(t)} + O\rbrac{n^{\pow(\bS, \be, z) - 4 + \eps}}.
\]
We get
\begin{align}
    &\E[\D \Ext^+(z, \phi) | \mc{H}_i]\nn\\
    & \le \frac{n^{s-2}|C|^{-r-k}}{\binom n2} a(t) ^{\binom{s}{2}-\ell-1} \sbrac{\frac{2p^2 + 2p^4\eps(1-t) \log n}{a(t)}   f_\Ava(t) + \rbrac{\frac{2p^2}{1-t} + 3p^2 \eps \log n}  f_\Ext- f_\Ext'(t)  }\nn\\
    & \qquad +O\rbrac{n^{\pow(\bS, \be, z) -2 - \m + (20p^4 + 60p^2 +2)\eps}}.\label{eqn:EDextplus}
\end{align}
 We verify that we have chosen $f_\Ava, f_\Ext$ that make the above line negative. Indeed, since $a(t) = e^{-H(t)}$ we have
 \[
  f_\Ext(t)  = n^{-\hm + 20p^2\eps} (1-t)^{-5p^2+1}a(t)^{-1}e^{10p^4 \eps t \log n} = n^{-\hm + 20p^2\eps} (1-t)^{-5p^2+1}e^{10p^4 \eps t \log n + H(t)}
 \]
and so
\begin{align*}
    f_\Ext'(t) &= n^{-\hm + 20p^2\eps} \Big( (5p^2-1)(1-t)^{-5p^2} e^{10p^4 \eps t \log n +H(t)}\\
    & \qquad \qquad \qquad + (1-t)^{-5p^2+1} (10p^4 \eps \log n + h(t))e^{10p^4 \eps t \log n +H(t)} \Big)\\
  & \ge n^{-\hm + 20p^2\eps} \Big((5p^2-1)(1-t)^{-5p^2} e^{10p^4 \eps t \log n +H(t)} \\
  & \qquad \qquad \qquad + 10p^4 \eps \log n \cdot (1-t)^{-5p^2+1}  e^{10p^4 \eps t \log n +H(t)} \Big).
\end{align*}
Thus we have that the expression in brackets on line \eqref{eqn:EDextplus} is 
\begin{align*}
    &\frac{2p^2 + 2p^4\eps(1-t) \log n}{a(t)}   f_\Ava(t) + \rbrac{\frac{2p^2}{1-t} + 3p^2 \eps \log n}  f_\Ext-  f'_\Ext \\
    &\le n^{-\hm + 20p^2\eps}\Bigg[\Big(2p^2 + 2p^4\eps(1-t) \log n\Big)   (1-t)^{-5p^2} e^{10p^4 \eps t \log n+H(t)} + 2p^2(1-t)^{-5p^2} e^{10p^4 \eps t \log n+H(t)} \\
    & \qquad +3p^2 \eps \log n \cdot (1-t)^{-5p^2+1} e^{10p^4 \eps t \log n+H(t)} \\
    & \qquad - (5p^2-1)(1-t)^{-5p^2} e^{10p^4 \eps t \log n +H(t)} - 10p^4  \eps \log n \cdot (1-t)^{-5p^2+1} e^{10p^4 \eps t \log n +H(t)}\Bigg]\\
    & \le  n^{-\hm + 20p^2\eps}(1-t)^{-5p^2} e^{10p^4 \eps t \log n+H(t)} \Bigg[2p^2 + 2p^4\eps(1-t) \log n     + 2p^2 
    +3p^2 \eps \log n \cdot (1-t)  - (5p^2-1)\\
    &\qquad \qquad -  10p^4 \eps \log n \cdot (1-t) \Bigg]\\
    &=n^{-\hm + 20p^2\eps}(1-t)^{-5p^2} e^{10p^4 \eps t \log n+H(t)} \Bigg[ 1 - p^2 + \rbrac{3p^2  - 8p^4 }\eps (1-t) \log n   \Bigg]\\
    & = -\Omega\rbrac{n^{-\hm + 20p^2\eps}}.
\end{align*}
Thus line \eqref{eqn:EDextplus} is 
\begin{align*}
& - \frac{n^{s-2}|C|^{-r-k}}{\binom n2} a(t) ^{\binom{s}{2}-\ell-1} \Omega\rbrac{n^{-\hm + 20p^2\eps}} +O\rbrac{n^{\pow(\bS, \be, z) -2 - \m + (20p^4 + 60p^2 +2)\eps}}\\
& = -\TOM \rbrac{n^{\pow(\bS, \be, z) - 2 -\hm + 10p^2\eps}} +O\rbrac{n^{\pow(\bS, \be, z) -2 - \m + (20p^4 + 60p^2 +2)\eps}}\\
& < 0. 
\end{align*}

\subsection{Applying Freedman's inequality to \texorpdfstring{$\Ext^+(z, \phi)$}{1}}

Note that by \eqref{eqn:D1bound}, \eqref{eqn:D2bound}, and \eqref{eqn:D3bound} we have
\[
|\D \Ext(z, \phi)| = |D_1 - D_2 - D_3|\le D_1 + D_2 + D_3 \le n^{\pow(\bS,\be,z)-\m + 30p^2 \eps}.
\]

 Thus  in our application of Lemma~\ref{lem:Freedman}, we will use $D=n^{\pow(\bS,\be,z)-\m + 30p^2 \eps}$. Also we have
\begin{align*}
  \Var[ \Delta\Ext^+(z, \phi)| \mc{H}_{k}] = \Var[ \Delta \Ext(z, \phi)| \mc{H}_{k}] &\le \E[ \rbrac{\Delta \Ext(z, \phi)}^2| \mc{H}_{k}]\\
 &\le \E\sbrac{ n^{\pow(\bS,\be,z)-\m + 30p^2 \eps} \cdot |\Delta \Ext(z, \phi)|\big| \mc{H}_{k}} \\
 &\le n^{\pow(\bS,\be,z)-\m + 30p^2 \eps} \cdot \frac{n^{s-2}|C|^{-r-k}}{\binom n2}6p^2 \epsilon \log n\\
 &= \frac{ \TO \rbrac{n^{2\pow(\bS,\be,z)-\m + 30p^2 \eps}}}{\binom n2},
\end{align*}
where the last inequality follows from using \eqref{equation ED1 2}, \eqref{equation D2 2} and \eqref{eqn:ED3part3} to write
\[
\E\sbrac{|\Delta \Ext(z, \phi)|\big| \mc{H}_{k}}\leq 3\max\{\E[D_1\big| \mc{H}_{k}],\E[D_2\big| \mc{H}_{k}],\E[D_3\big| \mc{H}_{k}]\}\leq \frac{n^{s-2}|C|^{-r-k}}{\binom n2}6p^2 \epsilon \log n.
\]
and so for our application of Freedman we have for all $i \le i_{max} < \binom n2$ that 
\[
V(i) = \sum_{0 \le k \le i} \Var[ \D \Ava^+(e,k)| \mc{H}_{k}] = \TO \rbrac{n^{2\pow(\bS,\be,z)-\m + 30p^2 \eps}},
\]
So, we will use $b = \TO \rbrac{n^{2\pow(\bS,\be,z)-\m + 30p^2 \eps}}.$ If $z$ has at least one colored edge then $\Ext(z, \phi, 0) = \ext_z(0) = 0$ and so we have 
 \[
 \Ext^+(z, \phi, 0) = \Ext(z, \phi, 0) - \ext_z(0) - n^{s-2}|C|^{-r-k} f_z(0) = - n^{s-2}|C|^{-r-k} f_z(0) .
 \]
 If $z$ has no colored edge then $r=k=0$ and $\Ext(z, \phi, 0) = (n)_{s-2}$ while $\ext_z(0)=n^{s-2}$ and so 
  \[
 \Ext^+(z, \phi, 0) = \Ext(z, \phi, 0) - \ext_z(0) - n^{s-2} f_z(0) = - n^{s-2} f_z(0) +O(n^{s-3}).
 \]
 Either way we have $\Ext^+(z, \phi, 0) = -(1+o(1)) n^{s-2}|C|^{-r-k} f_z(0)$.
 For $\Ext^+(z, \phi)$ to become positive would therefore be a positive change of $\lambda := (1+o(1))n^{s-2}|C|^{-r-k} f_z(0)=\TT\rbrac{n^{\pow(\bS, \be, z) - \hm + 20p^2 \eps}}$. Freedman's inequality gives us a failure probability of at most
\[
\exp\left(-\frac{\lambda^2}{2(b+D\lambda) }\right) \le \exp\rbrac{-n^{\Omega(\eps)}}
\]
which beats any polynomial union bound.

\section{Phase 2}\label{sec:finishing}

In this section we show that we can finish the coloring using a fresh set of colors and a simple random coloring on the remaining uncolored edges. We assume the event $\mc{E}_{i_{max}}$ holds. Let $C'$ be a set of colors disjoint from $C\cup \bC$ with $|C'|=|C|$. Then we independently color each edge not colored in Phase 1 with a color chosen uniformly at random from $C'$. We now show that this simple coloring scheme leaves us with a $(p,q)$-coloring with positive probability.

We use the asymmetric Lov\'asz local lemma found in \cite{AS}:

\begin{lem}[Lov\'asz local lemma] \label{lem:LLL}
Let $\mc{A}$ be a finite set of events in a probability space $\Omega$ and let $G$ be a dependency graph for $\mc{A}$. Suppose there is an assignment $x:\mc{A} \rightarrow [0, 1)$ of real numbers to $\mc{A}$ such that for all $A \in \mc{A}$ we have
\begin{equation}\label{eqn:LLLcond}
    \Pr(A) \le x(A) \prod_{B \in N(A)} (1-x(B)).
\end{equation}
Then we have that the probability none of the events in $\mc{A}$ happen is 
\begin{equation*}
    \Pr\rbrac{\bigwedge_{A \in \mc{A}} \overline{A}} \ge \prod_{A \in \mc{A}} (1-x(A)) >0.
\end{equation*}
\end{lem}

Let $S$ be a set of $p$ vertices. Say the process has left $S$ with $r \le \binom p2 -q$ repeats and $\ell$ colored edges. If $\ell \ge q+r-1$ then there is no chance that $S$ could end up with too many repeats. Indeed, even if all of the $\binom p2 - \ell \le \binom p2 -q- r + 1$ uncolored edges of $S$ all got the same color, this would be at most $\binom p2 -q- r $ new repeats for a total of at most $\binom p2 -q$. So we assume $\ell \le q+r-2$. To each such set $S$ we associate the number 
\[
x(S) := n^{\frac12 \eps -\f \rbrac{\binom p2 -q- r + 1}} = n^{\frac12 \eps  -p+2+r \cdot \f}
\]

We bound the number of such sets $S$ containing a fixed edge $e$. Each $S$ must fit some legal type using only Platonic colors on its colored edges, and so $S$ fits a trackable type. Since it only uses Platonic colors, the number of choices for this trackable type is $O\rbrac{|C|^{\binom p2 - \ell - 1}}$. For any such type $z$  and a $\phi$ of order 2 we have at step $i=i_{max}$ that (below $t=t_{max}$)
\begin{align*}
  \Ext(z, \phi) &\le \ext_z(t) + n^{p-2}|C|^{-r}f_{z}(t)  \\
  &= n^{p-2}|C|^{-r} a(t) ^{\binom{p}{2}-\ell-1}\sbrac{  t^\ell (1-t)^{\binom{p}{2}-\ell-1}  +  n^{-\hm + 20p^2\eps} (1-t)^{-5p^2+1}a(t)^{-1} e^{10p^4 \eps t \log n}}\\
  &= \TO\rbrac{n^{p-2-r \cdot \f - (\eps + \delta) \rbrac{\binom{p}{2}-\ell-1}}},
\end{align*}
where on the last line we have used that there exists some constant $\delta$ (depending on our other constants) such that $a(t_{max}) \le n^{-\d}$.
Indeed, from \eqref{eqn:hdef}, we can see that $h(t)\geq h_p(t)\geq ct^{\binom p2 - 1}\log n$  for some constant $c>0$, and so
\[
H(t_{max})=\int_0^{t_{max}} h(\tau) d\tau \geq \int_0^{t_{max}} c\tau^{\binom p2 - 1} \log n \;d\tau = \frac{c}{\binom p2} (1-n^{-\epsilon})^{\binom p2}\log n>\frac c{p^2}\log n.
\]
Thus, $a(t_{max})=\exp\cbrac{-H(t_{max}) }\leq n^{-\frac c{p^2}}$, so $\delta=\frac c{p^2}$ suffices.

If we sum the above over all $O\rbrac{|C|^{\binom p2 - \ell - 1}}$ choices for $z$, each set $S$ is counted $(1+o(1))(|C|a(t))^{\binom p2 - \ell - 1} \ge (1+o(1))(|C|n^{-\eps})^{\binom p2 - \ell - 1}$ times. Thus the number of $S$ is 
\[
\TO\rbrac{\frac{n^{p-2-r \cdot \f  - (\eps + \delta) \rbrac{\binom{p}{2}-\ell-1}} \cdot |C|^{\binom p2 - \ell - 1}}{(|C|n^{-\eps})^{\binom p2 - \ell - 1}}} = \TO\rbrac{n^{p-2-r \cdot \f  -   \delta \rbrac{\binom{p}{2}-\ell-1}} } = \TO\rbrac{n^{p-2-r \cdot \f  -   \delta } },
\]
where we have used that $\ell \le q+r-2$ and $r \le \binom p2 - q$. 
We apply Lemma \ref{lem:LLL}, where we have a bad event $A_S$ for each set $S$ of $p$ vertices, $r$ repeats and $\ell \le p+r-2$ colored edges. In particular the bad event $A_S$ is that $S$ gets an additional $\binom p2 - q - r +1$ repeats. We check condition \eqref{eqn:LLLcond}. We have
\[
\Pr(A_S) = O\rbrac{\frac{1}{|C|^{\binom p2 - q - r +1}}} = \TO\rbrac{n^{-\f \rbrac{\binom p2 - q - r +1}}}= \TO\rbrac{n^{-p+2 +r \cdot \f}}
\]
which is smaller than 
\[
x_S \prod_{r} \rbrac{1- n^{\frac12 \eps  -p+2+r \cdot \f}}^{\TO\rbrac{n^{p-2-r \cdot \f-    \delta}}  } = (1-o(1))x_S.
\]

 \section{Some properties of the functions  \texorpdfstring{$ \edgcoi, \coi$}{1} and \texorpdfstring{$\pow$}{1}}\label{sec:platonic}
 
Our goal of this section is to prove the technical lemmas (namely lemmas \ref{lem:preforbpowbound}, \ref{lem:appendpreforb} and~\ref{lem:union}) which were used in previous sections. We first will need to build up some machinery involving the functions $\edgcoi$, $\coi$ and $\pow$.

We have the following additive property. 

\begin{obs}[Additive property]\label{obs:additive} For any type $y$ on $\bT$ and $\bE_3 \subseteq \bE_2 \subseteq \bE_1 \subseteq \binom \bT 2 $ we have
\[
\edgcoi(\bE_1, \bE_3, y) = \edgcoi(\bE_1, \bE_2, y) + \edgcoi(\bE_2, \bE_3, y).
\]
Also for any $\bT_3 \subseteq \bT_2 \subseteq \bT_1 \subseteq \bT$ we have
\[
\coi(\bT_1, \bT_3, y) = \coi(\bT_1, \bT_2, y) + \coi(\bT_2, \bT_3, y)
\]
and
\[
\pow(\bT_1, \bT_3, y) = \pow(\bT_1, \bT_2, y) + \pow(\bT_2, \bT_3, y).
\]
\end{obs}

\begin{obs}[Monotonicity]\label{obs:monotonic}
 For any type $y$ on $\bT$ and $\bE_4 \subseteq \bE_3 \subseteq \bE_2 \subseteq \bE_1 \subseteq \binom \bT 2 $ we have
\[
\edgcoi(\bE_2, \bE_4, y) \le  \edgcoi(\bE_1, \bE_4, y) 
\]
and 
\[
\edgcoi(\bE_2, \bE_4, y) \ge  \edgcoi(\bE_2, \bE_3, y).
\]
In other words, the function $\edgcoi$ is increasing in its first argument and decreasing in its second argument. 
Similarly for any $\bT_4 \subseteq \bT_3 \subseteq \bT_2 \subseteq \bT_1 \subseteq \bT$ we have
\[
\coi(\bT_2, \bT_4, y) \le  \coi(\bT_1, \bT_4, y) 
\]
and 
\[
\coi(\bT_2, \bT_4, y) \ge  \coi(\bT_2, \bT_3, y).
\]
\end{obs}

\begin{defn}\label{def:phi}
For any type $y$ on $\bT$ and $ \bE_3 \subseteq  \bE_2 \subseteq \bE_1 \subseteq \binom{\bT}{2}$  we let 
\[
\nu(\bE_1, \bE_2, \bE_3, y) := \Big|\free(\bE_1 \sm \bE_3, \bE_2 \sm \bE_3, y) \sm \free(\bE_1, \bE_2,  y) \Big|.
\]
In other words, $\nu(\bE_1, \bE_2, \bE_3, y)$ is the number of Platonic colors $c$ such that $c$ appears (on colored edges under the labeling $y$) in $\bE_1 \sm \bE_2$ and in $\bE_3$ but not in $\bE_2 \sm \bE_3$.
\end{defn}

\begin{obs}\label{obs:edgeexcision}
     For any type $y$ on $\bT$ and $ \bE_3 \subseteq  \bE_2 \subseteq \bE_1 \subseteq \binom{\bT}{2}$ we have 
    \begin{align}
         \edgcoi(\bE_1, \bE_2, y)&= \edgcoi(\bE_1 \sm \bE_3, \bE_2 \sm \bE_3, y) + \nu(\bE_1, \bE_2, \bE_3, y) \label{eqn:edgeexcision1}\\
         &\ge \edgcoi(\bE_1 \sm \bE_3, \bE_2 \sm \bE_3, y)\label{eqn:edgeexcision2}.
    \end{align}

\end{obs}
\begin{proof}
 Obviously \eqref{eqn:edgeexcision1} implies \eqref{eqn:edgeexcision2}. To see \eqref{eqn:edgeexcision1}, first recall that  
  \[
      \edgcoi(\bE_1, \bE_2, y) = |M(\bE_1, \bE_2, y)| - |\free(\bE_1, \bE_2, y)|
  \]
  and 
   \[
      \edgcoi(\bE_1 \sm \bE_3, \bE_2 \sm \bE_3, y) = |M(\bE_1 \sm \bE_3, \bE_2 \sm \bE_3, y)| - |\free(\bE_1 \sm \bE_3, \bE_2 \sm \bE_3, y)|.
  \]
   Now observe that $M(\bE_1, \bE_2, y) = M(\bE_1 \sm \bE_3, \bE_2 \sm \bE_3, y)$ and $\free(\bE_1, \bE_2, y) \subseteq \free(\bE_1 \sm \bE_3, \bE_2 \sm \bE_3, y)$, and \eqref{eqn:edgeexcision1} follows. 
\end{proof}

The next observation looks just like line \eqref{eqn:edgeexcision2} except we do not assume $\bE_3 \subseteq \bE_2$.

\begin{obs}\label{obs:sloppyexcision}
For any type $y$ on $\bT$ and edge sets $\bE_2, \bE_3 \subseteq \bE_1 \subseteq \binom{\bT}2$ we have
\[
\edgcoi\rbrac{\bE_1, \bE_2, y } \ge \edgcoi\rbrac{\bE_1 \sm \bE_3, \bE_2 \sm \bE_3, y }
\]
\end{obs}

\begin{proof}
 We have
\begin{align*}
 \edgcoi\rbrac{\bE_1, \bE_2, y } &\ge \edgcoi\rbrac{\bE_1 \sm (\bE_2 \cap \bE_3), \bE_2 \sm (\bE_2 \cap \bE_3), y } \\
 & \ge \edgcoi\rbrac{\bE_1 \sm \bE_3, \bE_2 \sm \bE_3, y }. 
\end{align*}
The first line follows from Observation \ref{obs:edgeexcision}. The second line follows by monotonicity of $\edgcoi$ (from the first line to the second, we decreased the first argument and kept the second argument the same). 

\end{proof}

Given two vertex sets $A$ and $B$, we will denote the collection of edges with one endpoint in $A$ and one in $B$ as $[A,B]$.

\begin{lem}[Excision property] \label{lem:excision}
    For any type $y$ on $\bT$ and $ \bT_3 \subseteq \bT_2 \subseteq \bT_1 \subseteq \bT$ we have
    \begin{equation}
    \coi(\bT_1, \bT_2, y) = \coi(\bT_1 \sm \bT_3 ,\; \bT_2 \sm \bT_3 ,\; y) +   \edgcoi\rbrac{ \bF_1, \bF_2, y} +\nu\rbrac{\bF_2, \bF_3, \bF_4, y}\label{eqn:excisioncoi}
    \end{equation}
    where
    \begin{equation*}
           \bF_1 := \binom{\bT_1}{2},\qquad \bF_2 := \binom{\bT_1 }{2} \sm [\bT_1 \sm \bT_2, \bT_3],\qquad \bF_3 := \binom{\bT_2}2,\qquad \bF_4:=\binom{\bT_3}{2} \cup [\bT_2 \sm \bT_3, \bT_3].
    \end{equation*}

    We also have
    \begin{equation}
    \pow(\bT_1, \bT_2, y)= \pow(\bT_1 \sm \bT_3 ,\; \bT_2 \sm \bT_3 ,\; y) - \f \sbrac{\edgcoi\rbrac{ \bF_1, \bF_2, y} +\nu\rbrac{\bF_2, \bF_3, \bF_4, y}} \label{eqn:excisionpow}.
    \end{equation}
\end{lem}

We will often use Lemma \ref{lem:excision} and the fact that $\edgcoi, \nu$ are nonnegative to get inequalities by dropping the $\edgcoi$ or the $\nu$ term from \eqref{eqn:excisioncoi} or \eqref{eqn:excisionpow}. For example by dropping both terms we get the simpler inequalities
   \begin{equation*}
          \coi(\bT_1 , \bT_2, y) \ge \coi(\bT_1 \sm \bT_3 ,\; \bT_2 \sm \bT_3 ,\; y)
    \end{equation*}
    and
     \begin{equation*}
         \pow(\bT_1 ,\; \bT_2,\; y) \le  \pow(\bT_1 \sm \bT_3 ,\; \bT_2 \sm \bT_3 ,\; y).
     \end{equation*}
     
Recall that $\nu(\bF_2, \bF_3, \bF_4, y)$ is the number of Platonic colors appearing in $\bF_2 \sm \bF_3$ and $\bF_4$ but not in $\bF_3 \sm \bF_4$. Thus it will be helpful for future reference to write
\begin{equation*}
  \bF_2 \sm \bF_3 = \binom{\bT_1 \sm \bT_2}{2} \cup [\bT_1 \sm \bT_2, \bT_2 \sm \bT_3], \qquad \bF_3 \sm \bF_4 = \binom{\bT_2 \sm \bT_3}{2}, \qquad \bF_4=\binom{\bT_3}{2} \cup [\bT_2 \sm \bT_3, \bT_3]
\end{equation*}

\begin{proof}[Proof of Lemma \ref{lem:excision}]
    Since \eqref{eqn:excisioncoi} implies \eqref{eqn:excisionpow} it suffices to prove \eqref{eqn:excisioncoi}. Note that $\bF_4 \subseteq \bF_3 \subseteq \bF_2 \subseteq \bF_1.$
    We have
    \begin{align*}
       \coi(\bT_1 , \bT_2, y) &= \edgcoi\rbrac{\binom{\bT_1}2, \binom{\bT_2}2, y} = \edgcoi\rbrac{ \bF_1, \bF_3, y} \\ 
       & = \edgcoi\rbrac{ \bF_1, \bF_2, y} + \edgcoi\rbrac{ \bF_2, \bF_3, y}\\
      & = \edgcoi\rbrac{ \bF_1, \bF_2, y}+ \edgcoi\rbrac{ \bF_2\sm \bF_4, \bF_3\sm \bF_4, y}+ \nu\rbrac{ \bF_2, \bF_3, \bF_4, y}\\
      & =  \edgcoi\rbrac{\binom{\bT_1 \sm \bT_3}{2} , \binom{\bT_2 \sm \bT_3}2 , y}  +\edgcoi\rbrac{ \bF_1, \bF_2, y}+ \nu\rbrac{ \bF_2, \bF_3, \bF_4, y}\\
      & =\coi(\bT_1 \sm \bT_3, \bT_2 \sm \bT_3, y) +  \edgcoi\rbrac{ \bF_1, \bF_2, y}+ \nu\rbrac{ \bF_2, \bF_3, \bF_4, y}.
    \end{align*}
    Indeed, the first line is by definition of $\coi$, and the second is by additivity. The third line follows by applying Observation \ref{obs:edgeexcision} to $\edgcoi\rbrac{ \bF_2, \bF_3, y}$. On the fourth line we have rearranged terms and used the simple combinatorial identities 
    \[
    \bF_2\sm \bF_4 = \binom{\bT_1 \sm \bT_3}{2}, \qquad \bF_3\sm \bF_4 = \binom{\bT_2 \sm \bT_3}{2}.
    \]
    On the last line we use the definition of $\coi$ again. 
\end{proof}

\begin{obs}\label{obs:colportpow}
For any type $y$, and any $\be\subseteq \bT'\subseteq \bT$,
\[
\pow(\bT,\bT',y)=\pow(\bT,\bT',\colport(y)).
\]
\end{obs}

\begin{proof}
The uncolored edges of $y$ and $\colport(y)$ do not affect the value of $\pow(\bT,\bT',y)$ or $\pow(\bT,\bT',\colport(y))$, and $y$ and $\colport(y)$ are identical everywhere except on uncolored edges, so the result follows.
\end{proof}

\begin{lem}\label{lem:preforbpowbound}
Let $y$ be a $(\be',c',\be'',c'')$-preforbidder on $\bS$. Let $\a=\a(\bS')$ be the number of Platonic colors in $\{c', c''\}$ that appear in $\bS$ but not in $\bS'$.  Then for all $\bS'$ with $\be'\cup\be''\subseteq \bS'\subsetneq \bS$,
\begin{equation}\label{eqn:powdeltabound}
    \pow(\bS,\bS',y)\leq \f \a  - \m
\end{equation}
\end{lem}

\begin{proof}
    Consider the type $y_{\be',\be''}$ on $\bS$ formed from $y$ by replacing the label on $\be'$ with $\col{c'}$ and the label $\be''$ with $\col{c''}$ (and making no other changes). Then by Definition~\ref{def:ecpreforb}~\eqref{def:ecpreforb3}, $y_{\be',\be''}$ is not a legal type. In particular, the colored edges of $\bS$ under the labeling $y_{\be',\be''}$ have $R(|\bS|)+1$ repeats. Say that $x$ is the number of real colors among the colored edges of $\bS$ under $y_{\be',\be''}$. Then using \eqref{eqn:Rbounds} we have
    \[
    \pow(\bS, \emptyset, y_{\be',\be''}) \le |\bS| - 0 - \f \rbrac{\frac{(|\bS|-2)\rbrac{\binom{p}{2}-q+1}}{p-2}  + x} = 2 - \f x.
    \]
    Meanwhile, as a consequence of Definition~\ref{def:ecpreforb} \eqref{def:ecpreforb4}, the colored edges of $\bS'$ under the labeling $y_{\be',\be''}$ have at most $R(|\bS'|)$ repeats. There are at most  $x$ real colors among the colored edges of $\bS'$ under $y_{\be',\be''}$. Thus we have 
     \begin{align*}
    \pow(\bS', \emptyset, y_{\be',\be''}) &\ge |\bS'| - 0 - \f\rbrac{\frac{(|\bS'|-2)\rbrac{\binom{p}{2}-q+1}}{p-2} - \frac{1}{p-2}  + x}\\
    &= 2 + \m - \f x.
    \end{align*}
     Thus we have
    \begin{equation}\label{eqn:y''powbound}
    \pow(\bS, \bS', y_{\be',\be''}) = \pow(\bS, \emptyset, y_{\be',\be''}) - \pow(\bS', \emptyset, y_{\be',\be''}) \le  -\m.
    \end{equation}
    Now to get \eqref{eqn:powdeltabound}, we consider the effect of uncoloring the edges $\be, \be'$. In particular, recall that
    \[
\coi(\bS, \bS', y_{\be',\be''})=\edgcoi\left(\binom{\bS}{2},\binom{\bS'}2,y_{\be',\be''}\right)=\left|\textsc{M}\left(\binom{\bS}{2},\binom{\bS'}2,y_{\be',\be''}\right)\right|-\left|\free\left(\binom{\bS}{2},\binom{\bS'}2,y_{\be',\be''}\right)\right|,
\]
and
\[
\textsc{M}\left(\binom{\bS}{2},\binom{\bS'}2,y_{\be',\be''}\right)=\textsc{M}\left(\binom{\bS}{2},\binom{\bS'}2,y\right),
\]
while
\[
\left|\free\left(\binom{\bS}{2},\binom{\bS'}2,y_{\be',\be''}\right)\right|=\left|\free\left(\binom{\bS}{2},\binom{\bS'}2,y\right)\right|-\a,
\]
so
\[
\pow(\bS,\bS',y_{\be',\be''})=\pow(\bS,\bS',y)-\f \a,
\]
and this combined with \eqref{eqn:y''powbound} implies the result.
\end{proof}

\begin{lem}\label{lem:cutcoi}
    Let $y$ be a $(\be,c,\be',c')$-preforbidder on $\bT$. Let $\bT=\bR_1\cup \bR_2 \cup \bR_3$ be a partition such that $\bR_1,\bR_2\neq \emptyset$ and $|\bR_3| \le 2$ ($\bR_3$ can be empty). Then
    \[
    \edgcoi\left(\binom{\bT}{2},\binom{\bT}{2}\setminus [\bR_1,\bR_2],y\right)\ge 1.
    \]
\end{lem}

\begin{proof}
Let $|\bR_1|=r_1, |\bR_2|=r_2, |\bR_3|=r_3$. Let us assume without loss of generality that $r_1\geq r_2$ and so we have $1 \le r_2 \le p/2$. Note that since $y$ is a preforbidder we have
\begin{equation}\label{eqn:cutcoilower}
\coi(\bT,\be,y) \ge \frac{(|T|-2)\rbrac{\binom{p}{2}-q+1}}{p-2} - 2.
\end{equation}

On the other hand we have
\begin{align}
&  \coi(\bT,\be,y)=\edgcoi \left(\binom{\bT}{2},\emptyset,y\right)\nn\\
&= \edgcoi\left(\binom{\bT}{2},\binom{\bT}{2}\setminus [\bR_1,\bR_2],y\right)+\edgcoi\rbrac{\binom{\bT}{2}\setminus [\bR_1,\bR_2], \binom{\bR_1 \cup \bR_3}2, y } + \edgcoi\left(\binom{\bR_1 \cup \bR_3}2,\emptyset,y\right)\nn\\
& \le \edgcoi\left(\binom{\bT}{2},\binom{\bT}{2}\setminus [\bR_1,\bR_2],y\right)+ \binom {r_2}{2} + r_2r_3 + \frac{(r_1+r_3-2)\rbrac{\binom{p}{2}-q+1}}{p-2}+3.\label{eqn:cutcoiupper}
\end{align}
Indeed, the first line is by definition and the second line uses Observation \ref{obs:additive}. The last line follows since $y$ is legal ($\bR_1 \cup \bR_3$ can contain at most $\frac{(r_1+r_3-2)\rbrac{\binom{p}{2}-q+1}}{p-2}$ repeats and 3 real colors) and since there are only $\binom {r_2}{2} + r_2r_3$ many edges in $\binom{\bT}{2}\setminus [\bR_1,\bR_2]$ that are not in $\binom{\bR_1 \cup \bR_3}2.$ Now putting \eqref{eqn:cutcoilower} and \eqref{eqn:cutcoiupper} together we get
\begin{align*}
 \edgcoi\left(\binom{\bT}{2},\binom{\bT}{2}\setminus [\bR_1,\bR_2],y\right) &\ge \frac{(|T|-2)\rbrac{\binom{p}{2}-q+1}}{p-2}-2 - \binom {r_2}2 - r_2r_3- \frac{(r_1+r_3-2)\rbrac{\binom{p}{2}-q+1}}{p-2}-3\\
 & \ge \frac{r_2\rbrac{\binom{p}{2}-q+1}}{p-2}- \binom {r_2}2 - 2r_2 - 5\\
 & = r_2 \rbrac{\frac{\binom{p}{2}-q+1}{p-2}- \frac{r_2-1}{2} - 2}-5 \\
 & \ge r_2 \rbrac{\frac{\binom{p}{2}-q+1}{p-2}- \frac{p/2-1}{2} - 2}-5 \\
 & = r_2 \rbrac{\frac{p^2 - 6p +16 -4q}{4(p-2)}}-5 \\
 & \ge \frac{p^2 - 6p +16 -4q}{4(p-2)}-5 = \frac{p^2 - 26p +56 -4q}{4(p-2)}>0.\\
\end{align*}
In the second line we used that $|\bT|=r_1 + r_2 + r_3$ and that $r_3 \le 2$. On the fourth line we use $1 \le r_2 \le p/2$, and on the last line we use the same again, along with the bound $q \le \frac{p^2 - 26p +55}{4}$ in the hypothesis of Theorem~\ref{thm:main}. Since $\edgcoi$ only returns integer values we are done.
\end{proof}

\begin{lem}\label{lem:technical}
Let $y_1$ be a trackable type on $\bS_1$ with root $\be_1$, and let $y_2$ be a $(\be_2,c_2,\be_3,c_3)$-preforbidder on $\bS_2$ that is compatible with $y_1$. Let $\bS' \subseteq \bS_1 \cup \bS_2$. Assume the following:
\begin{enumerate}[(i)]
    \item \label{cond:tech1} $\bS_2 \not \subseteq \bS_1$
    \item \label{cond:tech2}$\be_1 \subseteq \bS'$
    \item \label{cond:tech3}$ \be_2 \cup \be_3 \subseteq (\bS' \cup \bS_1) \cap \bS_2$
\end{enumerate}
Then 
 \begin{equation}\label{eqn:techgoal}
    \pow(\bS_1 \cup \bS_2, \bS',  y_1\cup y_2) \le \pow(\bS_1, \be_1, y_1) - \m.
    \end{equation}
\end{lem}

\begin{proof}
    We will consider two cases. 
    
    \tbf{Case 1:} $\bS_2 \sm \bS_1 \not \subseteq \bS'.$ Equivalently, $(\bS' \cup \bS_1) \cap \bS_2 \neq \bS_2.$ In this case we have by additivity that
    \begin{equation}
       \pow(\bS_1 \cup \bS_2, \bS', y_1\cup y_2) = \pow(\bS' \cup \bS_1, \bS', y_1\cup y_2) +\pow(\bS_1 \cup \bS_2, \bS' \cup \bS_1, y_1\cup y_2)\label{eqn:beforeexcision}
\end{equation}
Now using Lemma \ref{lem:excision} 
we have   
\[
\pow(\bS' \cup \bS_1, \bS', y_1\cup y_2) \le \pow( \bS_1, \bS' \cap \bS_1, y_1\cup y_2).
\]
Using Lemma \ref{lem:excision} with $\bT_1 = \bS_1 \cup \bS_2,\; \bT_2 = \bS' \cup \bS_1$ and $\bT_3   = \bS_1 \sm \bS_2$ we have
\begin{align}
    \pow(\bS_1 \cup \bS_2, \bS' \cup \bS_1, y_1\cup y_2)  \le \pow( \bS_2, (\bS' \cup \bS_1) \cap \bS_2, y_2) - \f \nu\rbrac{\bF_2, \bF_3, \bF_4, y_1 \cup y_2},\label{eqn:append1}
\end{align}
where 
\[
           \bF_2 = \binom{\bT_1 }{2} \sm [\bT_1 \sm \bT_2, \bT_3],\qquad \bF_3 = \binom{\bT_2}2,\qquad \bF_4=\binom{\bT_3}{2} \cup [\bT_2 \sm \bT_3, \bT_3].
\]
 Let $\a$ be the number of Platonic colors in $\{c_2, c_3\}$ that appear in $\bS_2$, but which do not appear in $(\bS' \cup \bS_1) \cap \bS_2$. Note that this is precisely the value $\a=\a((\bS' \cup \bS_1) \cap \bS_2)$ we need to bound $\pow( \bS_2, (\bS' \cup \bS_1) \cap \bS_2, y_2)$ using Lemma \ref{lem:preforbpowbound}. We claim that $\nu\rbrac{\bF_2, \bF_3, \bF_4, y_1 \cup y_2}\ge \a$. Indeed, by Definition \ref{def:phi} this $\nu$ value counts the number of Platonic colors appearing in $\bF_2 \sm \bF_3 $ and in $\bF_4$ but not in $\bF_3 \sm \bF_4$. Suppose a Platonic color $c^* \in \{c_2, c_3\}$ is counted by $\a$. Then $c^*$ does not appear in $\binom{\bS_2 \cap (\bS' \cup \bS_1) }{2}=\bF_3 \sm \bF_4$. $c^*$ also appears in $\binom{\bS_2}{2}$, so it must appear in $\binom{\bS_2}{2} \sm \binom{\bS_2 \cap (\bS' \cup \bS_1) }{2} = \bF_2 \sm \bF_3.$
  Furthermore since $y_1$ is a trackable type assigning $\unc{c^*}$ to some edge, we know $c^*$ must appear in $ \binom{\bS_1}{2} \sm \binom{\bS_2 \cap (\bS' \cup \bS_1) }{2} \subseteq \bF_4.$ Thus $c^*$ is also counted by $\nu\rbrac{\bF_2, \bF_3, \bF_4, y_1 \cup y_2}$, and so $\nu\rbrac{\bF_2, \bF_3, \bF_4, y_1 \cup y_2} \ge \a.$
  
Using \eqref{eqn:append1} we now have 
\[
\pow(\bS_1 \cup \bS_2, \bS' \cup \bS_1, y_1\cup y_2)  \le \pow( \bS_2, (\bS' \cup \bS_1) \cap \bS_2, y_1\cup y_2) - \f \a
\]
and so \eqref{eqn:beforeexcision} becomes 
    \begin{align*}
       \pow(\bS_1 \cup \bS_2, \bS', y_1\cup y_2) &\le \pow( \bS_1, \bS' \cap \bS_1, y_1) + \pow( \bS_2, (\bS' \cup \bS_1) \cap \bS_2, y_2)-\f \a\\
        & \le \pow( \bS_1, \bS' \cap \bS_1, y_1) + \f \a  - \m-\f \a\\
        & = \pow( \bS_1, \be_1, y_1) - \pow(  \bS' \cap \bS_1, \be_1, y_1)-\m \\
        & \le \pow( \bS_1, \be_1, y_1) -\m 
\end{align*}
where on the second line we have used Lemma \ref{lem:preforbpowbound} (note that $(\bS' \cup \bS_1) \cap \bS_2\neq \bS_2$ follows from the case assumption) and the last inequality follows from $\pow(  \bS' \cap \bS_1, \be_1, y_1) \ge 0$ (since $y_1$ is trackable). So, \eqref{eqn:techgoal} holds and we are done in this case. 

    \tbf{Case 2:} $\bS_2 \sm \bS_1  \subseteq \bS'.$ In this case we have $\bS' \cup \bS_1 = \bS_1 \cup \bS_2.$ We have
        \begin{align}
       \pow(\bS_1 \cup \bS_2, \bS', y_1\cup y_2) &= \pow(\bS' \cup \bS_1, \bS', y_1\cup y_2) \nn\\
       & \le  \pow( \bS_1, \bS' \cap \bS_1, y_1)  - \f \edgcoi\rbrac{\bF_1, \bF_2 , y_1 \cup y_2}\nn\\
       & = \pow( \bS_1, \be_1, y_1) - \pow(  \bS' \cap \bS_1, \be_1,  y_1)  - \f \edgcoi\rbrac{\bF_1, \bF_2 , y_1 \cup y_2}, \label{eqn:append2}
\end{align}
where on the second line we have applied Lemma \ref{lem:excision} using 
\[
\bT_1 = \bS' \cup \bS_1,\qquad  \bT_2 = \bS', \qquad \bT_3=\bS' \sm \bS_1,
\]
and where
\[
           \bF_1 = \binom{\bT_1}{2},\qquad \bF_2 = \binom{\bT_1 }{2} \sm [\bT_1 \sm \bT_2, \bT_3].
\]
Now if $\bS' \cap \bS_1\neq\be_1$ then $\pow(  \bS' \cap \bS_1, \be_1,  y_1) \ge \m$ by Observation~\ref{obs:trackable}~\eqref{obs:trackable2} since $y_1$ is a trackable type, so \eqref{eqn:append2} would imply \eqref{eqn:techgoal}
 and we would be done.  Otherwise we have $\bS' \cap \bS_1 = \be_1$ and since we are also assuming $\bS_2 \sm \bS_1  \subseteq \bS'$ we have that $\bS' = \be_1 \cup (\bS_2 \sm \bS_1)$. Thus we have 
 \[
\bT_1 = \bS_1 \cup \bS_2,\qquad  \bT_2 = \be_1 \cup (\bS_2 \sm \bS_1), \qquad \bT_3=\bS_2 \sm \bS_1,
\]
and
\[
           \bF_1 = \binom{\bS_1 \cup \bS_2}{2},\qquad \bF_2 = \binom{\bS_1 \cup \bS_2 }{2} \sm [\bS_1 \sm \be_1, \bS_2 \sm \bS_1].
\]
Therefore
\begin{align*}
    \edgcoi\rbrac{\bF_1, \bF_2 , y_1 \cup y_2}& = \edgcoi\rbrac{\binom{\bS_1 \cup \bS_2}{2}, \binom{\bS_1 \cup \bS_2 }{2} \sm [\bS_1 \sm \be_1, \bS_2 \sm \bS_1] , y_1 \cup y_2}\\
    & \ge \edgcoi\rbrac{\binom{\bS_2}{2}, \binom{ \bS_2 }{2} \sm [(\bS_1 \cap \bS_2) \sm \be_1, \bS_2 \sm \bS_1] , y_1 \cup y_2}
\end{align*}
where the second line follows from Observation \ref{obs:sloppyexcision} by removing all edges outside of $\binom{\bS_2}{2}.$
But now by Lemma \ref{lem:cutcoi} with $\bT=\bS_2, \bR_1 = (\bS_1 \cap \bS_2) \sm \be_1, \bR_2 = \bS_2 \sm \bS_1, \bR_3 = \be_1 \cap \bS_2$ we have \newline $\edgcoi\rbrac{\binom{\bS_2}{2}, \binom{ \bS_2 }{2} \sm [(\bS_1 \cap \bS_2) \sm \be_1, \bS_2 \sm \bS_1] , y_1 \cup y_2} \ge 1$, and again \eqref{eqn:append2} implies \eqref{eqn:techgoal}, so we are done. 
\end{proof}

\begin{lem}\label{lem:appendpreforb}
Let $y_1$ be a trackable type on $\bS_1$ with root $\be_1$, and let $y_2$ be a $(\be_2,c_2,\be_3,c_3)$-preforbidder on $\bS_2$ that is compatible with $y_1$. Assume $\be_2,\be_3\subseteq \bS_1$ and that $y_1(\be_2)=\unc{c_2}, y_1(\be_3)=\unc{c_3}$. Assume $\bS_2\setminus\bS_1\neq \emptyset$. Then 
\[
\maxpow(\bS_1\cup\bS_2,\be_1 ,y_1\cup y_2)\leq \pow(\bS_1,\be_1,y_1)-\m.
\]
\end{lem}

\begin{proof}
 Consider $\bS'$ such that $\be_1 \subseteq \bS' \subseteq \bS_1 \cup \bS_2.$ We apply Lemma \ref{lem:technical}. It is easy to see that by our assumptions, Conditions \eqref{cond:tech1} and \eqref{cond:tech2} hold. Condition \eqref{cond:tech3} holds since $\be_2, \be_3$ are contained in both $\bS_1$ and $\bS_2$. Thus by Lemma \ref{lem:technical} we have
 \[
 \pow(\bS_1 \cup \bS_2, \bS',  y_1\cup y_2) \le \pow(\bS_1, \be_1, y_1) - \m.
 \]
 This completes the proof. 
\end{proof}

\begin{lem}\label{lem:union}
Let $y_1$ be a trackable type on $\bS_1$ with root $\be_1$, and let $y_2$ be a $(\be_2,c_2,\be_3,c_3)$-preforbidder on $\bS_2$ that is compatible with $y_1$. Assume $\be_3\subseteq \bS_1$, $\be_3 \neq \be_1$, and $\be_2 \not \subseteq \bS_1$. Assume $y_1(\be_3)=\unc{c_3}$. Assume  $c_2$ is a real color. Then 
    \[
    \maxpow(\bS_1 \cup \bS_2, \be_1 \cup \be_2, y_1 \cup y_2) \le  \pow(\bS_1, \be_1, y_1) - \m.
    \]
\end{lem}

\begin{proof}
 Consider $\bS'$ such that $\be_1 \cup \be_2 \subseteq \bS' \subseteq \bS_1 \cup \bS_2.$ We apply Lemma \ref{lem:technical}.  Condition \eqref{cond:tech1} holds since $\bS_2 \sm \bS_1$ contains $\be_2 \sm \bS_1 \neq \emptyset$. Condition \eqref{cond:tech2} holds by assumption. Condition \eqref{cond:tech3} holds since $\be_2 \subseteq \bS' \cap \bS_2$ and $\be_3 \subseteq \bS_1 \cap \bS_2$. Thus by Lemma \ref{lem:technical} we have
 \[
 \pow(\bS_1 \cup \bS_2, \bS',  y_1\cup y_2) \le \pow(\bS_1, \be_1, y_1) - \m.
 \]
 This completes the proof. 
\end{proof}

\section{The Forbidden Submatching Method}\label{sec:concluding}

Shortly after the authors uploaded the first draft of this paper, Delcourt, Li and Postle \cite{BDLP} (see version 1) extended Theorem \ref{thm:main} to all $q$ below the linear threshold, i.e. $q < \binom p2 - p + 3$, and gave a generalization of the result to list-coloring problems on hypergraphs. Shortly after that (see version 2 of \cite{BDLP}), the first author joined that paper to extend Theorem \ref{thm:main} to all $q$ between the linear and quadratic thresholds, i.e. $\binom p2 - p + 3 < q < \binom{p}2-\lfloor\frac{p}2\rfloor+2$. We summarize the differences and similarities between the two approaches. 

Delcourt and Postle \cite{DPmatch} proved several very general and powerful results that show the existence (under suitable technical conditions) of matchings in certain hypergraphs that do not contain any member of some family of forbidden submatchings. These results and their applications are called the {\bf forbidden submatchings method} (alternatively called {\bf conflict-free matchings} by Glock, Joos, Kim, K\"uhn and Lichev \cite{GJKKL} who independently proved similar results). In particular one of the results takes as input a hypergraph $\mc{H}$, a family of forbidden submatchings, and a subset $A\subseteq V(\mc{H})$, and guarantees that there is a matching that covers all of $A$ and avoids all the forbidden submatchings. Such a result can be used for $(p, q)$-coloring as follows. The hypergraph $\mc{H}$ will just be a bipartite graph with bipartition $A \cup B$, where $A=E(K_n)$ and $B=A \times C$ where $C$ is the set of colors. A vertex $a \in A$ is adjacent to all of the vertices $(a, c)\in B$ (and no other vertices). A matching $M \subseteq \mc{H}$ which covers $A$ corresponds to a coloring of $E(K_n)$ in the obvious way. To ensure that this coloring is a $(p, q)$-coloring, one must forbid our matching $M$ from containing certain submatchings that correspond to a $p$-clique having fewer than $q$ colors. Moreover, to ensure that $\mc{H}$ satisfies the necessary technical conditions to apply the forbidden submatchings method, one must forbid even more submatchings. Indeed, in \cite{BDLP} (as it applies to $(p, q)$-colorings) they forbid any submatching that would correspond to a coloring with some set of $s \le p$ vertices with more than $R(s)$ repeats, i.e. exactly the same number of repeats that we forbid on an $s$-clique in our process in this paper. 

The process we analyzed for phase 1 in this paper is equivalent to choosing one edge of $\mc{H}$ at each step uniformly at random from all edges that can be chosen without intersecting previous edges and without creating a forbidden submatching with the previous edges. The approach in \cite{BDLP} instead applies Theorem 1.16 from \cite{DPmatch} as a black box. The proof of this theorem first uses a random sparsification trick to delete many of the edges of $\mc{H}$, and then uses a nibble method, sometimes called a {\em semi-random} method, to find a matching using only the remaining edges. It is widely accepted that nibble processes are approximately equivalent to random greedy processes, but often the nibble version of a random greedy process turns out to be more amenable to analysis. Adding to the list of tools available for analyzing nibble processes, Delcourt and Postle \cite{DPmatch} proved and utilized a new version of Talagrand's inequality which was crucial in order to cover the sparse regime.




\bibliographystyle{abbrv}
\bibliography{refs}

\begin{thebibliography}{10}

\bibitem{AS}
N.~Alon and J.~H. Spencer.
\newblock {\em The probabilistic method}.
\newblock Wiley Series in Discrete Mathematics and Optimization. John Wiley \&
  Sons, Inc., Hoboken, NJ, fourth edition, 2016.

\bibitem{BEHK22}
J.~Balogh, S.~English, E.~Heath, and R.~A. Krueger.
\newblock Lower bounds on the {E}rdős-{G}yárfás problem via color energy
  graphs.
\newblock To appear in Journal of Graph Theory, 2022.

\bibitem{BCDP2022}
P.~Bennett, R.~Cushman, A.~Dudek, and P.~Prałat.
\newblock {The Erdős-Gyárfás function $f(n, 4, 5) = \frac 56 n + o(n)$ -- so
  Gyárfás was right}.
\newblock {\em arXiv:2207.02920}, 2022.

\bibitem{BDLP}
P.~Bennett, M.~Delcourt, L.~Li, and L.~Postle.
\newblock {On generalized Ramsey numbers in the non-integral regime}.
\newblock {\em arXiv:2212.10542}, 2022.

\bibitem{BD20}
P.~Bennett and A.~Dudek.
\newblock A gentle introduction to the differential equation method and dynamic
  concentration.
\newblock {\em Discrete Math.}, 345(12):Paper No. 113071, 17, 2022.

\bibitem{B2009}
T.~Bohman.
\newblock The triangle-free process.
\newblock {\em Adv. Math.}, 221(5):1653--1677, 2009.

\bibitem{BFL10}
T.~Bohman, A.~Frieze, and E.~Lubetzky.
\newblock A note on the random greedy triangle-packing algorithm.
\newblock {\em J. Comb.}, 1(3-4):477--488, 2010.

\bibitem{BFL15}
T.~Bohman, A.~Frieze, and E.~Lubetzky.
\newblock Random triangle removal.
\newblock {\em Adv. Math.}, 280:379--438, 2015.

\bibitem{BK2010}
T.~Bohman and P.~Keevash.
\newblock The early evolution of the {$H$}-free process.
\newblock {\em Invent. Math.}, 181(2):291--336, 2010.

\bibitem{BK2021}
T.~Bohman and P.~Keevash.
\newblock Dynamic concentration of the triangle-free process.
\newblock {\em Random Structures Algorithms}, 58(2):221--293, 2021.

\bibitem{BES1973}
W.~G. Brown, P.~Erd\H{o}s, and V.~T. S\'{o}s.
\newblock Some extremal problems on {$r$}-graphs.
\newblock In {\em New directions in the theory of graphs ({P}roc. {T}hird {A}nn
  {A}rbor {C}onf., {U}niv. {M}ichigan, {A}nn {A}rbor, {M}ich., 1971)}, pages
  53--63. Academic Press, New York, 1973.

\bibitem{CH18}
A.~Cameron and E.~Heath.
\newblock A {$(5,5)$}-colouring of {$K_n$} with few colours.
\newblock {\em Combin. Probab. Comput.}, 27(6):892--912, 2018.

\bibitem{CH20}
A.~Cameron and E.~Heath.
\newblock New upper bounds for the {E}rd{\H{o}}s-{G}y\'{a}rf\'{a}s problem on
  generalized {R}amsey numbers.
\newblock {\em Combinatorics, Probability and Computing}, page 1–14, 2022.

\bibitem{C2009}
D.~Conlon.
\newblock A new upper bound for diagonal {R}amsey numbers.
\newblock {\em Ann. of Math. (2)}, 170(2):941--960, 2009.

\bibitem{CF2021}
D.~Conlon and A.~Ferber.
\newblock Lower bounds for multicolor {R}amsey numbers.
\newblock {\em Adv. Math.}, 378:Paper No. 107528, 5, 2021.

\bibitem{CFLS15}
D.~Conlon, J.~Fox, C.~Lee, and B.~Sudakov.
\newblock The {E}rd{\H{o}}s-{G}y\'{a}rf\'{a}s problem on generalized {R}amsey
  numbers.
\newblock {\em Proc. Lond. Math. Soc. (3)}, 110(1):1--18, 2015.

\bibitem{CGLS2023}
D.~Conlon, L.~Gishboliner, Y.~Levanzov, and A.~Shapira.
\newblock {A new bound for the {B}rown-{E}rd\H{o}s-{S}\'{o}s problem}.
\newblock {\em J. Combin. Theory Ser. B}, 158(part 2):1--35, 2023.

\bibitem{DPmatch}
M.~Delcourt and L.~Postle.
\newblock {Finding an almost perfect matching in a hypergraph avoiding
  forbidden submatchings}.
\newblock {\em arXiv:2204.089815}, 2022.

\bibitem{DP2022}
M.~Delcourt and L.~Postle.
\newblock {The limit in the $(k+2, k)$-Problem of Brown, Erdős and Sós exists
  for all $k\geq 2$}.
\newblock {\em arXiv:2210.01105}, 2022.

\bibitem{EG97}
P.~Erd\H{o}s and A.~Gy\'{a}rf\'{a}s.
\newblock A variant of the classical {R}amsey problem.
\newblock {\em Combinatorica}, 17(4):459--467, 1997.

\bibitem{E75}
P.~Erd{\H{o}}s.
\newblock Problems and results on finite and infinite graphs.
\newblock In {\em Recent advances in graph theory ({P}roc. {S}econd
  {C}zechoslovak {S}ympos., {P}rague, 1974)}, pages 183--192. (loose errata),
  1975.

\bibitem{ES35}
P.~Erd{\H{o}}s and G.~Szekeres.
\newblock A combinatorial problem in geometry.
\newblock {\em Compositio Math.}, 2:463--470, 1935.

\bibitem{FPS20}
S.~Fish, C.~Pohoata, and A.~Sheffer.
\newblock Local properties via color energy graphs and forbidden
  configurations.
\newblock {\em SIAM J. Discrete Math.}, 34(1):177--187, 2020.

\bibitem{FS2008}
J.~Fox and B.~Sudakov.
\newblock Ramsey-type problem for an almost monochromatic {$K_4$}.
\newblock {\em SIAM J. Discrete Math.}, 23(1):155--162, 2008/09.

\bibitem{G2019}
S.~Glock.
\newblock Triple systems with no three triples spanning at most five points.
\newblock {\em Bull. Lond. Math. Soc.}, 51(2):230--236, 2019.

\bibitem{GJKKL}
S.~Glock, F.~Joos, J.~Kim, M.~K\"{u}hn, and L.~Lichev.
\newblock Conflict-free hypergraph matchings.
\newblock In {\em Proceedings of the 2023 {A}nnual {ACM}-{SIAM} {S}ymposium on
  {D}iscrete {A}lgorithms ({SODA})}, pages 2991--3005. SIAM, Philadelphia, PA,
  2023.

\bibitem{GJKKLP2022}
S.~Glock, F.~Joos, J.~Kim, M.~Kühn, L.~Lichev, and O.~Pikhurko.
\newblock {On the $(6,4)$-problem of Brown, Erdős and Sós}.
\newblock {\em arXiv:2210.01105}, 2022.

\bibitem{GPW20}
H.~Guo, K.~Patton, and L.~Warnke.
\newblock Prague dimension of random graphs.
\newblock To appear in Combinatorica.

\bibitem{KM2008}
A.~Kostochka and D.~Mubayi.
\newblock When is an almost monochromatic {$K_4$} guaranteed?
\newblock {\em Combin. Probab. Comput.}, 17(6):823--830, 2008.

\bibitem{Kurtz1970}
T.~G. Kurtz.
\newblock Solutions of ordinary differential equations as limits of pure jump
  {M}arkov processes.
\newblock {\em J. Appl. Probability}, 7:49--58, 1970.

\bibitem{L87}
H.~Lefmann.
\newblock A note on {R}amsey numbers.
\newblock {\em Studia Sci. Math. Hungar.}, 22(1-4):445--446, 1987.

\bibitem{M04}
D.~Mubayi.
\newblock An explicit construction for a {R}amsey problem.
\newblock {\em Combinatorica}, 24(2):313--324, 2004.

\bibitem{PS19}
C.~Pohoata and A.~Sheffer.
\newblock Local properties in colored graphs, distinct distances, and
  difference sets.
\newblock {\em Combinatorica}, 39(3):705--714, 2019.

\bibitem{RS1976}
I.~Z. Ruzsa and E.~Szemer\'{e}di.
\newblock Triple systems with no six points carrying three triangles.
\newblock In {\em Combinatorics ({P}roc. {F}ifth {H}ungarian {C}olloq.,
  {K}eszthely, 1976), {V}ol. {II}}, volume~18 of {\em Colloq. Math. Soc.
  J\'{a}nos Bolyai}, pages 939--945. North-Holland, Amsterdam-New York, 1978.

\bibitem{Sah2022}
A.~Sah.
\newblock Diagonal ramsey via effective quasirandomness.
\newblock To appear in Duke Mathematical Journal, 2022.

\bibitem{S2022}
W.~Sawin.
\newblock An improved lower bound for multicolor {R}amsey numbers and a problem
  of {E}rd{\H{o}}s.
\newblock {\em J. Combin. Theory Ser. A}, 188:Paper No. 105579, 11, 2022.

\bibitem{BES1973no2}
V.~T. S\'{o}s, P.~Erd\H{o}s, and W.~G. Brown.
\newblock On the existence of triangulated spheres in {$3$}-graphs, and related
  problems.
\newblock {\em Period. Math. Hungar.}, 3(3-4):221--228, 1973.

\bibitem{Warnke2020}
L.~Warnke.
\newblock On {W}ormald's differential equation method.
\newblock To appear in Combin. Probab. Comput..

\bibitem{W2021}
Y.~Wigderson.
\newblock An improved lower bound on multicolor {R}amsey numbers.
\newblock {\em Proc. Amer. Math. Soc.}, 149(6):2371--2374, 2021.

\bibitem{W1995}
N.~C. Wormald.
\newblock Differential equations for random processes and random graphs.
\newblock {\em Ann. Appl. Probab.}, 5(4):1217--1235, 1995.

\bibitem{W1999}
N.~C. Wormald.
\newblock The differential equation method for random graph processes and
  greedy algorithms.
\newblock In {\em {Lectures on Approximation and Randomized Algorithms, PWN,
  Warsaw}}, pages 73--155, 1999.

\end{thebibliography}

\end{document}